\documentclass{siamart190516}
\usepackage{amsfonts}
\usepackage{graphicx}
\ifpdf
  \DeclareGraphicsExtensions{.eps,.pdf,.png,.jpg}
\else
  \DeclareGraphicsExtensions{.eps}
\fi
\usepackage[mathscr]{eucal}
\usepackage{amssymb}
\usepackage{enumerate}
\usepackage{cases}
\usepackage{pgf,pgfplots,tikz}
\usepackage{color}
\usepackage{csquotes}
\usepackage{stmaryrd}

\usepgfplotslibrary{groupplots}
\usepackage[show]{ed}
\definecolor{dhcol}{rgb}{0,0.5,0}

\definecolor{sccol}{rgb}{0,0,0.5}

\definecolor{amcol}{rgb}{0.5,0,0}

\newcommand{\dsp}{\displaystyle}

\newcommand{\Arg}{\mathrm{Arg}}

\numberwithin{theorem}{section}

\newcommand{\TheTitle}{The complex-scaled half-space Matching Method}
\newcommand{\TheAuthors}{A.-S. Bonnet-Ben Dhia et al.}%, Simon N. Chandler-Wilde, Sonia Fliss, Christophe Hazard, Karl-Mikael Perfekt, Yohanes Tjandrawidjaja}

\headers{\TheTitle}{\TheAuthors}

\title{{\TheTitle}\thanks{Submitted to the editors DATE.
\funding{{The research of K.-M.~Perfekt was partially supported by UK Engineering and Physical Sciences Research Council (EPSRC) grant EP/S029486/1. The work of Y.~Tjandrawidjaja was supported by the French Commissariat \`a l'Energie Atomique et aux Energies Alternatives (CEA), through funding for PhD studies.}}}}

\author{Anne-Sophie Bonnet-Ben Dhia\thanks{POEMS (CNRS-INRIA-ENSTA Paris), Institut Polytechnique de Paris, Palaiseau, France (\email{anne-sophie.bonnet-bendhia@ensta-paris.fr}, \email{sonia.fliss@ensta-paris.fr}, \email{christophe.hazard@ensta-paris.fr})} \and %828 Boulevard des Mar\'echaux, 91762 Palaiseau Cedex, France} \and
  Simon N. Chandler-Wilde\thanks{Department of Mathematics and Statistics, University of Reading, Whiteknights PO Box 220, Reading RG6 6AX, UK
    (\email{s.n.chandler-wilde@reading.ac.uk}, \email{k.perfekt@reading.ac.uk})}
    \and Sonia Fliss\footnotemark[2]
    \and
    Christophe Hazard\footnotemark[2]
    \and
    Karl-Mikael Perfekt\thanks{Department of Mathematical Sciences, Norwegian University of Science and Technology (NTNU), NO-7491 Trondheim, Norway \email{karl-mikael.perfekt@ntnu.no})}
    \and
    Yohanes Tjandrawidjaja\thanks{Lehrstuhl I, Fakult\"at f\"ur Mathematik, Technische Universit\"at Dortmund, Vogelpothsweg 87,
44227 Dortmund, Germany (\email{Yohanes.Tjandrawidjaja@math.tu-dortmund.de})}
}
\ifpdf
\hypersetup{
  pdftitle={\TheTitle},
  pdfauthor={\TheAuthors}
}
\fi
\graphicspath{{Figs/}}
\begin{document}

\newcommand{\rf}[1]{(\ref{#1})}
\newcommand{\mmbox}[1]{\fbox{\ensuremath{\displaystyle{ #1 }}}}	% my box to highlight main formulas
%-------------------------------------------------------------------------------------------------------------------------------
\newcommand{\hs}[1]{\hspace{#1mm}}
\newcommand{\vs}[1]{\vspace{#1mm}}
%-------------------------------------------------------------------------------------------------------------------------------
\newcommand{\ri}{{\mathrm{i}}}
\newcommand{\re}{{\mathrm{e}}}
\newcommand{\rd}{\mathrm{d}}
%-------------------------------------------------------------------------------------------------------------------------------
\newcommand{\R}{\mathbb{R}}
\newcommand{\Q}{\mathbb{Q}}
\newcommand{\N}{\mathbb{N}}
\newcommand{\Z}{\mathbb{Z}}
\newcommand{\C}{\mathbb{C}}
\newcommand{\K}{{\mathbb{K}}}
%-------------------------------------------------------------------------------------------------------------------------------
\newcommand{\cA}{\mathcal{A}}
\newcommand{\cB}{\mathcal{B}}
\newcommand{\cC}{\mathcal{C}}
\newcommand{\cS}{\mathcal{S}}
\newcommand{\cD}{\mathcal{D}}
\newcommand{\cH}{\mathcal{H}}
\newcommand{\cI}{\mathcal{I}}
\newcommand{\cItilde}{\tilde{\mathcal{I}}}
\newcommand{\cIhat}{\hat{\mathcal{I}}}
\newcommand{\cIcheck}{\check{\mathcal{I}}}
\newcommand{\cIstar}{{\mathcal{I}^*}}
\newcommand{\cJ}{\mathcal{J}}
\newcommand{\cM}{\mathcal{M}}
\newcommand{\cP}{\mathcal{P}}
\newcommand{\cV}{{\mathcal V}}
\newcommand{\cW}{{\mathcal W}}
\newcommand{\scrD}{\mathscr{D}}
\newcommand{\scrS}{\mathscr{S}}
\newcommand{\scrJ}{\mathscr{J}}
\newcommand{\sD}{\mathsf{D}}
\newcommand{\sN}{\mathsf{N}}
\newcommand{\sS}{\mathsf{S}}
 \newcommand{\sT}{\mathsf{T}}
 \newcommand{\sH}{\mathsf{H}}
 \newcommand{\sI}{\mathsf{I}}
%-------------------------------------------------------------------------------------------------------------------------------
\newcommand{\bs}[1]{\mathbf{#1}}
\newcommand{\bb}{\mathbf{b}}
\newcommand{\bd}{\mathbf{d}}
\newcommand{\bn}{\mathbf{n}}
\newcommand{\bp}{\mathbf{p}}
\newcommand{\bP}{\mathbf{P}}
\newcommand{\bv}{\mathbf{v}}

\newcommand{\bxi}{\boldsymbol{\xi}}
\newcommand{\boldeta}{\boldsymbol{\eta}}	%Added by Dave
\newcommand{\ts}{\tilde{s}}
\newcommand{\tGamma}{{\tilde{\Gamma}}}
 \newcommand{\tbx}{\tilde{\bx}}
 \newcommand{\tbd}{\tilde{\bd}}
 \newcommand{\txi}{\xi}
%-------------------------------------------------------------------------------------------------------------------------------
\newcommand{\done}[2]{\dfrac{d {#1}}{d {#2}}}
\newcommand{\donet}[2]{\frac{d {#1}}{d {#2}}}
\newcommand{\pdone}[2]{\dfrac{\partial {#1}}{\partial {#2}}}
\newcommand{\pdonet}[2]{\frac{\partial {#1}}{\partial {#2}}}
\newcommand{\pdonetext}[2]{\partial {#1}/\partial {#2}}
\newcommand{\pdtwo}[2]{\dfrac{\partial^2 {#1}}{\partial {#2}^2}}
\newcommand{\pdtwot}[2]{\frac{\partial^2 {#1}}{\partial {#2}^2}}
\newcommand{\pdtwomix}[3]{\dfrac{\partial^2 {#1}}{\partial {#2}\partial {#3}}}
\newcommand{\pdtwomixt}[3]{\frac{\partial^2 {#1}}{\partial {#2}\partial {#3}}}
\newcommand{\bnabla}{\boldsymbol{\nabla}}
\newcommand{\dive}{\boldsymbol{\nabla}\cdot}
\newcommand{\curl}{\boldsymbol{\nabla}\times}
\newcommand{\Phixy}{\Phi(\bx,\by)}
\newcommand{\PhiOxy}{\Phi_0(\bx,\by)}
\newcommand{\dxPhixy}{\pdone{\Phi}{n(\bx)}(\bx,\by)}
\newcommand{\dyPhixy}{\pdone{\Phi}{n(\by)}(\bx,\by)}
\newcommand{\dxPhiOxy}{\pdone{\Phi_0}{n(\bx)}(\bx,\by)}
\newcommand{\dyPhiOxy}{\pdone{\Phi_0}{n(\by)}(\bx,\by)}
%-------------------------------------------------------------------------------------------------------------------------------
\newcommand{\eps}{\varepsilon}
%-------------------------------------------------------------------------------------------------------------------------------
\newcommand{\real}[1]{{\rm Re}\left[#1\right]} %Changed by Steve from \re to avoid conflict with \re below
\newcommand{\im}[1]{{\rm Im}\left[#1\right]}
\newcommand{\ol}[1]{\overline{#1}}
\newcommand{\ord}[1]{\mathcal{O}\left(#1\right)}
\newcommand{\oord}[1]{o\left(#1\right)}
\newcommand{\Ord}[1]{\Theta\left(#1\right)}
%-------------------------------------------------------------------------------------------------------------------------------
\newcommand{\hsnorm}[1]{||#1||_{H^{s}(\bs{R})}}
\newcommand{\hnorm}[1]{||#1||_{\tilde{H}^{-1/2}((0,1))}}
\newcommand{\norm}[2]{\left\|#1\right\|_{#2}}
\newcommand{\normt}[2]{\|#1\|_{#2}}
\newcommand{\on}[1]{\Vert{#1} \Vert_{1}}
\newcommand{\tn}[1]{\Vert{#1} \Vert_{2}}
%\newcommand{\norm}[2]{\left\|#1\right\|_{#2}}
%-------------------------------------------------------------------------------------------------------------------------------
\newcommand{\xt}{\mathbf{x},t}
\newcommand{\PhiF}{\Phi_{\rm freq}}
\newcommand{\cone}{{c_{j}^\pm}}
\newcommand{\ctwo}{{c_{2,j}^\pm}}
\newcommand{\cthree}{{c_{3,j}^\pm}}
%-------------------------------------------------------------------------------------------------------------------------------
\newtheorem{thm}{Theorem}[section]
\newtheorem{lem}[theorem]{Lemma}
\newtheorem{defn}[theorem]{Definition}
\newtheorem{prop}[theorem]{Proposition}
\newtheorem{cor}[theorem]{Corollary}
\newtheorem{rem}[theorem]{Remark}
\newtheorem{conj}[theorem]{Conjecture}
\newtheorem{ass}[theorem]{Assumption}
\newtheorem{example}[theorem]{Example} %Added by Steve
%-------------------------------------------------------------------------------------------------------------------------------
\newcommand{\tH}{\widetilde{H}}
\newcommand{\Hze}{H_{\rm ze}} 	% Zero extension space
\newcommand{\uze}{u_{\rm ze}}		% Zero extension of a function u
\newcommand{\dimH}{{\rm dim_H}}
\newcommand{\dimB}{{\rm dim_B}}
\newcommand{\IntClosOm}{\mathrm{int}(\overline{\Omega})}
\newcommand{\IntClosOmOne}{\mathrm{int}(\overline{\Omega_1})}
\newcommand{\IntClosOmTwo}{\mathrm{int}(\overline{\Omega_2})}
\newcommand{\Ccomp}{C^{\rm comp}}
\newcommand{\tCcomp}{\tilde{C}^{\rm comp}}
\newcommand{\uC}{\underline{C}}
\newcommand{\utC}{\underline{\tilde{C}}}
\newcommand{\oC}{\overline{C}}
\newcommand{\otC}{\overline{\tilde{C}}}
\newcommand{\capcomp}{{\rm cap}^{\rm comp}}
\newcommand{\Capcomp}{{\rm Cap}^{\rm comp}}
\newcommand{\tcapcomp}{\widetilde{{\rm cap}}^{\rm comp}}
\newcommand{\tCapcomp}{\widetilde{{\rm Cap}}^{\rm comp}}
\newcommand{\hcapcomp}{\widehat{{\rm cap}}^{\rm comp}}
\newcommand{\hCapcomp}{\widehat{{\rm Cap}}^{\rm comp}}
\newcommand{\tcap}{\widetilde{{\rm cap}}}
\newcommand{\tCap}{\widetilde{{\rm Cap}}}
\newcommand{\ccap}{{\rm cap}}
\newcommand{\ucap}{\underline{\rm cap}}
\newcommand{\uCap}{\underline{\rm Cap}}
\newcommand{\cCap}{{\rm Cap}}
\newcommand{\ocap}{\overline{\rm cap}}
\newcommand{\oCap}{\overline{\rm Cap}}
\DeclareRobustCommand
{\mathringbig}[1]{\accentset{\smash{\raisebox{-0.1ex}{$\scriptstyle\circ$}}}{#1}\rule{0pt}{2.3ex}}
\newcommand{\cirH}{\mathringbig{H}}
\newcommand{\cirHs}{\mathringbig{H}{}^s}
\newcommand{\cirHt}{\mathringbig{H}{}^t}
\newcommand{\cirHm}{\mathringbig{H}{}^m}
\newcommand{\cirHzero}{\mathringbig{H}{}^0}
\newcommand{\deO}{{\partial\Omega}}
\newcommand{\OO}{{(\Omega)}}
\newcommand{\Rn}{{(\R^n)}}
\newcommand{\Id}{{\mathrm{Id}}}
\newcommand{\gap}{\mathrm{Gap}}
\newcommand{\ggap}{\mathrm{gap}}
\newcommand{\isom}{{\xrightarrow{\sim}}}
\newcommand{\half}{{1/2}}
\newcommand{\mhalf}{{-1/2}}
\newcommand{\inter}{{\mathrm{int}}}

\newcommand{\Hsp}{H^{s,p}}
\newcommand{\Htq}{H^{t,q}}
\newcommand{\tHsp}{{{\widetilde H}^{s,p}}}
\newcommand{\SP}{\ensuremath{(s,p)}}
\newcommand{\Xsp}{X^{s,p}}

\newcommand{\dd}{{d}}\newcommand{\pp}{{p_*}}

\newcommand{\Rnn}{\R^{n_1+n_2}}
\newcommand{\Tr}{{\mathrm{Tr}}}
\newcommand{\sO}{\mathsf{O}}
\newcommand{\sC}{\mathsf{C}}
\newcommand{\sA}{\mathsf{A}}
\newcommand{\sM}{\mathsf{M}}
\newcommand{\sF}{\mathsf{F}}
\newcommand{\sG}{\mathsf{G}}
\newcommand{\mS}{\Gamma}
\newcommand{\omS}{{\overline{\mS}}}
\newcommand{\sumpm}[1]{\{\!\!\{#1\}\!\!\}}
\newcommand{\bsx}{\boldsymbol{x}}
\newcommand{\bsy}{\boldsymbol{y}}
\newcommand{\bsz}{\boldsymbol{z}}
\newcommand{\bsw}{\boldsymbol{w}}
\newcommand{\bx}{\bsx}
\newcommand{\by}{\bsy}
\newcommand{\bz}{\bsz}
\maketitle

% REQUIRED
\begin{abstract}
The Half-Space Matching (HSM) method has recently been developed as a new method for the solution of 2D scattering problems with complex backgrounds, providing an alternative to {Perfectly Matched Layers (PML)} or other artificial boundary conditions. Based on  half-plane representations for the solution, the scattering problem is rewritten as a system of integral equations in which the unknowns are restrictions of the solution to the boundaries of a finite number of  overlapping half-planes contained in the domain: this integral equation system is coupled to a standard finite element discretisation localised around the scatterer. While satisfactory numerical results have been obtained for real wavenumbers, wellposedness and equivalence to the original scattering problem have been established only for complex wavenumbers. In the present paper, by combining the HSM framework with a complex-scaling technique, we provide a new formulation for real wavenumbers  which is provably well-posed and has the attraction for computation that the complex-scaled solutions of the integral equation system decay exponentially at infinity. The analysis requires the study of double-layer  potential integral operators on intersecting infinite lines, and their analytic continuations.  The effectiveness of the method is validated by preliminary numerical results.	
\end{abstract}

% REQUIRED

% REQUIRED

\section{Introduction and the scattering problem}
\label{sec-introduction}
The mathematical and numerical analysis of scattering by bounded obstacles and/or inhomogeneities  in a homogeneous background is a mature research area, and there are many effective numerical schemes, at least for low to moderately high frequencies. However, when the background is heterogeneous (stratified, periodic,...) and/or anisotropic, especially when electromagnetic or elastic waves are considered, {many} theoretical questions are still open and the design of efficient numerical methods {remains} a significant challenge.

It is well known that, for homogeneous backgrounds, the Sommerfeld radiation condition ensures well-posedness of the scattering problem \cite{CoKr:98}. The extension of this standard  radiation condition to the aforementioned more complex backgrounds is really intricate (see, e.g., \cite{Bon-Dak-Haz-Cho-2009, Bon-Gou-Haz-2011, Cut-Haz-1998, Fli-Jol-2016, Jer-Ned-2012}). Moreover, {a} Green's function or tensor, which could be used to derive an integral {equation} formulation of the problem, is in general not available or hard to compute. Finally, Perfectly Matched Layer (PML) techniques, which are popular in homogeneous backgrounds because they are easy to implement, {can} produce spurious effects {for complex backgrounds,} as is well-known for anisotropic backgrounds in relation to instabilities in the time domain \cite{Bec-Fau-Jol-2003}.

To overcome these difficulties a new method, called the Half-Space Matching (HSM) method, has been developed {as an (exact) artificial boundary condition} for two-dimensional scattering problems. This method is based on explicit or semi-explicit expressions for the outgoing solutions of radiation problems in half-planes, these expressions established by using Fourier, generalized Fourier, or Floquet transforms when the background is, respectively, homogeneous \cite{BbFT18, Bon-Fli-Tja-2019} (and possibly anisotropic \cite{tonnoir2015,Bar-Bon-Fli-Tja-2018,tjandrawidjaja2019}), stratified {\cite{Ott}}, or periodic \cite{Fli-Jol-2009}. The domain exterior to a bounded region enclosing the scatterers is covered by a finite number of half-planes (at least three). The unknowns of the formulation are the traces of the solution on the boundaries of these half-planes and the restriction of the solution to the bounded region. The system of equations which couples these unknowns is derived by writing compatibility conditions between the different representations of the solution. This {coupled} system {includes} second-kind integral equations on the infinite boundaries of the half-planes.

{This new formulation is attractive and versatile as a method to truncate computational domains in problems of scattering by localised inhomogeneities in complex backgrounds (including backgrounds that may be different at infinity in different directions). It has been employed successfully in numerical implementations for a range of problems, namely } periodic \cite{Fliss:2010}  and stratified media {(including cases with different stratifications in different parts of the background domain)} \cite{Ott}, and anisotropic scalar {and} elastic problems \cite{tonnoir2015,tjandrawidjaja2019}.
%\textcolor{purple}{But}

Up to now the theoretical and numerical analysis of the method has remained an open question {in the challenging, and more practically relevant, non-dissipative case when waves radiate out to infinity.} But a rather complete analysis has been carried out in the {simpler} dissipative case, when the solution (and its traces) decay exponentially at infinity. In that case the analysis can be done using an $L^2$ framework for the traces, and the associated formulation has been shown to be of Fredholm type and well-posed {in a number of cases} \cite{BbFT18,Bon-Fli-Tja-2019}, {with the sesquilinear form of the weak formulation coercive plus compact, enabling standard numerical analysis arguments \cite{Bon-Fli-Tja-2019}.} One of the main difficulties in the %\textcolor{purple}{(more interesting and practically relevant)}
non-dissipative case is the slow decay at infinity of the solution which results in non-$L^2$ traces. The possibility, to address this, of working in the framework initially introduced in \cite{Bon-Til-2000, Bon-Til-2001} was investigated by the authors, but it seems to be inappropriate for the present analysis.

The objective of this paper is to propose a new formulation of the HSM method which is well-suited for theoretical and numerical analysis (and practical computation)  in the non-dissipative case. For the sake of clarity and as a first step, we restrict ourselves in the present paper to a {relatively simple model} problem for which the justification of the method is based on the simple form of the associated Green's function. (Let us mention that the extension of this formulation to anisotropic backgrounds has been already validated and will be the subject of a forthcoming paper.  See Section \ref{sec:Perspectives} for a brief discussion of extensions to other more complex configurations.) This new formulation exploits a fundamental property of the solution in the spirit of the ideas behind complex-scaling methods (e.g. the pioneering works of \cite{balslev:1971},\cite{aguilar:1971}): in any given direction, the solution, as a function of the associated curvilinear abscissa, has an analytic continuation into the complex plane which is exponentially decaying in the upper complex plane. This enables us to replace the system of equations for the traces by similar equations for exponentially-decaying analytical {continuations} of these traces. This recovers {well-posedness in} an $L^2$ framework, {with coercive plus compact sesquilinear forms;  moreover, attractive for computation, the rate of exponential decay of these analytically-continued traces increases with increasing wavenumber.}
Let us mention that in \cite{Lu-Lu-Qian-2018} a similar {integral-equation-based} complex-scaling idea{, namely a boundary integral equation formulation of PML, is used to compute 2D scattering by localised perturbations in a straight interface between different media.}

In the present paper we consider the {rather} simple model case of a scalar equation, the isotropic Helmholtz equation
\begin{equation}\label{pb:probleme_diff_initial}
			-\Delta u -  k^2\,\rho u  = f \; \text{in} \; \Omega,\\[4pt]
\end{equation}
deduced from the wave equation assuming a time-dependence $\re^{-\ri\omega t}$ for a given angular frequency $\omega>0$.
Here $\rho$ is a function in $L^\infty(\Omega)$, bounded from below by a positive constant, and such that $\rho-1$ is compactly supported, and the constant $k=\omega/c$ is the wavenumber, where $c$ is the wave speed outside the support of $\rho-1$, so that \eqref{pb:probleme_diff_initial} models propagation in a domain with a local perturbation in wave speed.

The propagation domain $\Omega$ is $\R^2$, or $\R^2$ minus a set of obstacles which are included in a bounded region. We assume that, for a positive constant $a$,
\[
\quad\partial\Omega\subset \overline{\Omega}_a\mbox{ where } \Omega_a:=  (-a,a)^2.
\]
In the presence of obstacles, boundary conditions have to be added to the model.
 The source term $f$ is supposed to be a function in $L^2(\Omega)$ with  compact support included in $\Omega_a$,  and we assume that
\begin{equation*}%\label{eq:ac-milieu-def-rho}
 	\text{supp}(\rho-1)\subset \Omega_a.
 \end{equation*}
As already mentioned, in order to get well-posedness, one has to prescribe in addition the Sommerfeld radiation condition, that, {for $\bsx:=(x_1,x_2)\in \R^2$,} %, which can be written in polar coordinates $(r,\theta)$ as
 \begin{equation}\label{eq:Sommerfeld}
	 \frac{\partial u}{\partial r}{(\bsx)}- {\rm i} k u{(\bsx)} = o\left(r^{-1/2}\right)\quad\text{as}\;r:=|\bsx|\rightarrow +\infty,
	 \end{equation}
	 uniformly with respect to $\widehat \bsx := \bsx/r$.

In the sequel, we will consider two configurations in order to focus first on the Half-Space Matching formulation and then on its coupling with a variational formulation in a bounded region. In Sections \ref{section-HSMcomplexfreq} to \ref{sec:unique} we consider a Dirichlet problem set in the exterior of the square $\Omega_a$. {The application of the analysis in Sections \ref{section-HSMcomplexfreq}-\ref{sec:unique}  to general configurations}, with source terms, inhomogeneities, {and/or}  obstacles contained in $\Omega_a$, is the object of Section \ref{sec-CHSMgeneralcase}. %Numerical results are shown in Section \ref{sec:numerical} for both the model  Dirichlet problem and a  more general configuration.%\scednote{Am not quite sure what this last sentence means.}
% \begin{enumerate}
% 	\item The propagation medium is  $\Omega=\R^2\setminus\Omega_a$ and non homogeneneous Dirichlet boundary conditions are imposed on its boundary. We assume $f=0$ and $\rho=1$ in $\Omega$.\\
% 	\item We consider then the general case in presence of source terms and possible arbitrary obstacles and perturbations of the coefficient $\rho$. The idea here is that any perturbation can be taken into account using a classical variational formulation, as soon as it is contained in a bounded region.\scednote{Edit this: we don't have obstacles so delete reference to obstacles. And is ``general problem'' the right terminology?}
% \end{enumerate}
% For the sake of brevity, we will refer to the first case as the exterior Dirichlet problem outside a square and to the second case as the general problem.

The outline of the paper is as follows.  In Section \ref{section-HSMcomplexfreq} we recall the main results concerning the HSM formulation in the dissipative case (that is with a complex wavenumber $k$). In previous papers the HSM formulation has been derived using Fourier representations for the solution in half-planes. Here we introduce a new formulation using double-layer potential integral representations %This is the formulation that we will use, for sake of simplicity,
 to derive the so-called complex-scaled version of the method, valid for a real wavenumber $k$.

 The derivation and the analysis of this new formulation is the object of Section \ref{sec-CHSM}{, which finishes with a statement of our main well-posedness result for the new method.} In Section \ref{sec:recon} we establish properties of the solution that can be reconstructed a posteriori, from knowledge of the complex-scaled traces{, notably elements of the far-field pattern.}
These properties are used in Section \ref{sec:unique} to establish the uniqueness result for the complex-scaled HSM problem.

The HSM method approach is extended to {general configurations} in Section \ref{sec-CHSMgeneralcase}{; the analysis in this section depends throughout on the well-posedness and other results obtained in Sections \ref{sec-CHSM}-\ref{sec:unique}.} {In Section \ref{sec:numerical} the implementation of a finite element discretization of the complex-scaled HSM formulation is described} and numerical results are presented, for both the model  Dirichlet problem and more general configurations. {The paper finishes with three appendices to which we defer certain technical details of the analysis.}
% In order to focus on the Half-Space Matching formulation, the main part of the paper will be devoted to the exterior Dirichlet problem. We recall the method for the dissipative case in Section\ based on the Fourier transform, and introduce a new version of the method based on the Green's function. This is this formulation that we will use, for sake of simplicity, to derive the so called complex-scaled version of the method in Section and analyse it.

\section{The HSM method for {complex wavenumber}}
\label{section-HSMcomplexfreq}
In this section, as preparation for studying the HSM method for real wavenumber, we first recall what is known about the method in the dissipative case. It is enough for this purpose to consider the Dirichlet problem for complex wavenumber ($\Im(k)>0$, $\Re(k)>0$) in the exterior of the square $\Omega_a$, i.e.
\begin{equation}\label{pb:probleme_Dir_diss}
	\begin{cases}
			-\Delta u -  k^2\,u  = 0 \; \text{in} \; \Omega:=\R^2\setminus\overline{\Omega_a},\\[4pt]
			u=g\;\text{on}\;\Sigma_a:=\partial\Omega_a,
			\end{cases}
\end{equation}
for a given $g\in H^{1/2}(\Sigma_a)$. It is well-known that Problem \eqref{pb:probleme_Dir_diss} has a unique solution $u\in H^1(\Omega)$.

Let us first {recall} the main results of  \cite{BbFT18}.

The domain $\Omega$ is the union of $4$ overlapping half-planes $\Omega^j$ that {abut} the $4$ edges of the square $\Omega_a$. We introduce the following local coordinates for all $j \in \llbracket0,3\rrbracket:=\{0,1,2,3\}$:
\begin{equation}\label{def:coord-loc-HP-ac-milieu}
	 \begin{pmatrix} x^j_1\\ x^j_2 \end{pmatrix} : =  \begin{pmatrix} \text{cos}(j\pi/2) & \text{sin}(j\pi/2)\\ {-}\text{sin}(j\pi/2) & \;\;\text{cos}(j\pi/2) \end{pmatrix}\, \begin{pmatrix} x_1\\ x_2 \end{pmatrix}.
\end{equation}
The half-planes and their boundaries are defined as follows for all $j \in \llbracket0,3\rrbracket$:
\begin{equation}\label{eq:notations1_1}
	\Omega^j := \{(x^j_1,x^j_2):  x^j_1 > a, x_2^j \in \R\}, \quad \Sigma^j:=\{(x^j_1,x^j_2):  x^j_1 = a, x_2^j \in \R\}.
\end{equation}
Finally, we denote
\begin{equation}\label{eq:notations1_2}
	\Sigma_{a}^j := \Sigma_{a}\cap\Sigma^j.
\end{equation}
These notations are summarised in Figure \ref{Fig:notations_Case1}.
%\begin{figure}[H]
%	\begin{center}
%  	\includegraphics[width=7cm]{Img/Notations1}
%  	\caption{The notations defined in (\ref{def:coord-loc-HP-ac-milieu}-\ref{eq:notations1_1}-\ref{eq:notations1_2}).}\label{Fig:notations_Case1}
%	\end{center}
%\end{figure}
\begin{figure}[!h]
	\definecolor{ffffff}{rgb}{0.2,0.2,0.2}
	\definecolor{yqyqyq}{rgb}{0.5,0.5,0.5}
	\begin{center}
		\begin{tikzpicture}[line cap=round,line join=round,x=0.4cm,y=0.4cm,scale=1,yshift=2]
		\begin{scope}[shift={(-15,0)}]
		\draw [dashed,->] (6,6)--(7.5,6);
		\draw[below] (7.5,6) node {$x_1$};
		\draw[left] (6,7.2) node {$x_2$};
		\draw [dashed,->] (6,6)--(6,7.5);
		\draw (5.5,5.5) node {$\Omega_a$};
		\draw[above] (14,8) node {$\Sigma^1$};
			\draw[right] (8,14) node {$\Sigma^0$};
				\draw[left] (4,-2) node {$\Sigma^2$};
					\draw[below] (-2,4) node {$\Sigma^3$};
		\draw[right] (8,6) node {$\Sigma^0_a$};
		\draw[above] (6,8) node {$\Sigma^1_a$};
		\draw[left] (4,6) node {$\Sigma^2_a$};
		\draw[below] (6,4) node {$\Sigma^3_a$};
		\draw  (14,6) node {$\Omega^0$};
		\draw  (6,14) node {$\Omega^1$};
		\draw  (-2,6) node {$\Omega^2$};
		\draw  (6,-2) node {$\Omega^3$};
		\draw [line width=1.pt] (-2.,4)--(14.,4);
		\draw [line width=1.pt] (-2.,8)--(14.,8);
		\draw [line width=1.pt] (4,-2)--(4,14.);
		\draw [line width=1.pt] (8,-2)--(8,14.);
		\end{scope}
		\end{tikzpicture}
	\end{center}
	\caption{The notations defined in (\ref{def:coord-loc-HP-ac-milieu}-\ref{eq:notations1_1}-\ref{eq:notations1_2}).}\label{Fig:notations_Case1}
\end{figure}
As explained in the introduction, the formulation uses the representation of the solution in each half-plane $\Omega^j$ in terms of its trace on $\Sigma^j$. More precisely, let us denote
\begin{equation}
	\label{eq:def_traces}
	 \varphi^j:=u\big|_{\Sigma^j}\quad \text{for}\;j \in \llbracket0,3\rrbracket
\end{equation}
so that
\begin{equation}
	\label{eq:def_traces2}
	 u\big|_{\Omega^j}=U^j(\varphi^j)\quad \text{for}\;j \in \llbracket0,3\rrbracket
\end{equation}
where, for any $\psi\in H^{1/2}(\Sigma^j)$,  $U^j(\psi)\in H^1(\Omega^j)$ is the unique solution of
\begin{equation}\label{eq:half-space}
	\begin{array}{|lcr}
			-\Delta U^j -  k^2\, U^j = 0 & \text{in} & \Omega^j,\\[4pt]
			\;U^j=\psi&\text{on}&\Sigma^j.
		\end{array}
\end{equation}
{In the sequel, we identify any function defined on $\Sigma^j$, in particular the function $\varphi^j$, with a function of the real variable $x_2^j$.}

We can express $U^j(\psi)$ explicitly in terms of its trace $\psi$ in two manners: using the Fourier transform or using a Green's function representation.
First, using the Fourier transform in the $x_2^j-$direction, it is easy to see that the solution of \eqref{eq:half-space} is given by \\
\begin{equation}\label{eq:expr-Pj}
	 U^j (\psi) ({\bsx}^j)= \dfrac{1}{\sqrt{2\pi}}\int_{\R} \widehat{\psi}(\xi) \re^{- \sqrt{ \xi^2- k^2}(x_1^j-a)} \re^{\ri \xi x_2^j} \rd\xi,\quad {\bsx}^j:=(x^j_1,x_2^j) \in \Omega^j,
\end{equation}
where the square root is defined with the convention $\Re(\sqrt{z})\geq 0$, {for} $z\in \mathbb{C}\setminus\mathbb{R}^-$ (with $\mathbb{R}^-:=(-\infty,0]$) and $\widehat{\psi}$ is the Fourier transformation of $\psi$ using the convention
\begin{equation} \label{eq:FTdef}
	\widehat{\psi}(\xi):= \dfrac{1}{\sqrt{2\pi}}\int_{\R} \psi(x_2^j)\,\re^{-\ri \xi x_2^j} \rd x_2^j,\quad \xi\in\R.
\end{equation}
Secondly, using a Green's function representation, we can show that
\begin{equation} \label{eq:hprGreen}
 U^j (\psi) ({\bsx}^j) = \int_{\Sigma^j}\frac{\partial G^j({\bsx}^j,{\bsy}^j)}{\partial n({\bsy}^j)}\,\psi({\bsy}^j)\, \rd s({\bsy}^j),\quad {\bsx}^j \in \Omega^j,
\end{equation}
where $G^j({\bsx}^j,{\bsy}^j)$ is the Dirichlet Green's function for $\Omega^j$ and $n({\bsy}^j)$ is the unit normal to $\Sigma^j$ that points into $\Omega^j$. Explicitly, $G^j({\bsx}^j,{\bsy}^j)= \Phi({\bsx}^j,{\bsy}^j)-\Phi(\tilde{\bsx}^j,{\bsy}^j)$, with $\tilde{\bsx}^j$ the image of ${\bsx}^j$ in $\Sigma^j$, where
$\Phi({\bsx},{\bsy})$ is the standard fundamental solution of the Helmholtz equation defined by
\begin{equation} \label{eq:Greenfct}
\Phi({\bsx},{\bsy}) := \frac{\ri}{4}H_0^{(1)}( k|\bsx-\bsy|), \quad {\bsx,\,\bsy\in \R^2, \;\;\bsx\neq \bsy},
\end{equation}
so that, equivalently,
\begin{equation} \label{eq:hprPhi}
 U^j (\psi) ({\bsx}^j) =2 \int_{\Sigma^j}\frac{\partial \Phi({\bsx}^j,{\bsy}^j)}{\partial n({\bsy}^j)}\,\psi({\bsy}^j)\, \rd s({\bsy}^j),\quad {\bsx}^j \in \Omega^j.
\end{equation}
 This leads to
 \begin{equation} \label{eq:hpr}
 U^j (\psi) ({\bsx}^j) = \int_{\R}h(x_1^j-a,x_2^j-y_2^j)\,\psi(y_2^j)\, \rd y_2^j,\quad {\bsx}^j \in \Omega^j,
 \end{equation}
 where
 \begin{equation} \label{eq:hpr_kernel}
 h(x_1,x_2):=\frac{\ri kx_1}{2}\frac{H^{(1)}_1( k R(x_1,x_2))}{R(x_1,x_2)},
 \end{equation}
 and
 \begin{equation} \label{eq:R(s,t)}
 R(x_1,x_2):=(x_1^2+x_2^2)^{1/2},\quad x_1,x_2\in\R.
 \end{equation}
Let us remark that the two representations \eqref{eq:expr-Pj} and \eqref{eq:hpr} of $U^j (\psi)$ can be derived the one from the other by using simply a Plancherel formula (e.g.\ \cite[p.~821]{CW-Impedance-1997}).

To derive the system of equations whose unknowns are the traces $\varphi^j$ of the solution, it suffices to write that the half-plane representations must coincide where they coexist.
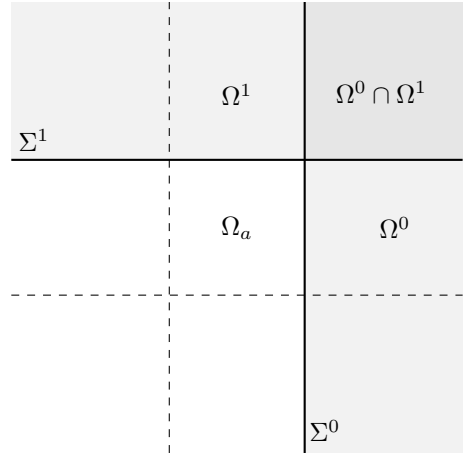
\begin{figure}[!h]
\begin{center}
	\begin{tikzpicture}[scale=.6]
\fill[color=gray!10] (1.5,-5)rectangle(5,5);
\fill[color=gray!10] (-5,1.5)rectangle(5,5);
\fill[color=gray!20] (1.5,1.5)rectangle(5,5);
\draw (0,0) node {$\Omega_a$};
\draw(-4.5,1.5) node[above] {${\Sigma^1}$};
\draw (1.4,-4.5) node[right] {${\Sigma^0}$};
\draw (3.5,0) node  {${\Omega^0}$};
\draw (0,2.5) node[above]  {${\Omega^1}$};
\draw (3.2,3) node  {${\Omega^0}\cap{\Omega^1}$};
\draw[thick] (1.5,-5)--(1.5,5);
\draw[thick] (-5,1.5)--(5,1.5);
\draw[dashed] (-5,-1.5) -- (5,-1.5) ;
\draw[dashed] (-1.5,-5) -- (-1.5,5) ;
%\draw[thick,blue] (1.5,1.5)--(5,1.5);
%\draw[thick,blue] (1.5,1.5)--(1.5,5);
\end{tikzpicture}
\end{center}
\caption{Construction of the compatibility condition.}
\end{figure}
For instance, in the quarter plane $\Omega^0\cap\Omega^1$ we have
\begin{equation}\label{eq:comp_quarter}
	u=U^0(\varphi^0)=U^1(\varphi^1)\quad \text{in}\;\Omega^0\cap\Omega^1,
\end{equation}
and in particular
\begin{equation}\label{eq:comp1}
	\varphi^1=U^0(\varphi^0)\quad \text{on}\;\Omega^0\cap\Sigma^1,
\end{equation}
which leads to a first integral equation linking $\varphi^0$ and $\varphi^1$. Indeed, for any point of $\Sigma^1\cap\Omega^0$, represented by $(x_1^0,x_2^0)$ and $(x_1^1,x_2^1)$ in local coordinates systems, we have $x_2^1=-x_1^0<-a$ and $x_1^1=x_2^0=a$ {(indeed, the point is on $\Sigma^1\cap \Omega^0$ if and only if these equations hold)}. {Thus} the compatibility relation
\eqref{eq:comp1} can be rewritten (identifying $\Sigma^1$ with $\R$ in the way noted above) as
\begin{equation} \label{eq:comp1a}
	\varphi^1(x_2^1)=U^0(\varphi^0)(-x_2^1,a),\quad x_2^1<-a,
\end{equation}
where we can use either of the two integral representations \eqref{eq:expr-Pj} and \eqref{eq:hpr} for the half-plane solution $U^0(\varphi^0)$.
From \eqref{eq:comp_quarter} we have also
\[
	\varphi^0=U^1(\varphi^1)\quad \text{on}\;\Sigma^0\cap\Omega^1,
\]
which leads to another integral equation linking $\varphi^0$ and $\varphi^1$:
\begin{equation*} %\label{eq:comp2a}
	\varphi^0(x_2^0)=U^1(\varphi^1)(x_2^0,-a),\quad x_2^0>a.
\end{equation*}

{Repeating this for each quarter plane} we get 8 equations linking the 4 traces.
% \[
% 	\varphi^1(x_2^1)=\dfrac{1}{\sqrt{2\pi}}\int_{\R} \widehat{\varphi}^0(\xi) \re^{\ri \sqrt{ k^2 - \xi^2}(x_1^0-a)} \re^{\ri \xi a} d\xi
% \]
In order to write the system of equations in a condensed manner, we shall use the same notation $t$ instead of the various variables $x_2^j$ so that the 8 equations become
\begin{equation}\label{eq:syst_comp}
	\forall j\in \llbracket0,3\rrbracket,\quad\begin{array}{|lcl}
		\varphi^j(t)=S\,D\,\varphi^{j-1}(t),\quad t<-a,\\
		\varphi^j(t)=D\,S\,\varphi^{j+1}(t),\quad t>a,
		\end{array}
\end{equation}
where we have set $\varphi^{-1}:=\varphi^3$ and $\varphi^{4}:=\varphi^0$, and where the operators $S$ and $D$ are defined as follows.  For any $\psi\in L^2(\R)$,
\begin{equation}
	\label{eq:sym}
	S\psi(t):=\psi(-t),\quad t\in\mathbb{R},
\end{equation}
and the integral operator $D$ is defined by
\begin{equation}
	\label{eq:opD_comp_def}
	D\psi (t) := U^0(\psi)(t,a),\quad t>a.
\end{equation}
$D$ can be given explicitly by either of the following two expressions:
\begin{equation}
	\label{eq:opD_comp_exprF}
	D\psi (t)=\dfrac{1}{\sqrt{2\pi}}\int_{\R} \widehat{\psi}(\xi) \re^{\ri \sqrt{ k^2 - \xi^2}(t-a)} \re^{\ri \xi a}\, \rd \xi,\quad t>a,
\end{equation}
or
\begin{equation}
	\label{eq:opD_comp_exprG}
	D\psi (t)=\int_{\R}h(t-a,s-a)\,\psi(s)\, \rd s,\quad t>a,
\end{equation}
where the kernel $h$ is defined in \eqref{eq:hpr_kernel}\footnote{Note that, as is clear from (\ref{eq:hprPhi}), $D$ is precisely a double-layer potential operator (in the sense, e.g., of \cite{CoKr:83} or \cite{Ch:84}) from $\{(a,t): t\in\mathbb{R}\}$ to $\{(t,a): t\geq a\}$.  }. The system of equations has to be completed with the Dirichlet boundary condition rewritten as
\begin{equation}
	\label{eq:syst_comp_BC}
	\varphi^j(t)=g|_{\Sigma_a^j}(t),\quad -a<t<a,\quad j \in \llbracket0,3\rrbracket.
\end{equation}

One can easily check that \eqref{eq:syst_comp}-\eqref{eq:syst_comp_BC} is equivalent to the original problem \eqref{pb:probleme_Dir_diss}. More precisely, if $\{\varphi^0,\varphi^1,\varphi^2,\varphi^3\}\in (H^{1/2}(\R))^4$ is a solution to \eqref{eq:syst_comp}-\eqref{eq:syst_comp_BC} then one can recover the solution $u$ to \eqref{pb:probleme_Dir_diss}, from the knowledge of the $\varphi^j$'s, thanks to the half-plane representations  \eqref{eq:expr-Pj} or \eqref{eq:hpr}. Indeed, by uniqueness of Dirichlet quarter-plane problems, two half-plane representations $U^j(\varphi^j)$ and  $U^{j+1}(\varphi^{j+1})$ coincide on the quarter-plane $\Omega^j\cap\Omega^{j+1}$ since the compatibility conditions \eqref{eq:syst_comp} imply that they coincide on its boundary.

For the analysis and the computation it is convenient to consider the formulation in an $L^2$-framework:
\begin{equation}\label{eq:HSMM_complexe}
	\begin{array}{c}
		\text{Find}\; \{\varphi^0,\varphi^1,\varphi^2,\varphi^3\}\in (L^2(\R))^4\;\text{such that}\\[5pt]
	\begin{array}{|lcl}
		\varphi^j(t)=S\,D\,\varphi^{j-1}(t),\quad t<-a,\\[3pt]
		\varphi^j(t)=g|_{\Sigma_a^j}(t),\quad-a<t<a,\\[3pt]
		\varphi^j(t)=D\,S\,\varphi^{j+1}(t),\quad t>a.
		\end{array}\quad j\in \llbracket0,3\rrbracket
		\end{array}
\end{equation}
One attraction of this $L^2$-framework is that it allows the use of elementary operations of restriction and extension. More precisely, for any function of $L^2(\R)$, its restriction to an open interval $I\subset \R$ is in $L^2(I)$. More {significantly}, any function of $L^2(I)$ extended by 0 belongs to $L^2(\R)$. For simplicity, any function defined on a part of $\R$ is identified hereafter with its extension by 0.  With this convention we can write
\begin{equation}\label{eq:decomp_L2}
	L^2(\R)=L^2(-\infty,-a)\oplus L^2(-a,a)\oplus L^2(a,+\infty),
\end{equation}
which will be extensively used hereafter. {In line with this convention}, we define $D\psi(t)$ for all $t\in\mathbb{R}$ by setting
\begin{equation} \label{eq:DdefExt}
D\psi(t):=0,\quad t\leq a.
\end{equation}

With these various conventions (\ref{eq:HSMM_complexe}) can be rewritten shortly as
\begin{equation}
\label{eq:HSMM_uneligne}
\varphi^j=SD	\varphi^{j-1}+DS	\varphi^{j+1}+g|_{\Sigma_a^j},\quad  j\in \llbracket0,3\rrbracket.
\end{equation}	
Also, with the above conventions, results proved in \cite{Ch:84,BbFT18} can  be stated as follows:
\begin{proposition}\label{prop:D_complexe}
\noindent
	\begin{itemize}
		\item[(i)]$D$ is a continuous operator {on} $L^2(\R)$, with range in $L^2(a,+\infty)\subset L^2(\R)$.
		\item[(ii)]As an operator on $ L^2(a,+\infty)$, ${D}$ is the sum of an operator of norm $\leq 1/\sqrt{2}$ and a compact operator.
		\item[(iii)]As an operator from $ L^2(-\infty,-a)$ to $L^2(a,+\infty)$, ${D}$ is compact.
	\end{itemize}
\end{proposition}
\begin{proof}
	{Note first that \emph{(i)} {follows immediately from} \emph{(ii)} and \eqref{eq:DdefExt}, since \emph{(ii)}, together with a symmetry argument with respect to $a$, implies that $D \colon L^2(-\infty,a) \to L^2(a,+\infty)$ is also the sum of an operator of norm $\leq 1/\sqrt{2}$ and a compact operator.}
	Consider the expression \eqref{eq:opD_comp_exprG} for $D$. Because of the dissipation (${\Im}(k)>0$), the kernel $h$ is exponentially decaying at infinity (i.e. as $t$ or $s$ tends to $\infty$). Further, the mapping $(t,s)\mapsto h(t-a,s-a)$ is continuous except at $t=s=a$. Thus \emph{(iii)} is clear since the kernel of $D$ is Hilbert-Schmidt, i.e.
	\[
	(t,s)\mapsto h(t-a,s-a)\;\in L^2((-\infty,-a)\times(a,+\infty)).
	\]
	
	{To} show \emph{(ii)}, the only difficulty comes from the singularity of the kernel $h$ at $t=s=a$. As in Proposition \ref{prop:DLPquad} and Remark \ref{rem:D0bounded} {in Appendix \ref{sec:ComplIntOper}}, let $h_0$ and $D_0$ denote $h$ and $D$, respectively, when $k=0$.
For $b>a$ let %$P_{(a,b)}\psi = \chi_{(a,b)}\psi$ be the product of $\psi$ with
$\chi_{(a,b)}$ denote the characteristic function of $(a,b)$.
Then it is straightforward to see that, for every $b>a$, $D-\chi_{(a,b)}D_0$ is an integral operator with  kernel $h(t-a,s-a)-\chi_{(a,b)}({t})h_0(t-a,s-a)$ that is Hilbert-Schmidt, so $D-\chi_{(a,b)}D_0$ is compact. Further (Remark \ref{rem:D0bounded}), as an operator on $L^2(a,+\infty)$, $\|D_0\|=1/\sqrt{2}$, so also $\|\chi_{(a,b)}D_0\|\leq \|\chi_{(a,b)}\|\, \|D_0\| =1/\sqrt{2}$.
\end{proof}

 Now the system \eqref{eq:HSMM_complexe} can be formulated in an operator form. Let us introduce
\begin{equation} \label{L20def}
	L^2_0(\R):=\{\psi\in L^2(\R):\;\psi(t)=0\;\text{for}\;-a<t<a\}
\end{equation}
and $$\Phi:=\{\varphi^0,\varphi^1,\varphi^2,\varphi^3\}\in (L^2(\R))^4,$$
\begin{equation} \label{eq:Phig}
	\Phi_g:=\{g|_{\Sigma_a^0},g|_{\Sigma_a^1},g|_{\Sigma_a^2},g|_{\Sigma_a^3}\}\in (L^2(-a,a))^4 \subset (L^2(\R))^4.
\end{equation}
Then, noting (\ref{eq:HSMM_uneligne}), the system \eqref{eq:HSMM_complexe} can be rewritten as:
\begin{equation}\label{eq:HSMM_complexe_matrice}
	\begin{array}{c}
		\text{Find}\; \Phi\in (L^2(\R))^4\;\text{such that}\; \Phi-\Phi_g\in  (L^2_0(\R))^4\;\text{and}\\[5pt]
	(\mathbb{I}-\mathbb{D})(\Phi-\Phi_g)=\mathbb{D}\,\Phi_g,
		\end{array}
\end{equation}
where
\begin{equation}\label{eq:def_Dmat}
	\mathbb{D}:=\left[\begin{matrix}
	 0&\,{D}\,S&0&\,S\,{D}\\
	 \,S\,{D}&0&\,{D}\,S&0\\
	 0&S\,{D}&0&{D}\,S\\
	 {D}\,S&0&\,S\,{D}&0
	 \end{matrix}\right].
\end{equation}
In \cite{BbFT18}, the following result is proven:

\newpage

\begin{theorem}\label{th:HSMM_complexe}
\noindent
	\begin{itemize}
		\item[(i)]$\mathbb{D}$ is a continuous operator on $(L^2(\R))^4$ and  $\mathbb{D}((L^2(\R))^4)\subset (L^2_0(\R))^4$.
		\item[(ii)]As an operator on $ (L^2_0(\R))^4$, $\mathbb{D}$ is the sum of an operator of norm $\leq 1/\sqrt{2}$ and a compact operator.
		\item[(iii)]
Problem \eqref{eq:HSMM_complexe_matrice} is well-posed.
	\end{itemize}
\end{theorem}
Let us give some ideas of the proof which will be relevant for the following sections. {The} property \emph{(iii)} is {largely} a consequence of \emph{(ii)} since \emph{(ii)} gives that $(\mathbb{I}-\mathbb{D})$, as an operator acting on $(L^2_0(\R))^4$, is the sum of a coercive operator\footnote{Recall that, {given a Hilbert space $\mathcal{H}$ with inner product $(\cdot,\cdot)$,} we call a bounded linear operator {$A$} on $\mathcal{H}$ {\em coercive} if the corresponding sesquilinear form $a(\cdot,\cdot)$, defined by $a(\phi,\psi)=(A\phi,\psi)$, $\forall \phi,\psi\in \mathcal{H}$, is coercive, i.e., if, for some constant $\gamma>0$, $\Re(a(\phi,\phi))\geq \gamma \|\phi\|^2$, $\forall \phi\in \mathcal{H}$.} and a compact one. By Fredholm theory, it suffices then to show uniqueness (which is not straightforward in the $L^2$ framework, see \cite{BbFT18} for more details). {The} properties \emph{(i)} and \emph{(ii)} are consequences (see Appendix \ref{Appendix_Matrices}) of Proposition \ref{prop:D_complexe}.

\section{The complex-scaled HSM method for real wavenumber}
\label{sec-CHSM}
{
In this section we consider the Dirichlet problem of the previous section, but now with real {wavenumber} ($k>0$). Where $H^1_{\mathrm{loc}}(\Omega):=\{v|_\Omega:v\in H_{\mathrm{loc}}^1(\R^2)\}$, we seek $u\in H^1_{\mathrm{loc}}(\Omega)$ such that
\begin{equation}\label{pb:probleme_Dir_real}
\begin{cases}
-\Delta u -  k^2\,u  = 0 \; \text{in} \; \Omega:=\R^2\setminus\overline{\Omega_a},\\[4pt]
u=g\;\text{on}\;\Sigma_a,
\end{cases}
\end{equation}
for a given $g\in H^{1/2}(\Sigma_a)$, and such that the radiation condition \eqref{eq:Sommerfeld} holds. It is well known that this problem has a unique solution.
% \footnote{{This problem} is well-posed even for $g\in L^2(\partial \Omega_a)$ if the boundary condition is understood in terms of maximal functions and non-tangential convergence (instead of the usual trace sense) \cite{ToWe93}. {And in fact one} can extend Theorem ?? below to show equivalence of the boundary value problem with our {complex-scaled HSM method} formulation \eqref{eq:HSMM_matrice_cmplx} {also for $g\in L^2(\partial \Omega_a)$.}}

The half-plane representations \eqref{eq:def_traces2}, with $U^j(\varphi^j)$ given by \eqref{eq:hprGreen} (equivalently \eqref{eq:hprPhi} or \eqref{eq:hpr}), still hold for $k$ real, and can be derived via Green's theorem using the radiation condition \eqref{eq:Sommerfeld} (cf.\ \cite[Theorem 2.1]{CWZ:98}).
	As a consequence, the traces $\varphi_j$, $j\in  \llbracket0,3\rrbracket$, of the solution $u$ on $\Sigma^j$, still satisfy the system of equations \eqref{eq:HSMM_complexe} when $k>0$.
%	
%	
%	 However, as the solution of (\ref{pb:probleme_Dir_real}-\ref{eq:Sommerfeld}) decays only slowly, like $r^{-1/2}$ as $r\to \infty$, we cannot expect that $\varphi^j\in L^2(\R)$ only that $\varphi^j\in L^2_{\mathrm{loc}}(\R)$.
%	
	 We note that, although the solution of (\ref{pb:probleme_Dir_real})-(\ref{eq:Sommerfeld}) decays only slowly, like $r^{-1/2}$ as $r\to +\infty$, the integrals \eqref{eq:hpr} still make sense. The bound \eqref{eq:Hbound1}, that follows from asymptotics of the Hankel function $H_1^{(1)}$,  implies that, for some constant $C>0$ depending only on $k$,
		\begin{equation*}% \label{eq:Phibound1}
		\left|h(x_1,x_2)\right| \leq C x_1(R^{-3/2} + R^{-2}), \quad x_1>0, \; x_2\in \R,
		\end{equation*}
		where $R:=(x_1^2+x^2_2)^{1/2}$, so that \eqref{eq:hpr} is well-defined, even when $k$ is real, for every $\bsx^j\in \Omega^j$ and $\psi\in L^2_{\mathrm{loc}}(\R)$ with $\psi(t)=\mathcal{O}(1)$ as $|t|\to+\infty$. Moreover, though we shall not need this, it is still possible to rewrite \eqref{eq:hpr} equivalently as \eqref{eq:expr-Pj}, provided that care is taken in interpreting the right hand side of \eqref{eq:expr-Pj}; see the discussion in \cite{ArensHohage05,CWElschner10} and \cite{Bon-Til-2000,Bon-Til-2001}.
	
	From a numerical point of view, the HSM method for real $k$ works well \cite{BbFT18}. However, from a theoretical point of view, the formulation does not make sense in an $L^2$ setting. Indeed, as the solution of (\ref{pb:probleme_Dir_real})-(\ref{eq:Sommerfeld}) decays only  like $r^{-1/2}$ as $r\to +\infty$,   we cannot expect that its traces on $\Sigma^j$, $j\in \llbracket0,3\rrbracket$, are in $L^2(\R)$. In parallel work \cite{bibid} we have shown that the HSM formulation for real $k$ is equivalent with the original problem \eqref{pb:probleme_Dir_real} if we supplement it with radiation conditions analogous to the Sommerfeld condition \eqref{eq:Sommerfeld} (see \eqref{eq:RadCond} below). However, there are still significant gaps in our understanding of this formulation when $k$ is real. In particular, while \eqref{eq:HSMM_complexe} can be written formally in operator form as \eqref{eq:HSMM_complexe_matrice}, just as in the dissipative case, in the case when $k$ is real {there is no obvious function space setting (e.g., $L^p$, or a weighted $L^p$ space)} for which this formulation makes sense, with $\mathbb{D}$ a well-defined bounded linear operator. Consequently, we are not able to justify the numerical method and neither provide a priori error estimates.

%As a consequence, the properties and well-posedness of the formulation are obscure, in particular its equivalence with the original problem \eqref{pb:probleme_Dir_real} (the analogue of Theorem \ref{th:HSMM_complexe}) is unclear.

 These difficulties with the standard formulation for real $k$ are part of the motivation for  the method proposed in this paper that we term the {\em complex-scaled HSM method}.   The idea behind this method is to ``complexify''. Since the pioneering works of Aguilar, Balslev and Combes \cite{aguilar:1971,balslev:1971}, complex-scaling methods have been used intensively to construct the analytic continuation of resolvents in mathematical physics, see for instance \cite{DyZw:19} and the references therein. The complex-scaling method in \cite{aguilar:1971,balslev:1971} is closely related to the idea behind PML (see
for instance \cite{Col-Mon-1998}). Our use of complex-scaling is somewhat different since we consider only analytic continuation of traces of the solution on particular infinite half-lines and apply that complex scaling in an integral equation context. This is similar to manipulations made to
understand analyticity of boundary traces in high frequency scattering problems in
\cite[\S4.1]{NonConvex}. Precisely our plan is as follows:
\begin{enumerate}
	\item From properties of the solution $u$ of  (\ref{pb:probleme_Dir_real})-(\ref{eq:Sommerfeld}) %(\ref{pb:probleme_diff_initial}-\ref{eq:Sommerfeld}),
we deduce that the traces $\varphi^j$, $j\in \llbracket0,3\rrbracket$, have analytic continuations into the complex plane from $(-\infty,-a)$ and from $(a,+\infty)$. Further, we introduce paths in the complex plane on which the $\varphi^j$'s are $L^2$ (in fact, decay exponentially), see Proposition \ref{cor:varphiL2}. The objective of the next steps is to derive an equivalent of the HSM formulation for these ``complex-scaled'' traces.
	\item  For real wavenumbers  equations (\ref{eq:def_traces2}) and (\ref{eq:hpr}) %-(\ref{eq:hpr_kernel})
 provide half-plane representations of the solution $u$ in terms of the traces $\varphi^j$, $j\in \llbracket0,3\rrbracket$.
	The magic result is that the solution $u$ can also be represented in terms of the complex-scaled traces, see Theorem \ref{th:deformed_representation}. The price to pay is that the $j$th representation, for $j\in \llbracket0,3\rrbracket$, is valid only in a part, which depends on the path chosen in step 1, of the corresponding half-plane $\Omega^j$. These new representations are deduced from the initial ones (\ref{eq:def_traces2}) % and (\ref{eq:hpr})}%(\ref{eq:hpr})-(\ref{eq:hpr_kernel})
by applying Cauchy's integral theorem.
	\item Fortunately, the part of $\Omega^j$ where the new representation holds contains the half-lines $\Sigma^{j\pm 1}\cap\Omega^j$. By complexifying  -- by which we mean analytically continuing --  this new representation into the complex plane from the half-lines $\Sigma^{j\pm 1}\cap\Omega^j$, we can derive compatibility conditions for the complex-scaled traces which constitute a complex-scaled version of the HSM formulation in an $L^2$ setting. Fredholmness of this formulation can be proven using similar arguments as for the standard HSM for complex wavenumbers, see Theorem \ref{th:HSMM_cmplxed}.
	\item Once the complex-scaled traces are computed the solution can be reconstructed using the new representations in terms of the complex-scaled traces established in step 2.
\end{enumerate}

Let us mention that, while the initial motivation for the complexification was a theoretical one, it turns out that the new formulation is very attractive computationally, because of the fast decay at infinity of the complex-scaled traces.  Let us note also that this idea of complexification is potentially valuable for computation in the dissipative case too, and it is likely that the formulation in the non-dissipative case could be derived from a complex-scaled formulation in the dissipative case by a limiting absorption argument. (Something similar has been done in the context of scattering by wedges in \cite{CrLe:99,Ka:05}.) An attraction of such a derivation would be that traces are in $L^2$ at each step of the derivation.

\subsection{The complex-scaled traces}
The construction of the so-called \emph{complex-scaled traces} is based on an analyticity property of any solution $u\in H^1_{\mathrm{loc}}(\Omega)$ to \eqref{eq:Sommerfeld}-\eqref{pb:probleme_Dir_real}, which can be derived as follows. Following \cite{BraWer:65} (and see \cite[Theorem 2.27, Corollary 2.28]{ChGrLaSp:11}), $u$ can be expressed as a combined single- and double-layer potential on $\Sigma_a$, i.e. as
\begin{equation} \label{eq:ansatz}
u(\bsx) = \cD\phi (\bsx) - \ri k \cS\phi(\bsx), \quad \bsx\in \Omega,
\end{equation}
for some $\phi\in H^{1/2}(\Sigma_a)$ (specified below). Here $\cS\phi$ and $\cD\phi$ are the (acoustic) single- and double-layer potentials, respectively, with density $\phi$, defined for $\phi\in L^2(\Sigma_a)$ by
$$
\cS\phi (\bsx) := \int_{\Sigma_a}\Phi(\bsx,\bsy)\phi(\bsy)\, \rd s(\bsy), \quad \cD\phi (\bsx) := \int_{\Sigma_a}\frac{\partial \Phi(\bsx,\bsy)}{\partial n(\bsy)}\phi(\bsy)\, \rd s(\bsy), \quad \bsx\in \Omega,
$$
where the normal $n$ is directed into $\Omega$ and $\Phi$ is the outgoing Green's function of the Helmholtz equation given in \eqref{eq:Greenfct}.
The function $u$ defined in \eqref{eq:ansatz} satisfies the Helmholtz equation \eqref{pb:probleme_Dir_real} and the Sommerfeld radiation condition \eqref{eq:Sommerfeld} for any choice of $\phi\in H^{1/2}(\Sigma_a)$ ({in fact}, any $\phi\in L^2(\Sigma_a)$),  and (see \cite[\S2.6]{ChGrLaSp:11}), satisfies the boundary condition $u=g$  on $\Sigma_a$ provided
\begin{equation} \label{eq:bie}
A\phi = g,
\end{equation}
where $A\phi$ is defined for $\phi\in L^2(\Sigma_a)$ and almost all $\bsx\in {\Sigma_a}$ by
\[
	A\phi(\bsx) := {\frac{\phi(\bsx) }{2}} + \int_{\Sigma_a}\left[\frac{\partial \Phi(\bsx,\bsy)}{\partial n(\bsy)}\,-\ri k  \Phi(\bsx,\bsy)\right]\,\phi(\bsy) \rd s(\bsy),
\]
with the integral  understood as a Cauchy principal value. Since $A:H^s({\Sigma_a})\to H^s({\Sigma_a})$  is invertible  for $0\leq s\leq 1$ \cite[Corollary  2.8]{Convex}, in particular for $s=1/2$, \eqref{eq:bie} has a unique solution $\phi\in H^{1/2}({\Sigma_a})$.

For a given $j \in \llbracket0,3\rrbracket,$ we apply \eqref{eq:ansatz} for $\bsx \in \Sigma^j \setminus \Sigma^j_a,$ and we use the coordinate system $(x_1^j,x_2^j)$ defined in \eqref{def:coord-loc-HP-ac-milieu}. Defining $\bsx^j(t) := (a,t)$ with $|t| > a,$ this yields, by definition of $\varphi^j$,
\begin{equation} \label{eq:phijrep}
\varphi^j(t) = \cD^j\phi(t)-\ri k \cS^j\phi(t),
\end{equation}
for real $t$ with $|t|>a$, where
\begin{equation} \label{eq:SDjdef}
\cD^j\phi(t) := \int_{\Sigma_a} \frac{\partial \Phi(\bsx^j(t),\bsy^j)}{\partial n(\bsy^j)}\phi(\bsy^j)\rd s(\bsy^j), \;\; \cS^j\phi(t):=\int_{\Sigma_a} \Phi(\bsx^j(t),\bsy^j) \phi(\bsy^j)\, \rd s(\bsy^j),
\end{equation}
and we recall from \eqref{eq:Greenfct} that
\[
 \Phi(\bsx^j(t),\bsy^j) = \frac{\ri}{4}\,H^{(1)}_0\big(  k R(a-y_1^j,t-y_2^j) ),
\]
where $R$ is defined in \eqref{eq:R(s,t)}.

Let us use \eqref{eq:phijrep} to prove that the function $\varphi^j$, defined by \eqref{eq:def_traces}, can be continued analytically into the complex plane from $(-\infty,-a)$ and from $(a,+\infty)$. Consider a fixed $\bsy^j \in \Sigma_a$. The function {$z \mapsto R(a-y_1^j,z-y_2^j)$} has an analytic continuation from {$(a,+\infty)$} (respectively, {$(-\infty,-a)$}) to the complex half-plane $\Re(z) > a$ (respectively, $\Re(z) < -a$). Indeed, to obtain this analytic continuation we simply have to use, in the {definition \eqref{eq:R(s,t)} of $R(a-y_1^j,z-y_2^j)$} for real $z$, the principal square root of a complex number, which we will denote by $z^{1/2}$ or $\sqrt{z}$, defined as
\begin{equation*}
\sqrt{z} := |z|^{1/2} \ \re^{\ri \Arg(z) / 2}
\quad\text{with } \Arg(z) \in (-\pi,+\pi],
\end{equation*}
which is analytic in $\C \setminus \R^-$, where $\R^-:=(-\infty,0]$. The analyticity of  $z \mapsto R(a-y_1^j,z-y_2^j)$ follows by noticing that $(y_1^j-a)^2 + (y_2^j-z)^2 \in \C \setminus \R^-$ if $\Re(z) > a$ (respectively, if $\Re(z) < -a$), since $y_2^j \in [-a,+a]$. Since also $z \mapsto H^{(1)}_0(z)$ is analytic in $\Re(z) > 0,$ we conclude that the function $z \mapsto \Phi(\bsx^j(z),\bsy^j)$ is analytic in ${|\Re(z)| > a}$. And the same arguments and conclusion apply also for $z \mapsto \partial \Phi(\bsx^j(z),\bsy^j) / \partial n(\bsy^j)$. Finally, using standard results about analyticity of functions defined as integrals (e.g.\ \cite[Corollary X.3.19]{Amann:09}), we conclude that $z\mapsto \cS^j\phi(z)$, $z\mapsto \cD^j\phi(z)$, and so also $z \mapsto \varphi^j(z)$, have analytic continuations from {$(a,+\infty)$} (respectively, {$(-\infty,-a)$}) to the complex half-plane $\Re(z)> a$ (respectively, $\Re(z) < -a$).

The behavior of these analytic continuations as $|z| \to +\infty$ depends on $\Im(z).$ Indeed, for $m \in \N:=\{0,1,...\},$ we have \cite[Equation (9.2.30)]{AbramowitzStegun}
\begin{equation} \label{eq:HankAsym}
	H^{(1)}_m(z) = \sqrt{\frac{2}{\pi z}}\ \re^{\ri (z - m\pi/2 - \pi/4)} \Big( 1 + \mathcal{O}\big( |z|^{-1}\big) \Big)
	\quad\text{as } |z| \to +\infty,
\end{equation}
uniformly in $\Arg(z)$ for $|\Arg(z)|< \zeta$, for every $\zeta<\pi$. Further, as a consequence of Lemma \ref{lem:Im(R-z)} and Remark \ref{rem:Im(R+z)},
\begin{eqnarray} \label{eq:Rasymp1}
{R(a-y_1^j,z-y_2^j)}  &=& z-y_2^j + \mathcal{O}(|z|^{-1}) \quad \mbox{ in } \Re(z) >a, \;\;\mbox{ and }\\ \label{eq:Rasymp2}
{R(a-y_1^j,z-y_2^j)}  &=& y_2^j-z + \mathcal{O}(|z|^{-1}) \quad \mbox{ in } \Re(z) <-a,
\end{eqnarray}
as $|z|\to+\infty$, uniformly in $\Arg(z)$ {and $\bsy^j$, for $\bsy^j\in \Sigma_a$}. {Using \eqref{eq:phijrep}} it follows from the above asymptotics that
\begin{eqnarray} \label{eq:ExpDecay1}
\varphi^j(z) &=& \mathcal{O}\left(\re^{-k\Im(z)}\, |z|^{-1/2}\right), \quad \mbox{ in } \Re(z)>a, \mbox{ and }\\  \label{eq:ExpDecay2}
\varphi^j(z) &=& \mathcal{O}\left(\re^{k\Im(z)}\, |z|^{-1/2}\right), \quad \mbox{ in } \Re(z)<-a,
\end{eqnarray}
as $|z|\to +\infty$, uniformly with respect to $\Arg(z)$, and we note that $\Im({R(a-y_1^j,z-y_2^j)})$ and $\Im(z)$ have the same sign if $\Re(z)> a$ (opposite signs if $\Re(z) < -a$).
Thus $\varphi^j(z)$ is exponentially decreasing in the quadrants $\big\{\Re(z) > a \text{ and  } \Im(z) > 0\big\}$ and $\big\{\Re(z) < -a \text{ and  } \Im(z) < 0\big\}.$

Our idea is to choose half-lines in these quadrants and consider the analytic continuations of the $\varphi^j$ on these half-lines as new unknowns instead of the initial traces. In other words, choosing some $\theta \in (0,\pi/2)$, we introduce the complex path (Figure \ref{fig:path}) parameterized by
\begin{equation}
\label{eq:parametrage_t_theta}
	z = \tau_\theta(s) := \left\{
	\begin{array}{ll}
		-a + (s+a) \, \re^{\ri\theta} & \text{if }s < -a, \\
		s & \text{if } -a \le s \le +a, \\
		a + (s-a) \, \re^{\ri\theta} & \text{if }s>a,		
	\end{array}
    \right.
\end{equation}
and we define the complex-scaled traces by
\begin{equation}
\label{eq:definition_phi_theta}
	\varphi^j_\theta (s) := \varphi^j \big(\tau_\theta(s)\big)
	\quad\text{for } s \in \R \text{ and } j \in \llbracket0,3\rrbracket.
\end{equation}

\begin{rem}
We have chosen particular complex paths in the quadrants $\big\{\Re(z) > a \text{ and  } \Im(z) > 0\big\}$ and $\big\{\Re(z) < -a \text{ and  } \Im(z) < 0\big\}$, given by \eqref{eq:parametrage_t_theta}, that move into the complex plane already from $\pm a$, the corners of $\Sigma_a$.  Note that it is possible, alternatively, to start to complexify at a positive distance from the corners, i.e.\ from $\pm b$, for some $b>a$. It is also possible to choose a smoother complex change of variable as usually done in PML methods.
\end{rem}
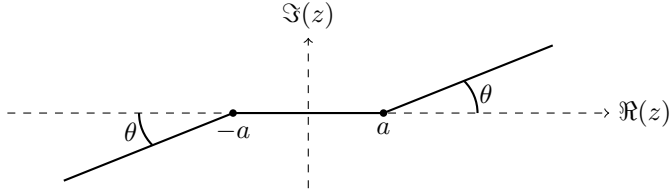
\begin{figure}[h!]
\begin{center}
	\begin{tikzpicture}[scale=.5]
	\draw[->,dashed] (-8,0) -- (8,0) ;
	\draw[->,dashed] (0,-2) -- (0,2) ;
	%\draw[thick] (-7,0)--(7,0);
	\draw[thick] (-6.5,-1.8)--(-2,0) -- (2,0) -- (6.5,1.8);
	\draw[thick]  (4.5,0) arc (0:45:1.2);
	\draw (4.3,.5) node[right]{$\theta$};
		\draw[thick]  (-4.5,0) arc (180:225:1.2);
	\draw (-4.3,-.5) node[left]{$\theta$};
	\fill (2,0) circle (0.1);
	\draw (2,0) node[below] {$a$};
	\fill (-2,0) circle (0.1);
	\draw (-2,0) node[below] {$-a$};
	\draw(8,0) node[right] {$\Re (z)$};
	\draw(0,2) node[above] {$\Im (z)$};
	\end{tikzpicture}
\end{center}
\caption{The complex path $s\rightarrow \tau_\theta(s)$}
\label{fig:path}
\end{figure}

We note that, by definition, $\varphi^j_\theta|_{(-a,a)}\in L^2(-a,a)$, since $\varphi^j_\theta(t) = \varphi^j(t)=u(\bsx^j(t))$ $=g(\bsx^j(t))$, for $-a<t<a$, and $g\in H^{1/2}(\Sigma_a)\subset L^2(\Sigma_a)$. Moreover, thanks to the analyticity of $\varphi^j$, it is clear that {the restrictions of $\varphi^j_\theta$ to $(a,+\infty)$ and $(-\infty,-a)$ are continuous.} %|_{(a,\infty)}\in C(a,\infty)$ and $\varphi^j_\theta|_{(-\infty,-a)}\in C(-\infty,a)$.
Further, it follows from (\ref{eq:ExpDecay1}-\ref{eq:ExpDecay2}) that
\begin{equation} \label{eq:varphiAsymp}
\varphi^j_\theta (s) = \mathcal{O}\left(\re^{-k|s|\sin\theta}|s|^{-1/2}\right), \quad \mbox{as } |s|\to+\infty.
\end{equation}
 Thus whether or not $\varphi^j_\theta\in L^2(\R)$ depends on the behaviour of $\varphi^j_\theta(t)$ as $t\to a^+$ and $t\to -a^-$.

 {The following propositions} bound $\varphi^j_\theta(t)$ on $|t|>a$, in particular near $\pm a$, and show that $\varphi^j_\theta \in L^2(\R)$ for $0<\theta<\pi/2$. {We relegate some of the technical details to Appendix \ref{sec:ComplIntOper}. Applying the bounds from Proposition \ref{lem:SDbounds} to (the analytic continuation of) \eqref{eq:phijrep} we obtain the following proposition}, on observing, from \eqref{eq:bie}, that $\phi = A^{-1}g$ and that \cite[Corollary  2.8]{Convex} $A^{-1}$ is bounded as an operator on $L^2(\Sigma_a)$.
\begin{prop} \label{cor:localbound}
For every $\theta\in (0,\pi/2)$ and $j\in\llbracket0,3\rrbracket$ there exists a constant $C>0$, that depends only on $a$, $k$, and $\theta$, such that
\begin{eqnarray*}
|\varphi^j(z)| \leq  \left\{\begin{array}{ll}
                                C|z-a|^{-1/2}\re^{-k\Im(z)}\,\|g\|_{L^2(\Sigma_a)}, & \mbox{if  }\Re(z) > a  \mbox{ with }|\mathrm{Arg}(z-a)|\leq \theta, \\
                                C|z+a|^{-1/2}\re^{k\Im(z)}\,\|g\|_{L^2(\Sigma_a)}, & \mbox{if } \Re(z) < -a \mbox{ with }|\mathrm{Arg}(-z-a)|\leq \theta.
                              \end{array}\right.
\end{eqnarray*}
\end{prop}

The bound in the above proposition implies that $\varphi^j_\theta\in L^1(\R)$, for $0<\theta<\pi/2$, and is sufficient for our Cauchy's integral formula arguments below in \S\ref{sec:deformed halfspace}. But it is not quite strong enough to establish $\varphi^j_\theta\in L^2(\R)${, as discussed in the proof of  Proposition \ref{prop:SDbounded} in Appendix \ref{sec:ComplIntOper}.}

{It follows from (the analytic continuation of) \eqref{eq:phijrep} that, in the notation of Proposition \ref{prop:SDbounded}, $\varphi^j_\theta(s) = \cD^j_\theta \phi(s)-\ri k \cS^j_\theta\phi(s)$, for $|s|>a$, where $\phi= A^{-1}g$. Thus, and arguing as above Proposition \ref{cor:localbound}, we deduce the following result from the above proposition and Proposition \ref{prop:SDbounded}.}
\begin{prop} \label{cor:varphiL2}
For $0<\theta<\pi/2$ and $j\in\llbracket0,3\rrbracket$, $\varphi^j_\theta \in L^2(\R)$. Further, for some constant $C>0$ depending only on $\theta$, $a$, and $k$,
\begin{equation*} %\label{eq:varphidecay}
|\varphi^j_\theta(s)| \leq C(|s|-a)^{-1/2}\re^{-k\sin(\theta)(|s|-a)}\|g\|_{L^2(\Sigma_a)}, \quad |s|>a,
\end{equation*}
and $\|\varphi^j_\theta\|_{{L^2(\R)}}\leq C\|g\|_{L^2(\Sigma_a)}$.
\end{prop}

\subsection{The deformed half-plane representations} \label{sec:deformed halfspace}
We have introduced in the previous section the complex-scaled traces $\varphi_{\theta}^j$ and proved that they belong to $L^2(\R)$. The objective now is to derive a HSM formulation for these new unknowns. The first step is to establish new representation formulas for the solution $u$ of (\ref{pb:probleme_diff_initial})-(\ref{eq:Sommerfeld}) in the half-planes using these complex-scaled traces instead of the original traces $\varphi^j$.

We recall that the solution $u$ of (\ref{pb:probleme_diff_initial}-\ref{eq:Sommerfeld}) can be represented in each half-plane in terms of its traces as
\begin{equation} \label{eq:hs}
 u({\bsx}^j) = \int_{\R}h(x_1^j-a,x_2^j-y_2^j)\,\varphi^j(y_2^j)\, \rd y_2^j,\quad  {\bsx}^j \in \Omega^j,
\end{equation}
where the kernel $h$ is defined by (\ref{eq:hpr_kernel}).
%\[
% h(x_1,x_2,y_2)=\frac{i\omega(x_1-a)}{2}\frac{H^{(1)}_1(\omega R)}{R},\quad R=[(x_1-a)^2+(x_2-y_2)^2]^{1/2}.
% \]
 Our objective is to derive a similar formula using the complex-scaled trace $\varphi^j_\theta$ instead of the trace $\varphi^j$. This can be done by deforming the path of integration into the complex plane. This leads to the following crucial result.
\begin{theorem}\label{th:deformed_representation}
Let $u$ be the solution of (\ref{pb:probleme_diff_initial})-(\ref{eq:Sommerfeld}) and $\varphi^j_\theta$ be defined as in (\ref{eq:definition_phi_theta}). For $0<\theta<\pi/2$  we have
\begin{equation}
\label{eq-repinOmegatheta}
u({\bsx}^j) = \int_{\R}h(x_1^j-a,x_2^j-\tau_\theta(s))\,\varphi^j_\theta(s)\,\tau'_\theta(s)\, \rd s,\quad {\bsx}^j\in \Omega^j_\theta,\quad j\in \llbracket0,3\rrbracket,
\end{equation}
where
\begin{equation}\label{eq:Omegaj_theta}
	\Omega^j_\theta:=\{{\bsx}^j=(x_1^j,x_2^j)\in \Omega^j: \;x_1^j-a>(|x_2^j|-a)\tan(\theta)\}.
\end{equation}
\end{theorem}
\begin{rem}
%	\label{rem-Omega_theta}
	Let us point out that the new representation formula (\ref{eq-repinOmegatheta}) is valid only in a subdomain $\Omega^j_\theta$ of the half-plane $\Omega^j$ (Figure \ref{fig:domainsvalid}). The larger the angle $\theta$, the faster the decay of the complex-scaled traces $\varphi^j_\theta$ (Proposition \ref{cor:varphiL2}) but the smaller the domain of validity $\Omega^j_\theta$ of the new representation (\ref{eq-repinOmegatheta}).
\end{rem}
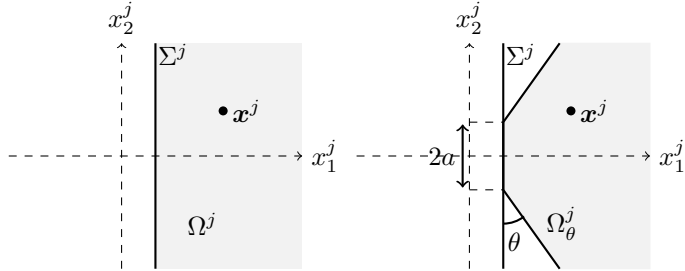
\begin{figure}[h!]
	\begin{center}
		\begin{tikzpicture}[scale=.3]
		\fill[color=gray!10] (1.5,-5)rectangle(8,5);
		\draw[->,dashed] (-5,0) -- (8,0) ;
		\draw[->,dashed] (0,-5) -- (0,5) ;
		\draw(8,0) node[right] {$x_1^j$};
		\draw(0,5) node[above] {$x_2^j$};
		\draw (1.2,4.5) node[right] {${\Sigma^j}$};
		\draw (2.5,-3) node[right] {${\Omega^j}$};
		\draw[thick] (1.5,-5)--(1.5,5);
		%		\draw[<->] (.1,-3)--(1.4,-3);
		%		\draw(.7,-3) node[below] {$a$};
		\fill (4.5,2) circle (0.2);
		\draw (4.5,2) node[right] {$\bsx^j$};
		\end{tikzpicture}	\begin{tikzpicture}[scale=.3]
		\filldraw[color=gray!10] (1.5,1.5)--(4,5)--(8,5)--(8,-5)--(4,-5)--(1.5,-1.5)--(1.5,1.5);
		\draw[->,dashed] (-5,0) -- (8,0) ;
		\draw[->,dashed] (0,-5) -- (0,5) ;
		\draw(8,0) node[right] {$x_1^j$};
		\draw(0,5) node[above] {$x_2^j$};
		\draw (1.2,4.5) node[right] {${\Sigma^j}$};
		\draw (3,-3) node[right] {${\Omega^j_\theta}$};
		\draw[thick] (1.5,-5)--(1.5,5);
		\draw[<->,thick] (-.3,-1.4)--(-.3,1.4);
		\draw(-.2,0) node[left] {$2a$};
		\draw[-,dashed] (0,1.5) -- (1.5,1.5) ;
		\draw[-,dashed] (0,-1.5) -- (1.5,-1.5) ;
		%		\draw[<->] (.1,-3)--(1.4,-3);
		%		\draw(.7,-3) node[below] {$a$};
		\fill (4.5,2) circle (0.2);
		\draw (4.5,2) node[right] {$\bsx^j$};
		\draw[thick] (4,-5)--(1.5,-1.5)--(1.5,1.5)--(4,5);
		\draw[thick]  (1.5,-3) arc (-90:-52:1.5);
		\draw  (2,-3) node[below]{$\theta$};
		\end{tikzpicture}
	\end{center}
	\caption{The domains of validity of the half-plane representations: left is without complex-scaling ($\Omega^j$); right is with complex-scaling ($\Omega_\theta^j$).}
\label{fig:domainsvalid}
\end{figure}
\begin{proof}
To derive \eqref{eq-repinOmegatheta} from \eqref{eq:hs}
 it suffices to show (note the symmetry $\tau_\theta(-s)=-\tau_\theta(s)$) that, for all $ {\bsx}^j\in \Omega^j_\theta$,
	\[
 \int_{a}^{+\infty}h(x_1^j-a,x_2^j-y_2^j)\,\varphi^j(y_2^j)\, \rd y_2^j= \int_{a}^{+\infty}h(x_1^j-a,x_2^j-\tau_\theta(s))\,\varphi^j_\theta(s)\,\tau'_\theta(s)\, \rd s
	\]
and%\scednote{I think to make clear where the definition of $\Omega^j_\theta$ comes from it is really helpful to have this extra equation}
	\[
 \int_{a}^{+\infty}h(x_1^j-a,x_2^j+y_2^j)\,\varphi^j(-y_2^j)\, \rd y_2^j= \int_{a}^{+\infty}h(x_1^j-a,x_2^j+\tau_\theta(s))\,\varphi^j_\theta(-s)\,\tau'_\theta(s)\, \rd s.
	\]
%\textcolor{red}{and a similar identity on $(-\infty,-a)$ which follows analogously (note the symmetry $\tau_\theta(-s)=\tau_\theta(s)$).} \scednote{I couldn't quite see what was meant here by "by symmetry and by the definition of $\tau_\theta$" this is my interpretation.}
For $0<\delta<M$ and $0<\theta<\pi/2$ we introduce the complex domains
\begin{eqnarray*}
		T_\theta &:=& \{ z=a+r \re^{\ri \alpha}:\;0<\alpha<\theta,\;r>0\} \;\;\text{ and }\\ %\nonumber
 T_{\theta}^{\delta,M}&:=&\{ z=a+r \re^{\ri \alpha}:\;0<\alpha<\theta,\;\delta<r<M\}
 \end{eqnarray*}
	(see Figure \ref{fig:contour}). First we show that, if ${\bsx}^j=(x_1^j,x_2^j)\in \Omega^j$, then
	\begin{equation} \label{eq:equivalence}
	{\bsx}^j\in \Omega^j_\theta\;\Leftrightarrow \; \left\{\begin{array}{ll}z\mapsto h(x_1^j-a,x_2^j\pm z) \text{ is analytic in } T_{\theta}\\ \text{and continuous in }\overline{T_{\theta}}.\end{array}\right.
	\end{equation}
	Indeed, for each ${\bsx}^j$,  the function $z\mapsto h(x_1^j-a,x_2^j-z)$ has two branch points $z_\pm$, the points where $(x_1^j-a)^2+(x_2^j-z)^2$ vanishes, given by $z_\pm=x_2^j\pm \ri(x_1^j-a)$. %, by definition of $T_{\delta,M}$,
If ${\bsx}^j=(x_1^j,x_2^j)\in \Omega^j$ these branch points are outside $\overline{T_{\theta}}$ if and only if $x_1^j-a>(x_2^j-a)\tan(\theta)$; similarly, the branch points of the mapping $z\mapsto h(x_1^j-a,x_2^j+z)$ are outside $\overline{T_{\theta}}$ if and only if $x_1^j-a>(-x_2^j-a)\tan(\theta)$.

Thus, and since also $z \mapsto \varphi^j(\pm z)$ is analytic in   $T_\theta^{\delta,M}$ and continuous in $\overline{T_\theta^{\delta,M}}$, for $0<\delta<M$,
 applying Cauchy's integral theorem we have
	\[
		 \int_{\partial T_\theta^{\delta,M}}h(x_1^j-a,x_2\mp z)\,\varphi^j( \pm z)\, \rd z=0,\quad{\bsx}^j\in \Omega^j_\theta.
	\]
To complete the proof, we have to show that
	\[
		\lim_{\delta\rightarrow0}\int_{0}^\theta h(x_1^j-a,x_2^j\mp {(a+\delta \re^{\ri \alpha})})\varphi^j(\pm{(a+\delta \re^{\ri \alpha})})\ri \delta \re^{\ri \alpha}\,\rd \alpha =0
		\]
		and
	\[
		\lim_{M\rightarrow+\infty}\int_{0}^\theta h(x_1^j-a,x_2^j\mp {(a+M\re^{\ri \alpha})})\varphi^j(\pm {(a+M\re^{\ri \alpha})})\ri M\re^{\ri \alpha}\,\rd \alpha =0.
		\]
These two limits are a consequence of Proposition \ref{cor:localbound},  since the constraint $\bsx^j\in \Omega_\theta^j$ ensures, by \eqref{eq:equivalence}, that $h(x_1^j-a,x_2^j \mp z)$ is a continuous {function of $z$} in $\overline{T_\theta}$, and the bound \eqref{eq:Hbound1} and the asymptotics \eqref{eq:Rasymp1}-\eqref{eq:Rasymp2} imply that $h(x_1^j-a,x_2^j \mp z) =\mathcal{O}\left({|z|^{-1/2}}\right)$ as $z\to +\infty$ in $T_\theta$, uniformly in $\Arg(z)$.
	\end{proof}

\begin{figure}[h!]
	\begin{center}
		\begin{tikzpicture}[scale=.6]
		\draw[->,dashed] (-1,0) -- (8,0) ;
		\draw[->,dashed] (0,-3) -- (0,5) ;
		\draw[-,thick] (2.5,0)--(7,0);
		\draw(8,0) node[right] {$\Re e\,z$};
		\draw(0,5) node[above] {$\Im m\,z$};
		\fill (4,0) circle (0.1);
		\draw (4,0) node[below] {$x_2^j$};
		\fill (4,2) circle (0.1);
		\fill (4,-2) circle (0.1);
		\draw[-,thick] (4,-1.8)--(4,-3.5);
		%	\draw (2.8,2.5) node[above] {Branch point of $h(x_1^j,x_2^j,z)$};
		%	\draw  (3,1.8) node[above] { of $h$};
		\draw[-,thick] (4,1.8)--(4,4.5);
		\draw[dashed] (4,0) -- (4,2) ;
		\draw[dashed] (0,2) -- (4,2) ;
		\draw[dashed] (4,-.8) -- (4,-2) ;
		\draw[dashed] (0,-2) -- (4,-2) ;
		\draw (0,2) node[left] {${x_1^j}-a$};
		\draw (0,-2) node[left] {$-(x_1^j-a)$};
		\draw[thick] (2.4,0.25) -- (6.5,1.8);
		\draw[dashed]  (4.5,0) arc (0:31:1.9);
		\draw (4.5,.5) node[right]{$\theta$};
		\draw[thick]  (7,0) arc (0:31:3.5);
		\draw[thick]  (2.5,0) arc (0:31:.5);
		\fill (2,0) circle (0.1);
		\draw (2,-.15) node[below] {$a$};
			\fill (7,0) circle (0.1);
		\draw (7,-.15) node[below] {$a+M$};
	%	\draw (2.4,-.15) node[below] {$a+\delta$};
		
		\end{tikzpicture}
	\end{center}
	\caption{The contour $\partial T_\theta^{\delta,M}$ and the branch points $z_\pm=x_2^j\pm \ri (x_1^j-a)$
of $h(x_1^j-a,x_2^j-z)$.}
\label{fig:contour}
\end{figure}
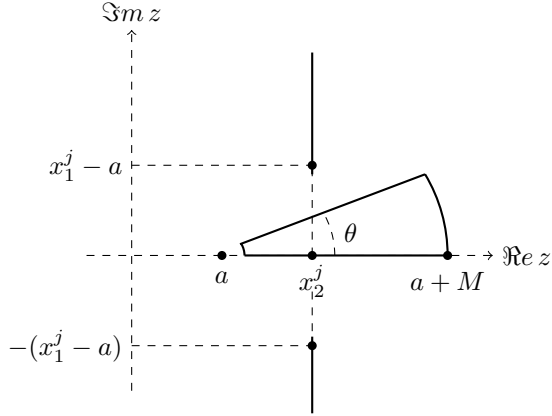

\subsection{Derivation and analysis of the complex-scaled HSM method}	\label{sec:Derivation}
	We derive now an analogue of \eqref{eq:syst_comp} for the complex-scaled traces $\varphi^j_\theta$. %We proceed in the same way to derive compatibility relations.
We know, thanks to Theorem \ref{th:deformed_representation}, that
	\[
		 u({\bsx}^0) = \int_{\R}h(x_1^0-a,x_2^0-\tau_\theta(s))\,\varphi^0_\theta(s)\,\tau'_\theta(s)\, \rd s,\quad {\bsx}^0\in \Omega^0_\theta.
	\]
	As $\Omega^0\cap\Sigma^1\subset \Omega^0_\theta$, for all $\theta\in(0,\pi/2)$, this holds, in particular, when $\bsx^0=(x^0_1,x^0_2)\in \Omega^0\cap\Sigma^1$, i.e.\ for $x_1^0=-x_2^1>a$ and $x_2^0=x_1^1= a$, so that (cf.\ \eqref{eq:comp1a}) %and  and  $x_1^1=a$, for we have in particular (since for any point of $\Sigma^1\cap\Omega^0$, $x_2^1=-x_1^0<-a$ and $x_1^1=x_2^0=a$)
	\begin{equation}\label{eq:HSMM_step1}
		\varphi^1(x_2^1)=\int_{\R}h(-x_2^1-a,a-\tau_\theta(s))\,\varphi^0_\theta(s)\,\tau'_\theta(s)\, \rd s,\quad x_2^1<-a.
	\end{equation}
% Since by definition,  $\varphi^1_\theta(t)=\varphi^1(\tau_\theta(t))$, we get the new complex scaled compatibility relation
% 	\[
% 		\varphi^1_\theta(t)=\int_{\R}h(-\tau_\theta(t),a,\tau_\theta(s))\,\varphi^0_\theta(s)\,\tau'_\theta(s)\, \rd s,\quad\text{for}\;t<-a.
% 	\]
% 	To justify the previous identity, one has to check that the function $$z\rightarrow \int_{\R}h(z,x_2^0,\tau_\theta(s))\,\varphi^0_\theta(s)\,\tau'_\theta(s)\, \rd s$$ is analytic in a neighboorhood of the complex path $[a,a+\re^{i\theta}\R^+]$. This amounts to check that the quantity $(\tau_\theta(t)-a)^2+(\tau_\theta(s)-a)^2$ never vanishes for $t>a$ and $s\in\R$.
%
	Remember that $\varphi^1(x^1_2)$  is analytic in $\Re (x^1_2)<-a$ and, by the definition \eqref{eq:definition_phi_theta}, $\varphi^1_\theta(t)$ for $t<-a$ is the value of the analytic continuation of $\varphi^1$ at $x^1_2=\tau_\theta(t)$.

To obtain the new complex-scaled compatibility relation, the idea is to complexify $x^1_2$ in \eqref{eq:HSMM_step1}, i.e.\ simply to substitute $x^1_2=\tau_\theta(t)$. By the uniqueness of analytic continuation, this is valid provided the right-hand side is analytic as a function of $x^1_2$ in a connected domain containing the half-line $x_2^1<-a$ and the half-line $x_2^1=\tau_\theta(t)$, $t<-a$. We check this in the next lemma. %To do that, one has first to check the analyticity of the right-hand side as a function of $x^1_2$.
	\begin{lemma}\label{lem:analyticity_utheta}
	Let $\psi\in L^2(\R)$. For $0<\theta<\pi/2$ the function
\begin{equation} \label{eq:hopeanalytic}
z\rightarrow \int_{\R}h(z-a,a-\tau_\theta(s))\psi(s)\tau'_\theta(s) ds
\end{equation}
	is analytic in  the domain $G_\theta:=\{z\in\C:\,z\neq a,\,-\pi/2+\theta<\Arg(z-a)<\pi/2\}$.
	\end{lemma}
	\begin{proof}
	We will prove that this function is analytic using standard result{s} about analyticity of functions defined as integrals based on the dominated convergence theorem {(e.g.,} \cite[Corollary X.3.18]{Amann:09}{)}.
	By the definition \eqref{eq:hpr_kernel}, for each $s\in\R$ the kernel $h(z-a,a-\tau_\theta(s))$ is locally an analytic function of $z$ wherever    the quantity $$R(z-a,\tau_\theta(s)-a)=[(z-a)^2+(\tau_\theta(s)-a)^2]^{1/2}$$ does not vanish. Lemma \ref{lem-Rbound-z-tau} shows that, for every $\theta_0\in (\theta,\pi/2)$,
\begin{equation} \label{eq:Rlower}
	|R(z-a,\tau_\theta(s)-a)|^2\geq \cos(\theta_0)|z-a|^2
\end{equation}
if $-\theta_0+\theta\leq\Arg(z-a)\leq\theta_0$.
	%where $m(\theta,z)=\min (\cos(\Arg(z-a)),\cos(\theta-\Arg(z-a)))$ is a continuous positive function of $z$, which
	%never vanishes in the domain $\{z\in\C\,z\neq a,\,-\pi/2+\theta<\Arg(z-a)<\pi/2\}$.
Consequently,  $R(z-a,\tau_\theta(s)-a)$ does not vanish in $G_\theta$. % so that $z\mapsto h(z-a,a-\tau_\theta(s))$ is analytic in $G\theta$, so every $s\in \R$. This implies the required result in a domain $U\subset G_\theta$ provided the integrand on the right hand side of \eqref{eq:hopeanalytic} is bounded by some $L^1$ function, for every $s\in\R$.%$\{z\in\C:z\neq a,\,-\pi/2+\theta<\Arg(z-a)<\pi/2\}$.}

{Further, \eqref{eq:Rlower} implies that,}  for every bounded subdomain $U$ of $G_\theta$ that is bounded away from $a$, %the domain $\{z\in\C:\,z\neq a,\,-\pi/2+\theta<\Arg(z-a)<\pi/2\}$,
there exists $m_U>0$ such that
	\begin{equation} \nonumber %\label{eq-lowerboundR}
|R(z-a,\tau_\theta(s)-a)|\geq m_U,\quad z\in U, \;\; s\in \R.	
\end{equation}
	%Let us first check that the quantity $(z-a)^2+(\tau_\theta(s)-a)^2$ never vanishes for $z$ in  $-\pi/2<\Arg(z-a)<\pi/2$  and $s\in\R$,  except when $z=s=a$. This can be done by considering first the imaginary part of $(z-a)^2+(\tau_\theta(s)-a)^2$. Because $\theta<\pi/2$, the imaginary part of $(z-a)^2$ is non negative, and it vanishes only for $z$ on the real axis. Likewise, the imaginary part of $(\tau_\theta(s)-a)^2$ is non negative, and it vanishes only for $|s|\leq a$. Finally for real $z$ and $s$, it is easy to conclude.
{Moreover,} for any such subdomain it follows from Lemma \ref{lem:Im(Rhat-z)}  that, for some constant $C>0$ independent of $s\in\R$ and $z\in U$,
$$\Im(R(z-a,\tau_\theta(s)-a))\geq |s|\sin(\theta)-C.$$
From the above bounds, and the bound \eqref{eq:Hbound1}, it follows that
$$
|h(z-a,a-\tau_\theta(s))|\leq C'\re^{- k \sin(\theta)|s|},
$$
for some constant $C'>0$ independent of $s\in\R$ and $z\in U$. Thus, for every $z\in U$ and $s\in \R$, the integrand in \eqref{eq:hopeanalytic} has modulus $\leq H(s)$, where $H\in L^1(\R)$ is defined by
$$
H(s):= C'\re^{- k \sin(\theta)|s|}|\psi(s)|, \quad s\in \R.
$$

This domination property, and the analyticity of $z\mapsto h(z-a,a-\tau_\theta(s))$ in $U\subset G_\theta$, imply {(e.g.,} \cite[Corollary X.3.18]{Amann:09}{)} that the function \eqref{eq:hopeanalytic} is analytic in $U$, for every $U$, and so analytic in $G_\theta$.
\end{proof}

%Now let us apply Lebesgue differentiation theorem in such a bounded subdomain $U$ of  the domain $\{z\in\C:\,z\neq a,\,-\pi/2+\theta<\Arg(z-a)<\pi/2\}$.
%To verify the domination property, we
%
%
%we first deduce from the previous calculation the following lower bound:
%$$|(z-a)^2+(\tau_\theta(s)-a)^2|\geq \Im m((\tau_\theta(s)-a)^2).$$
%As a consequence, since
%$$\frac{\partial}{\partial z} h(z-a,a-\tau_\theta(s))=\frac{i k}{2R}H^{(1)}_1( k R)-\frac{i k^2(z-a)^2}{2R^2}H^{(1)}_2( k R) $$
%with $R=R(z-a,\tau_\theta(s)-a)$,
%using the asymptotic behavior (\ref{eq:HankAsym}) and (\ref{eq-lowerboundR}),
%and the fact that
%$$\Im m((\tau_\theta(s)-a)^2)\sim \sin(2\theta)|s|^2$$ for large $s$,
%we get the following estimate for a constant $C$ independent of $s\in\R$ and $z\in U$:
%$$\left|\frac{\partial}{\partial z} h(z,a,\tau_\theta(s))\right|\leq C\frac{\re^{- k \sin\theta|s|}}{|s|^{3/2}}. $$
%$$\left|\frac{\partial}{\partial z} h(z-a,a-\tau_\theta(s))\right|\leq C\re^{- k \sin\theta|s|}. $$
%This ensures  the domination  by an integrable function independent of $z$, and concludes the proof.
%	\end{proof}
Noting that, with $G_\theta$ as defined in the above lemma, $-x_2^1\in G_\theta$ for $x_2^1<-a$ and $-\tau_\theta(t)\in G_\theta$ for $t<-a$, we see that we have justified the analytic continuation of \eqref{eq:HSMM_step1} from $x_2^1<-a$ to the path $x_2^1=\tau_\theta(t)$, $t<-a$. Thus we obtain finally the new complex-scaled compatibility relation
		\[
			\varphi^1_\theta(t)=\int_{\R}h(-\tau_\theta(t)-a,a-\tau_\theta(s))\,\varphi^0_\theta(s)\,\tau'_\theta(s)\, \rd s,\quad t<-a.
		\]
		
By applying  similar reasoning, {and noting that $\tau_\theta(-t)=-\tau_\theta(t)$,} we get 8 equations linking the four complex-scaled traces (cf.\ \eqref{eq:syst_comp}}), namely
	\begin{equation}\label{eq:syst_comp_cmplx}
		\forall j\in \llbracket0,3\rrbracket,\quad\begin{array}{|lcl}
			\varphi^j_\theta(t)=S\,D_\theta\,\varphi^{j-1}_\theta(t),\quad t<-a,\\
			\varphi^j_\theta(t)=D_\theta\,S\,\varphi^{j+1}_\theta(t),\quad t>a,
			\end{array}
	\end{equation}
	where we have set $\varphi^{-1}_\theta:=\varphi^3_\theta$ and $\varphi^{4}_\theta:=\varphi^0_\theta$. In this system the operator $S$ is defined as in \eqref{eq:sym}, and, where $h$ is as given in \eqref{eq:hpr_kernel}, $D_\theta$ is defined, for $0<\theta<\pi/2$, by
	\begin{equation}
		\label{eq:opD_comp_exprG_cmplx}
		D_\theta\psi (t):=\int_{\R}h(\tau_\theta(t)-a,a-\tau_\theta(s))\,\psi(s)\,\tau'_\theta(s)\, \rd s,\quad t>a,\;\;\psi\in L^2(\R),
	\end{equation}
	 {and, similarly to \eqref{eq:DdefExt}, by}
	\begin{equation*}
	%	\label{eq:opD_comp_exprG_cmplx2}
		D_\theta\psi (t):=0,\quad t\leq a,\;\;\psi\in L^2(\R).
	\end{equation*}
 The system \eqref{eq:syst_comp_cmplx} has to be completed with the Dirichlet boundary condition
	\begin{equation}
		\label{eq:syst_comp_BC_cmplx}
		\varphi^j_\theta(t)=g|_{\Sigma_a^j}(t),\quad -a<t<a,\quad j\in \llbracket0,3\rrbracket.
	\end{equation}
	
As in Section \ref{section-HSMcomplexfreq}, these equations can be formulated as a single operator equation. Introducing
	\[
		\Phi_\theta:=\{\varphi^0_\theta,\varphi^1_\theta,\varphi^2_\theta,\varphi^3_\theta\}\in (L^2(\R))^4,
	\]
	and recalling the definition \eqref{eq:Phig} of $\Phi_g$, the systems of equations \eqref{eq:syst_comp_cmplx} and \eqref{eq:syst_comp_BC_cmplx} can be rewritten as
	\begin{equation}\label{eq:HSMM_matrice_cmplx}
		\begin{array}{c}
			\text{Find}\; \Phi_\theta\in (L^2(\R))^4\;\text{such that}\; \Phi_\theta-\Phi_g\in  (L^2_0(\R))^4\;\text{and}\\[5pt]
		(\mathbb{I}-\mathbb{D}_\theta)(\Phi_\theta-\Phi_g)=\mathbb{D}_\theta\,\Phi_g,
			\end{array}
	\end{equation}
	where $\mathbb{D}_\theta$ is obtained by replacing $D$ by $D_\theta$ in \eqref{eq:def_Dmat}, i.e.
	\begin{equation}\label{eq:def_Dmat_cmplx}
		\mathbb{D}_\theta:=\left[\begin{matrix}
		 0&\,{D}_\theta\,S&0&\,S\,{D}_\theta\\
		 \,S\,{D}_\theta&0&\,{D}_\theta\,S&0\\
		 0&S\,{D}_\theta&0&{D}_\theta\,S\\
		 {D}_\theta\,S&0&\,S\,{D}_\theta&0
		 \end{matrix}\right].
	\end{equation}
	
	Our first main result is that, as proved for $\mathbb{D}$ in the dissipative case (Theorem \ref{th:HSMM_complexe}), the operator $\mathbb{I}-\mathbb{D}_\theta$, as an operator on $ (L^2_0(\R))^4$, is Fredholm of index zero, indeed (importantly for numerical analysis of Galerkin methods) is a compact perturbation of a coercive operator.
	\begin{theorem}\label{th:HSMM_cmplxed}
For $0<\theta<\pi/2$:
		\begin{itemize}
			\item[(i)]$\mathbb{D}_\theta$ is a continuous operator on $(L^2(\R))^4$, and  $\mathbb{D}_\theta((L^2(\R))^4)\subset (L^2_0(\R))^4$;
			\item[(ii)]as an operator on $ (L^2_0(\R))^4$, $\mathbb{D}_\theta$ is the sum of an operator of norm $\leq 1/\sqrt{2}$ and a compact operator.
		\end{itemize}
	\end{theorem}

	To show well-posedness of \eqref{eq:HSMM_matrice_cmplx}, it remains to prove a uniqueness result which is the subject of section \ref{sec:unique} (Theorem \ref{thm:uniqueHSM}). At the end of this section we will combine the above theorem with Theorem \ref{thm:uniqueHSM} to write down a result expressing this well-posedness and the equivalence of  \eqref{eq:HSMM_matrice_cmplx} with the original scattering problem \eqref{pb:probleme_Dir_real}-\eqref{eq:Sommerfeld}.

The proof of Theorem \ref{th:HSMM_cmplxed}, which we defer to Appendix \ref{Appendix_Matrices}, mirrors the proof of Theorem \ref{th:HSMM_complexe}, once we establish properties of the operator $D_\theta$ to mirror those proved for $D$ in Proposition \ref{prop:D_complexe}. Establishing these properties of $D_\theta$, in Propositions \ref{prop:D_realk_cmplx_plus}, \ref{prop:D_realk_cmplx_moins}, and \ref{prop:D_realk_cmplx_tot}, is the focus of most of the rest of this section: these propositions
%  As for Theorem \ref{th:HSMM_complexe}, Theorem \ref{th:HSMM_cmplxed} is a direct consequence of the properties of the operator $D_\theta$, that are given .
  %In contrast to the single PropositMore precisely,
  give the properties of $D_\theta$ when it acts on functions whose support is, respectively, in $(a,+\infty)$, $(-\infty,-a)$,  and the whole of $\R$. (This splitting is necessary because of the piecewise definition of the complex-scaling function $\tau_\theta$.)

  Let us point out a useful fact (see the proof of the following proposition): when $D_\theta$ acts on functions whose support is in $(a,+\infty)$, it is equal %, up to the factor $\re^{\ri\theta}$,
  to the operator $D$ defined in \eqref{eq:opD_comp_def} for the dissipative case with wavenumber $ k \re^{\ri\theta}$.
	\begin{proposition}\label{prop:D_realk_cmplx_plus}
			Suppose that $\theta\in (0,\pi/2)$. For all $\psi \in L^2(a,+\infty)$ we have
				\begin{equation}
					\label{eq:opD_comp_exprG_cmplx_plus}
					D_\theta\psi (t)=\frac{\ri  (t-a) k \re^{\ri\theta}}{2}\int_{a}^{+\infty}\frac{H_1^{(1)}( k \re^{\ri\theta}\,R(t-a,s-a))}{R(t-a,s-a)}\,\psi(s)\, \rd s,\quad t>a,
				\end{equation}
				with $R$ defined in \eqref{eq:R(s,t)}.
			{As a consequence,} as an operator on $ L^2(a,+\infty)$, ${D}_\theta$ is the sum of an operator of norm $\leq 1/\sqrt{2}$ and a compact operator.
	\end{proposition}
	\begin{proof}
		Using the definition \eqref{eq:opD_comp_exprG_cmplx} of $D_\theta$ and the expression \eqref{eq:hpr_kernel} for the kernel $h$, we easily see \eqref{eq:opD_comp_exprG_cmplx_plus}. This implies that $D_\theta$, when it acts on functions whose support is in $(a,+\infty)$, is exactly the operator $D$ defined in \eqref{eq:opD_comp_def} for the dissipative case if we set the wavenumber in the dissipative case to be $ k \re^{\ri\theta}$. The result is therefore a direct consequence of item (ii) of Proposition \ref{prop:D_complexe}.
		\end{proof}
		\begin{proposition}\label{prop:D_realk_cmplx_moins}
				${D}_\theta$ is compact as an operator from $ L^2(-\infty,-a)$ to  $ L^2(a,+\infty)$, for $\theta\in(0,\pi/2)$.
		\end{proposition}
		\begin{proof}
			Using the definition \eqref{eq:opD_comp_exprG_cmplx} of $D_\theta$ and the expression \eqref{eq:hpr_kernel} for the kernel $h$, we can show easily that, for all $\psi \in L^2(-\infty,-a)$, we have
							\[
								D_\theta\psi (t)=\frac{\ri  k (t-a)\re^{2\ri\theta}}{2}\int_{-\infty}^{-a}\frac{H_1^{(1)}( k R_\theta(t,s))}{R_\theta(t,s)}\,\psi(s)\, \rd s,\quad t>a,
							\]
							with $R_\theta(t,s):=R(\re^{\ri\theta}(t-a),2a-(s+a)\re^{\ri\theta})$ and $R$ defined in \eqref{eq:R(s,t)}.
							%$R_\theta(t,s)=\big[\re^{2\ri\theta}(t-a)^2+(2a-(s+a)\re^{\ri\theta})^2)\big]^{1/2}$.
							We prove the compactness of $D_\theta$, as we prove the compactness of $D$ in (iii) of Proposition \ref{prop:D_complexe}, by showing that $D_\theta$ is a Hilbert-Schmidt operator from $ L^2(-\infty,-a)$ to  $ L^2(a,+\infty)$. By a simple change of variable $t\mapsto t-a$ and $s\mapsto-(s+a)$, it suffices to show that
							\begin{equation}\label{eq:Dmoins_HS}
								\int_{0}^{+\infty}\int_{0}^{+\infty}|K(t,s)|^2\,\rd t\, \rd s<+\infty,
							\end{equation}
							where
							\[
								K(t,s):=t\,\frac{H_1^{(1)}( k \widetilde{R}_\theta(t,s))}{\widetilde{R}_\theta(t,s)}\quad\text{and}\quad \widetilde{R}_\theta(t,s):=R(\re^{\ri\theta}t,2a+\re^{\ri\theta}s).
							\]
							Where $\R^+:= [0,+\infty)$, $K$ is continuous on $\R^+\times\R^+$, since $\Re(\tilde{R}_\theta)>0$ on $\R^+\times\R^+$. % See for instance Lemma \ref{lem:Rbound}.
We want to use now the bound on the Hankel function given in \eqref{eq:Hbound1} which implies that, for some $c>0$,
\begin{equation} \label{eq:Hbound1_FF}
\left|\frac{ H_1^{(1)}(z)}{z}\right| \leq c|z|^{-3/2}\re^{-\Im(z)}, \quad \Re(z)>0,\;|z|\geq 1.
\end{equation}
$\widetilde{R}_\theta$ can be rewritten as
%\begin{equation}\label{eq:Rthetamoins}
$\widetilde{R}_\theta(t,s)=\hat{R}(\hat{z},z)$, where  $\hat{z}:={4ase^{\ri\theta}}+{4a^2}$, $z:=e^{\ri\theta}(t^2+s^2)^{1/2}$, and
%\end{equation}
$\hat{R}$ is defined in \eqref{eq-defRandRhat}. Thus we can use Lemma \ref{lem:Im(Rhat-z)} and deduce that, for some $C>0$ and all $s,t\in \R^+$ with $s^2+t^2$ sufficiently large, it holds that
        $$|\widetilde{R}_\theta(t,s)-\re^{\ri\theta}(t^2+s^2)^{1/2}|\leq C\frac{1+|s|}{(s^2+t^2)^{1/2}}$$
		and
		\[
			\Im (\widetilde{R}_\theta(t,s))\geq \sin(\theta)(t^2+s^2)^{1/2}-C\frac{1+|s|}{(s^2+t^2)^{1/2}}.
		\]
        Thus, for {some $C'>0$ and } all $s^2+t^2$ large enough, $|\widetilde{R}_\theta(t,s)|\geq \frac{1}{2}(t^2+s^2)^{1/2}$ and $\Im (\widetilde{R}_\theta(t,s))\geq \sin(\theta)(t^2+s^2)^{1/2}-{C'}$, so that, by \eqref{eq:Hbound1_FF},
\[
	|K(t,s)|=\mathcal{O}\left(\frac{\re^{- k\sin(\theta) \sqrt{t^2+s^2}}}{(t^2+s^2)^{1/4}}\right)\quad\text{as}\;\sqrt{t^2+s^2}\rightarrow+\infty,
\]
uniformly in $t$ and $s$. Thus \eqref{eq:Dmoins_HS} is clear.
			\end{proof}
			\begin{proposition}\label{prop:D_realk_cmplx_tot}
					 For $0<\theta<\pi/2$, ${D}_\theta$ is a continuous operator from $L^2(\R)$ to  $ L^2(a,+\infty)$.
			\end{proposition}
			\begin{proof}
From  Propositions \ref{prop:D_realk_cmplx_plus} and \ref{prop:D_realk_cmplx_moins}, it suffices to show that  ${D}_\theta$ is a continuous operator from $L^2(-a,a)$ to  $ L^2(a,+\infty)$. But this is immediate from the continuity of $\cD^j_\theta$ from $L^2(\Sigma_a)$ to $L^2(a,+\infty)$, established in Proposition \ref{prop:SDbounded}.
				\end{proof}

We finish this section with the promised statement of well-posedness of \eqref{eq:HSMM_matrice_cmplx}, and of its equivalence with the original scattering problem \eqref{pb:probleme_Dir_real}-\eqref{eq:Sommerfeld}.
\begin{theorem} \label{th:mainWP} For every $\theta \in (0,\pi/2)$, the operator $\mathbb{I}-\mathbb{D}_\theta$ is invertible on $(L^2_0(\R))^4$. Thus, for every $g\in L^2(\Sigma_a)$, \eqref{eq:HSMM_matrice_cmplx} has exactly one solution $\Phi_\theta\in (L^2(\R))^4$ such that $\Phi_\theta-\Phi_g\in  (L^2_0(\R))^4$. Moreover, for some constant $c>0$ {depending on $\theta$},
\begin{equation} \label{eq:stability}
\|\Phi_\theta\|_{(L^2(\R))^4} \leq c\|\Phi_g\|_{(L^2(\R))^4}=c\|g\|_{L^2(\Sigma_a)},
\end{equation}
for every $g\in L^2(\Sigma_a)$. Further, if $g\in H^{1/2}(\Sigma_a)$ and $\Phi_\theta=\{\varphi^0_\theta,\varphi^1_\theta,\varphi^2_\theta,\varphi^3_\theta\}\in (L^2(\R))^4$ is the solution of \eqref{eq:HSMM_matrice_cmplx}, then, for $j\in \llbracket0,3\rrbracket$:
\begin{itemize}
\item[(i)]  $\varphi_\theta^j$ satisfies \eqref{eq:definition_phi_theta}, i.e.\ $\varphi_\theta^j$ is the analytic continuation to the path $\tau_\theta$ of the restriction to $\Sigma^j$ of the solution $u$ of \eqref{pb:probleme_Dir_real}-\eqref{eq:Sommerfeld};
    \item[(ii)] the solution $u$ of \eqref{pb:probleme_Dir_real}-\eqref{eq:Sommerfeld} is given in terms of $\varphi_\theta^j$ in $\Omega_\theta^j$ by \eqref{eq-repinOmegatheta};
        \item[(iii)] for some constant $C>0$ that depends only on $a$, $k$, and $\theta$, $\varphi_\theta^j(s)$ satisfies the bounds of Proposition \ref{cor:varphiL2}, for $|s|>a$.
\end{itemize}
\end{theorem}
\begin{proof} Theorem \ref{th:HSMM_cmplxed} implies that, as an operator on $(L^2_0(\R))^4$, $\mathbb{I}-\mathbb{D}_\theta$ is Fredholm of index zero, and Theorem \ref{the-uniqueness} implies that it is injective, so that $\mathbb{I}-\mathbb{D}_\theta$ is invertible with a bounded inverse. This implies, since $\mathbb{D}_\theta$ is a bounded operator from $(L^2(\R))^4$ to  $(L^2_0(\R))^4$ by Theorem \ref{th:HSMM_cmplxed}(i), that \eqref{eq:HSMM_matrice_cmplx} has exactly one solution, and this solution satisfies the bound \eqref{eq:stability}.

In the case that $g\in H^{1/2}(\Sigma_a)$, that (i) holds follows from the derivation of \eqref{eq:HSMM_matrice_cmplx} from \eqref{pb:probleme_Dir_real}-\eqref{eq:Sommerfeld} in Section \ref{sec:Derivation}, and since \eqref{eq:HSMM_matrice_cmplx} has only one solution; that (ii) holds follows from (i) and Theorem \ref{th:deformed_representation}; that (iii) holds follows from (i) and Proposition \ref{cor:varphiL2}.
\end{proof}

  \section{Reconstruction of the solution and far-field formula} \label{sec:recon}
    Suppose that we have computed the solution $\Phi_\theta$ to (\ref{eq:HSMM_matrice_cmplx}). Then the solution $u$ of the problem (\ref{pb:probleme_Dir_real}-\ref{eq:Sommerfeld})  can be recovered a posteriori through the representation formulas (\ref{eq-repinOmegatheta}). More precisely, as we have observed in Theorem \ref{th:mainWP}(ii), it can be reconstructed in the union for $j\in \llbracket0,3\rrbracket$ of the domains $\Omega^j_\theta$ defined by  (\ref{eq:Omegaj_theta}). Let us point out that, if $\theta<\pi/4$, the  union of the $\Omega^j_\theta$ covers the whole domain $\Omega$, so that the whole solution $u$ can  be reconstructed a posteriori.

      It is well known (e.g., \cite[Lemma 2.5]{ChGrLaSp:11}) that the solution of (\ref{pb:probleme_Dir_real}-\ref{eq:Sommerfeld}) satisfies
   \begin{equation}\label{eq-farfield}
     u(\bsx) = \frac{\re^{\ri k r}}{r^{1/2}}\left(F(\widehat\bsx) + \mathcal{O}(r^{-1})\right), \quad \mbox{as} \quad r\to +\infty,
  \end{equation}
     uniformly in $\widehat\bsx := \bsx/r$, where $F\in C^\infty(S^1)$, with $S^1$ the unit circle, is the \emph{far-field pattern}. By analogy with classical boundary integral methods, one can wonder if this far-field pattern can also be recovered  from properties of  %as a function of
     the $\varphi_\theta^j$. A partial answer will be given in this section, by deriving far-field formulas in the four directions orthogonal to the edges of the square $\Omega_a$. The proof of this result requires first that we establish some properties of the solution $u$ that can be deduced from the representation formulas (\ref{eq-repinOmegatheta}). The far-field behaviour that we will establish, indeed all of the results of this section, will be ingredients in the
  proof of uniqueness  for problem (\ref{eq:HSMM_matrice_cmplx}) that will be the focus of the next section.
	
 From (\ref{eq-repinOmegatheta}), let us consider the representation formula for any $\psi\in L^2(\R)$
 	\begin{equation}\label{eq:expr_Ujtheta}
 		 U_\theta(\psi)({\bsx}^0) := \int_{\R}h(x_1^0-a,x_2^0-\tau_\theta(s))\,\psi(s)\,\tau'_\theta(s)\, \rd s,\quad {\bsx}^0\in \Omega^0_\theta,
 	\end{equation}
 	and the associated integral operator defined by
 	\begin{equation}\label{eq:opDthetainter}
 		\widetilde{D}_\theta\psi(t):=U_\theta(\psi)(t,a),\quad t>a.
 	\end{equation}
 	Let us note that in Lemma \ref{lem:analyticity_utheta} we have shown that, for any data $\psi\in L^2(\R)$, the function $t\mapsto \widetilde{D}_\theta\psi(t)$ can be continued analytically  from $(a,+\infty)$ into the domain of the complex plane $\{z\in\C:\,z\neq a,\,-\pi/2+\theta<\Arg(z-a)<\pi/2\}$ just replacing $t$ by $z$ in \eqref{eq:opDthetainter}.
 	The following proposition bounds $\widetilde{D}_\theta\psi(z)$, in particular near $a$ and for large values of $z$.
 	% First, due to Lemma \ref{lem-analytic}, we know that $D_\theta\psi$ is a  $\mathcal{C}^\infty$ function of $t$ in $(a,+\infty)$. Let us now derive some estimates of the behavior of $D_\theta\psi(t)$ respectively when $t\rightarrow +\infty$ and when $t\rightarrow a$. These results are very similar to the estimates summarized in Corollary \ref{cor:varphiL2}.
 	%This will be a component in the proof of uniqueness in the next section.

%Thinking about the above I realise that our Lemma 4.1 can be simplified and sharpened (and this is needed to see that D:X\to X).
%Let F(z) := \widetilde D_\theta \psi(z)(1+ (z-a)^{-1/2})^{-1} exp(-ikz), which is analytic in Re(z)>a with 0\leq Arg(z-a) \leq \theta. Then, by Lemma 4.1,
%|F(z)| \leq C\|\psi\|_L^2 B(z), where B(z) := exp(k Im(z)(1- \cos(\theta - \Arg(z-a)))). Now B(z)=1 when \Arg(z-a) = 0 or \theta, and

%|B(z)| \leq exp(k|z|), for all 0\leq Arg(z-a) \leq \theta.

%Thus, by a standard Phragmen-Lindelof argument (e.g.  Conway, Functions of One Complex Variable I, Chapter VI, Cor. 4.2),
%|F(z)| \leq C\|\psi\|_L^2 for all 0\leq Arg(z-a) \leq \theta,
%i.e. the bound in Lemma 4.1 can be simplified/sharpened to
%|\widetilde D_\theta \psi(z)| \leq C(|z-a|^{-1/2}+1) exp(-k Im(z)) \|\psi\|_L^2.
%I think we should make the small edits needed to do this sharpening/simplification. It’s just an extra Phragmen-Lindelof step at the end of the proof of Lemma 4.1.

 \begin{lem}
         \label{lem-Dtheta(t)infty}
         There exists a constant $C>0$, depending only on $\theta$, $a$ and $k$, such that, for all  $\psi\in L^2(\R)$ and $z\in\C$ with $\Re(z)>a$ and $\Arg(z-a)\in [0,\theta]$,
         $$|\widetilde{D}_\theta\psi(z)|\leq C(|z-a|^{-1/2}+1)\,{\exp(-k \Im(z))}\, \|\psi\|_{L^2(\R)}.$$
     \end{lem}
 \begin{proof} Throughout this proof $C$ will denote any positive constant, depending  only on $a$, $k$, and $\theta$, not necessarily the same at each occurrence.

     Using the expressions (\ref{eq:expr_Ujtheta}-\ref{eq:opDthetainter}), the definition \eqref{eq:hpr_kernel} of the kernel $h$, and the definition of the complex-scaling function (\ref{eq:parametrage_t_theta}), we see that  $ |\widetilde{D}_\theta\psi(z)|\leq \mathcal{I}(z-a)$, where
 \begin{equation*}% \label{eq:basicbound}
 \mathcal{I}(w) := \frac{k|w|}{2}\int_{\R} \left|\frac{H_1^{(1)}( k \,R(w,\tau_\theta(s)-a))}{R(w,\tau_\theta(s)-a)}\right| \, |\psi(s)|\,ds, \quad \Re(w)>0, \;\; \Arg(w)\in [0,\theta],
 \end{equation*}
 with $R$ defined in \eqref{eq:R(s,t)}. The bound on the Hankel function given in \eqref{eq:Hbound1}  implies that
 \begin{equation} \nonumber %\label{eq:Hbound1_bis}
 \left|\frac{ H_1^{(1)}(z)}{z}\right| \leq C\left( |z|^{-2}+|z|^{-3/2}\right)\re^{-\Im z}, \quad \Re(z)>0.
 \end{equation}
  Further, Lemma \ref{lem-Rbound-z-tau} (applied with $\theta_0=\theta$) gives, for $\Arg(w)\in [0,\theta]$, that
 \[
 	|R(w,\tau_\theta(s)-a)|^2\geq \cos(\theta)(|w|^2+|\tau_\theta(s)-a|^2)
 \]
 and
 \[
 	\Im R(w,\tau_\theta(s)-a)\geq \cos(\theta-\Arg(w))\Im(w)-C.
 \]
 Letting $t:=|w|$ and $\gamma:=\Arg(w)$, we deduce from the above bounds that
 \begin{eqnarray*}
 |\mathcal{I}(w)| \,\re^{k\cos(\theta-\gamma)\Im(w)}& \leq & C\int_\R\bigg(\frac{t}{t^2+|\tau_\theta(s)-a|^2}+
 \frac{t}{\left(t^2+|\tau_\theta(s)-a|^2\right)^{3/4}}\bigg)\,|\psi(s)|\,\rd s.
 \end{eqnarray*}
Applying the Cauchy-Schwarz inequality, and noticing that  $|\tau_\theta(s)-a|^2=(s-a)^2$, for $s\geq -a$, while $|\tau_\theta(s)-a|^2=(s+a)^2 + 4a^2(1-\cos(\theta))-4as\cos(\theta)\geq(s+a)^2$, for $s< -a$, yields
 	\begin{equation} \nonumber
 		|\mathcal{I}(w)|\, \re^{k\cos(\theta-\gamma)\Im(w)}\leq C \left[\int_{0}^{+\infty}\frac{t^2}{(t^2+s^2)^2}\, \rd s +\int_{0}^{+\infty}\frac{t^2}{(t^2+s^2)^{3/2}}\,\rd s\right]^{1/2}\, \|\psi\|_{L^2(\R)}.
 	\end{equation}
We see, by substituting $s=tp$, that the first and second integrals on the right hand side of this last inequality are $\leq Ct^{-1}$ and $\leq C$, respectively. Thus we have shown that
 $$|\widetilde{D}_\theta\psi(z)|\leq C(|z-a|^{-1/2}+1)\,\exp(-k\cos(\theta-\Arg(z-a))\Im(z-a))\, \|\psi\|_{L^2(\R)},$$
 for all  $\psi\in L^2(\R)$ and $z\in\C$ with $\Re(z)>a$ and $\Arg(z-a)\in [0,\theta]$. %This gives immediately the required bound in the cases that $\Arg(z-a)\in \{0,\theta\}$.
Now, defining
 $$
 \mathcal{F}(z) := (1+ (z-a)^{-1/2})^{-1}\, \exp(-\ri kz)\,\widetilde D_\theta \psi(z),
 $$
 this last bound implies that, for $\Re(z)>a$ with $0\leq \Arg(z-a) \leq \theta$, $\mathcal{F}$ is analytic and
$|\mathcal{F}(z)| \leq C\|\psi\|_{L^2(\R)}\, B(z)$, where $B(z) := \exp(k \Im(z)(1- \cos(\theta - \Arg(z-a))))$. Now $B(z)=1$ when $\Arg(z-a) = 0$ or $\theta$, and $|B(z)| \leq \exp(k|z|)$, for all $z$ with $\Re(z)>a$ and $0\leq \Arg(z-a) \leq \theta$. Thus, by a standard Phr\'agmen-Lindel\"of principle (e.g., \cite[Chapter VI, Cor. 4.2]{Con78}),
$|\mathcal{F}(z)| \leq C\|\psi\|_{L^2(\R)}$, for all $z$ with $\Re(z)>a$ and $0\leq \Arg(z-a) \leq \theta$, %for all 0\leq Arg(z-a) \leq \theta,
%i.e. the bound in Lemma 4.1 can be simplified/sharpened to
and the required bound on $|\widetilde{D}_\theta\psi(z)|$ follows.     \end{proof}% %and finally
 	%\[
 	%	\mathcal{I}(te^{\ri\gamma})\leq C[t^{-1/2}+1]\re^{-k\cos(\theta-\gamma)t\sin\gamma}
 	%\]
    %         where we have used that
 %            $$\int_{0}^\infty\frac{t^2}{(t^2+s^2)^2}\,\rd s=\frac{1}{t}\int_{0}^\infty\frac{1}{(1+s^2)^2}\,\rd s, \;\; \int_{0}^\infty\frac{t^2}{(t^2+s^2)^{3/2}}\, \rd s=\int_{0}^\infty\frac{1}{(1+s^2)^{3/2}}\, \rd s.$$
%     \end{proof}

 	% In the previous section, we proved that the operator $D_\theta$ is continuous from $L^2(\R)$ to  $ L^2(a,+\infty)$. Our first objective is to get more precise properties of the function $D_\theta\psi$ for some arbitrary function $\psi\in L^2(\R)$.

 %\sce%dnote{Deleted, because I don't think we need it at this point: One consequence of the analyticity  of the mapping  $z\mapsto \widetilde{D}_\theta\psi(z)$ is that if $\{\varphi^0_\theta,\varphi^1_\theta, \varphi^2_\theta,\varphi^3_\theta\}$ is a solution of \eqref{eq:HSMM_matrice_cmplx} then the four functions are $C^\infty$ in any interval which does not include $t=\pm a$, they are in $L^1$ and they are exponentially decaying at infinity.}

Defining $\phi(t):= \widetilde{D}_\theta\psi(t)$, for $t>a$, Lemma \ref{lem-Dtheta(t)infty} implies that  $\phi\in L^1(a,b)$, for every $b>a$, if $\psi\in L^2(\R)$. In the uniqueness proof in the next section we will need also the following stronger result.
\begin{lem}
         \label{lem-Dtheta(t)infty2}
         If $\psi\in L^2(\R)$ and $\phi(t):= \widetilde{D}_\theta\psi(t)$, for $t>a$, then  $\phi\in L^2(a,b)$, for every $b>a$.
     \end{lem}
 \begin{proof} %Lemma \ref{lem-Dtheta(t)infty} implies that, for every $b>a$, $\phi\in L^\infty(b,\infty)$. Further,
 Arguing as in the proof of Lemma \ref{lem-Dtheta(t)infty}, we see that %$|\phi(t)| \leq \mathcal{I}(t-a)$, for $t>a$, where $\mathcal{I}$ (defined by \eqref{eq:basicbound}) satisfies
 \begin{eqnarray*}
 |\phi(t)| & \leq & C\left(\int_\R \frac{(t-a)\,|\psi(s)|}{(t-a)^2+|\tau_\theta(s)-a|^2}\,\rd s +
\int_\R \frac{(t-a)\, |\psi(s)|}{\left((t-a)^2+|\tau_\theta(s)-a|^2\right)^{3/4}}\,\rd s\right),
 \end{eqnarray*}
 for $t>a$. Further, arguing as at the end of the proof of Lemma \ref{lem-Dtheta(t)infty}, using that $|\tau_\theta(s)-a|^2=(s-a)^2$, for $s\geq -a$, while $|\tau_\theta(s)-a|^2\geq(s+a)^2$, for $s< -a$, we see that the second integral in the above sum is bounded on $(a,+\infty)$, and so is in $L^2(a,b)$, for every $b>a$, while the first integral is
 \begin{eqnarray*}% \label{eq:RHSstatic}
 \leq  \int_{-a}^{+\infty} \frac{t-a}{(t-a)^2+(s-a)^2}|\psi(s)|\,\rd s + \int_{a}^{+\infty} \frac{t-a}{(t-a)^2+(s-a)^2}|\psi(-s)|\,\rd s.
 \end{eqnarray*}
 The right hand side of this last inequality is in $L^2(a,+\infty)$ by Remark \ref{rem:D0bounded}.
 %But, arguing as at the end of the proof of Proposition \ref{prop:SDbounded}, relating the above integrals to the static ($k=0$) version of the operator $D$ of Proposition \ref{prop:D_complexe}, we see that the right hand side of \eqref{eq:RHSstatic} is in $L^2(\textcolor{red}{0},\infty)$.
 \end{proof}

Let us remark that by definition \eqref{eq:opD_comp_exprG_cmplx} we have, for any $\psi\in L^2(\R)$,
 \begin{equation}
 	\label{eq:Dtilde-D}
 		D_\theta\psi(t)=\tilde{D}_\theta\psi(\tau_\theta(t))=U_\theta(\psi)(\tau_\theta(t),a), \quad t>a.
 \end{equation}
 	We deduce thus from Lemma \ref{lem-Dtheta(t)infty} the following result.
     \begin{cor}
         \label{cor-Dtheta(t)infty}
         There exists a constant $C>0$, depending only on $\theta$, $a$ and $k$, such that, for all  $\psi\in L^2(\R)$,
         $$|D_\theta\psi(t)|\leq C[(t-a)^{-1/2}+1]\,\re^{-k\sin (\theta)t}\|\psi\|_{L^2(\R)},\quad t>a.$$
     \end{cor}

 	% Besides, we know thanks to Theorem \ref{th:deformed_representation} that the outgoing solution of \eqref{pb:probleme_Dir_real} can be represented thanks to the complexifed traces $\varphi^j_\theta$ by
 % 	\begin{equation}\label{eq:u_Omega_j_theta}
 % 		 u({\bsx}^j) = \int_{\R}h(x_1^j-a,x_2^j-\tau_\theta(s))\,\varphi^j_\theta(s)\,\tau'_\theta(s)\, \rd s,\quad {\bsx}^j\in \Omega^j_\theta,
 % 	\end{equation}
 % 	where $\Omega^j_\theta$, a domain which is strictly included in $\Omega^j$, is defined in \eqref{eq:Omegaj_theta}. These representations are interesting in itself since they enable to construct the outgoing solution in several subdomains of the domain $\Omega$. Let us remark also that if $\theta>\pi/4$, this enables to construct the solution in the whole domain $\Omega$ since the $\Omega^j_\theta$ overlap. A question is: can one also recover the far-field from the knowledge of the  $\varphi_j^\theta$ ?The next proposition gives a partial answer :  it gives the far-field amplitudes in the directions orthogonal to the edges of the square $\Omega_a$.
 % This result will be also useful for the proof of uniqueness in the next section.
 Now we are able to prove the main result of this section, which provides far field formulas in the direction orthogonal to the edges of the square $\Omega^a$.
 	\begin{prop}\label{prop:FF_sol}
 		Let $\psi\in L^2(\R)$ be such that $(t^2+1)\psi(t)\in L^1(\R)$. Then the function $ U_\theta(\psi)$ defined by \eqref{eq:expr_Ujtheta}
 		has the following behaviour at infinity:
 		\[
 			U_\theta(\psi)(x_1,x_2) = C_\infty \frac{\re^{\ri k x_1}}{\sqrt{x_1}}\left(1+\mathcal{O}\left(x_1^{-1}\right)\right),\quad\text{as}\; x_1\rightarrow +\infty,\;\text{for each}\;x_2\in[-a,a]\;\text{fixed,}
 		\]
 		where
 		\[
 			C_\infty :=\sqrt{\frac{k}{\pi}}\frac{1-\ri}2\int_\R\psi(s)\tau'_\theta(s)\,ds.
 		\]
 		\end{prop}
 		\begin{proof}
 			Throughout this proof, $C$ will denote a positive constant, depending  only on $a$, $k$, and $\theta$, not necessarily the same at each occurrence.
 			
 			For any fixed $x_2\in[-a,a]$, let $\phi_{x_2}(t):=U_\theta(\psi)(t+a,x_2)$, for $t>0$. By the definition \eqref{eq:hpr_kernel} of $h$,  we have
 			\[
 				\phi_{x_2}(t)=\frac{\ri k t}{2}\int_{\R}\frac{H_1^{(1)}(kR(t,x_2-\tau_\theta(s)))}{R(t,x_2-\tau_\theta(s))}\,\psi(s)\,\tau'_\theta(s)\, \rd s.
 			\]
 		In order to use \eqref{eq:HankAsym} that says that
 		\begin{equation}\label{eq:FF_Hankel}
 		H^{(1)}_1(z) = -\frac{1+\ri}{\sqrt{\pi}}\ \frac{\re^{\ri z }}{\sqrt{z}} \Big( 1 + \mathcal{O}\big( |z|^{-1}\big) \Big)
 		\quad\text{as } |z| \to +\infty,
 		\end{equation}
 		uniformly in $\Arg(z)$ for $|\Arg(z)|<\zeta$, for every $\zeta<\pi$, we write
 			\begin{equation}\label{eq:decomp_phib}
 		\phi_{x_2}(t)-\dsp \sqrt{\frac{k}{\pi}}\frac{1-\ri}2\frac{\re^{\ri kt }}{\sqrt{t}}\int_{\R}\psi(s)\,\tau'_\theta(s)\, \rd s=\phi_{x_2}^1(t)+\phi_{x_2}^2(t)
 		\end{equation}
 		where
 		$$\begin{array}{rl}
 		\phi_{x_2}^1(t):=&\dsp \sqrt{\frac{k}{\pi}}\frac{1-\ri}2t\int_{\R}\left[ \frac{\re^{\ri k R(t,x_2-\tau_\theta(s)) }}{[R(t,x_2-\tau_\theta(s))]^{3/2}}- \frac{\re^{\ri kt }}{t^{3/2}}\right]\,\psi(s)\,\tau'_\theta(s)\, \rd s,\\[10pt]
 		\phi_{x_2}^2(t):=& \dsp\frac{\ri k t}{2}\int_{\R}\bigg[\frac{H_1^{(1)}(kR(t,x_2-\tau_\theta(s)))}{R(t,x_2-\tau_\theta(s))}\\[10pt]
 		& \dsp \hspace*{13ex} + \hspace*{1ex} \frac{1+\ri}{\sqrt{k\pi}}\ \frac{\re^{\ri kR(t,x_2-\tau_\theta(s)) }}{[R(t,x_2-\tau_\theta(s))]^{3/2}}\bigg]\,\psi(s)\,\tau'_\theta(s)\, \rd s.
 		\end{array}$$
 						
 		{To estimate these quantities we note first that, if $t>0$ and $z\in Z_\theta:=\{z= r\re^{\ri \alpha}: 0\leq \alpha \leq \theta,\, r\geq 0\}$ (or $-z\in Z_\theta$), then $0\leq \Arg(R(t,z))\leq \pi/2$, $|\cos(\Arg(z))|\geq \cos(\theta)$,  and, by Lemma \ref{lem:Rbound},
 						\begin{equation}
 						\label{Rbound-farfield2}
 					|R(t,z)|\geq \cos(\theta)\, t.
 						\end{equation}
 Since, {from the definition (\ref{eq:parametrage_t_theta}),} $x_2-\tau_\theta(s)$ or $\tau_\theta(s)-x_2$ is in $Z_\theta$ for all $x_2\in [-a,a]$ and all $s\in \R$, these observations hold in particular if $z=x_2-\tau_\theta(s)$, so that
 %$\phi_{x_2}^2(t)$, we notice that for $x_2\in[-a,a]$, $|\cos(\Arg(x_2-\tau_\theta(s)))|\geq \cos(\theta)$ for all $s\in\R$. Therefore, by  Lemma \ref{lem:Rbound},
% 						\begin{equation}
% 						\label{Rbound-farfield}
% 					|R(t,x_2-\tau_\theta(s))|\geq \cos(\theta)\, t,\qquad t>0,\ s\in\R.
% 						\end{equation}
% 						Since \scednote{Changed $|\Arg(R(t,x_2-\tau_\theta(s)))|\leq \pi/2$. It is true and, I think, essential, that the argument is non-negative, so that $|\exp(\ri kR)|\leq 1$} ${0 \leq } \Arg(R(t,x_2-\tau_\theta(s)))\leq \pi/2$,
 \eqref{eq:FF_Hankel} applies and implies that}
 						\begin{equation}
 					\label{phi2}
 					|\phi_{x_2}^2(t)|\leq \frac{C}{t^{3/2}}\int_\R|\psi(s)|\,\rd s,
 						\end{equation}
{for all sufficiently large $t>0$. 		To estimate $\phi_{x_2}^1(t)$ we observe that, for $t>0$ and $z\in Z_\theta$,}
 		$$\left| \frac{\re^{\ri kR(t,z) }}{[R(t,z)]^{3/2}}- \frac{\re^{\ri kt }}{t^{3/2}}\right|\leq \sup_{w\in U(z)}\left|\frac{\partial }{\partial w}\left[ \frac{\re^{\ri kR(t,w) }}{[R(t,w)]^{3/2}}\right]\right|\,  |z|  $$
 				where $U(z):= \{w\in Z_\theta:|w|\leq |z|\}$ %\{ w\in \mathbb{C}, \; |w|\leq |z|,\; 0\leq \Arg(w)\leq \theta\}$
 and where we have used that $R(t,0)=t$. Using again (\ref{Rbound-farfield2}), this yields {that}
		$$\left| \frac{\re^{\ri kR(t,z) }}{[R(t,z)]^{3/2}}- \frac{\re^{\ri kt }}{t^{3/2}}\right|\leq  C\frac{|z|^2}{t^{5/2}},\qquad %0\leq \Arg(z)\leq \theta, \;\;
{z\in Z_\theta,\;\; t>0.} $$		
				 Since $|x_2-\tau_\theta(s)|\leq C(1+|s|)$, it {follows that}
 			\begin{equation}
 			\label{phi1}
 			|\phi_{x_2}^1(t)|\leq \frac{C}{t^{3/2}}\int_\R(1+|s|)^2|\psi(s)|\, \rd s, \quad {t>0.}
 			\end{equation}		
 Combining (\ref{eq:decomp_phib}), (\ref{phi2}) and (\ref{phi1}), we get that
 						\[
 							\left|\phi_{x_2}(t)- \sqrt{\frac{k}{\pi}}\frac{1-\ri}2\frac{\re^{\ri kt }}{\sqrt{t}}\int_{\R}\psi(s)\,\tau'_\theta(s)\, \rd s\right|\leq \frac{C}{t^{3/2}}\int_{\R}(1+s^2)|\psi(s)|\, \rd s,
 						\]
 {for all sufficiently large $t>0$, }					which ends the proof.
 						\end{proof}
  	
  This proposition proves in particular that
  	\begin{equation}
  	\label{eq-farfieldtheta}
  	F\left(\cos (j\pi/2),\sin (j\pi/2)\right)=\sqrt{\frac{k}{\pi}}\frac{1-\ri}2\int_\R\varphi^j_\theta(s)\tau'_\theta(s)\,\rd s, \quad j\in \llbracket0,3\rrbracket,
  	\end{equation}
  	where the far-field $F$ is defined by \eqref{eq-farfield}  and  $\varphi^0_\theta,\varphi^1_\theta,\varphi^2_\theta,\varphi^3_\theta$ are the complex-scaled traces of $u$, which are exponentially decaying at infinity thanks to {Proposition \ref{cor:varphiL2}, or thanks to \eqref{eq:syst_comp_cmplx} and Corollary \ref{cor-Dtheta(t)infty}}.

 	\section{Uniqueness} \label{sec:unique}
 	%Suppose that $\{\varphi^0_\theta,\varphi^1_\theta,\varphi^2_\theta,\varphi^3_\theta\}\in [L^2(\R)]^4$ satisfies the complex-scaled version of the HSM \eqref{eq:syst_comp_cmplx} with the condition
 	%\[
 	%\varphi^j_\theta(t)=0,\;-a<t<a,\quad \;j\in \llbracket 0,3\rrbracket.
 	%\]
 	%We want to show that $\varphi^j_\theta\equiv 0$ for all $j$.
	
 	In this section we will prove uniqueness of solution for the complex-scaled HSM method \eqref{eq:HSMM_matrice_cmplx}. {This result, important in its own right, is also key to the proof of uniqueness for the complex-scaled HSM method for more general configurations; see Proposition \ref{lem-uniqueness-generalcase} below.} Our proof depends on the following uniqueness result for  the standard HSM \eqref{eq:HSMM_complexe} that we prove, {using completely different arguments,} in \cite{bibid}. %This result, which holds for real $k$, requires that the solutions $\varphi^j$ of \eqref{eq:HSMM_complexe} satisfy the radiation condition thta

\begin{theorem} \label{thm:uniqueHSM} If $k>0$ and $\{\varphi^0,\varphi^1,\varphi^2,\varphi^3\}\in (L^2_{\mathrm{loc}}(\R))^4$  satisfies \eqref{eq:HSMM_complexe} with $g=0$, and satisfies the radiation condition that
 \begin{equation}
 \label{eq:RadCond}
 \varphi^j(t) = \left\{\begin{array}{cc} c_+^j\,\re^{\ri k |t|}|t|^{-1/2}\left(1 + \mathcal{O}(|t|^{-1})\right), & \mbox{as} \quad t\to +\infty,\\
 c_-^j\,\re^{\ri k |t|}|t|^{-1/2}\left(1 + \mathcal{O}(|t|^{-1})\right), & \mbox{as} \quad t\to -\infty,\end{array}\right.
 \end{equation}
 for some constants $c_\pm^j\in\C$ and every $j\in \llbracket0,3\rrbracket$,
 then $\varphi^j=0$, for $j\in \llbracket0,3\rrbracket$.
 \end{theorem}

To see how this result is relevant to uniqueness for the complex-scaled HSM method \eqref{eq:HSMM_matrice_cmplx}, recall that to derive \eqref{eq:syst_comp_cmplx} and \eqref{eq:syst_comp_BC_cmplx} we started from \eqref{eq:HSMM_complexe}, satisfied for real $k$ by $\varphi^j$, the trace on $\Sigma^j$ of the solution $u$ of \eqref{pb:probleme_Dir_real}-\eqref{eq:Sommerfeld}, for $j\in \llbracket0,3\rrbracket$. We showed that $\varphi^j_\theta$, defined by \eqref{eq:definition_phi_theta} as the analytic continuation of $\varphi^j$ from the real line to the path $\tau_\theta$ of Figure \ref{fig:path}, satisfies \eqref{eq:syst_comp_cmplx} and \eqref{eq:syst_comp_BC_cmplx}, equivalently \eqref{eq:HSMM_matrice_cmplx}. A key component in this argument was to deform paths of integration from the real line to the path $\tau_\theta$ (Theorem \ref{th:deformed_representation}).

In the following uniqueness proof we reverse this derivation. We show that if the $\varphi^j_\theta$ satisfy \eqref{eq:syst_comp_cmplx} and \eqref{eq:syst_comp_BC_cmplx}, with $g=0$, then, in the sense that \eqref{eq:definition_phi_theta} holds, they are the analytic continuations onto the path $\tau_\theta$ of functions $\varphi^j$ that satisfy the system \eqref{eq:HSMM_complexe} with $g=0$. (A key component in this argument is a deformation of paths of integration from $\tau_\theta$ back to the real line, justified as in the proof of Theorem \ref{th:deformed_representation}.)   Moreover, by an application of Proposition \ref{prop:FF_sol} (justified by Corollary \ref{cor-Dtheta(t)infty}), the functions $\varphi^j$ satisfy the radiation conditions \eqref{eq:RadCond}, so that $\varphi^j=0$ by Theorem \ref{thm:uniqueHSM}. Thus $\varphi^j_\theta$, which is the analytic continuation of $\varphi^j$, is also zero.

%  the $\phi^j$'s in \eqref{eq:HSMM_complexe} have analytic continuations from the real line to the path $\tau_\theta$ of \eqref{eq:parametrage_t_theta} and Figure \ref{fig:path}, and defining $\phi^j_\theta$ in \eqref{eq:definition_phi_theta} to be that analytic continuation;

 		\begin{theorem}
 			\label{the-uniqueness}
 				Suppose that $\{\varphi^0_\theta,\varphi^1_\theta,\varphi^2_\theta,\varphi^3_\theta\}\in (L^2(\R))^4$  is a solution of \eqref{eq:syst_comp_cmplx} such that
 			\begin{equation} \label{eq:hom}
 			\varphi_j^\theta(t)=0,\;-a<t<a,\quad \;j\in \llbracket 0,3\rrbracket.
 			\end{equation}
 	Then $\varphi_j^\theta =0$ for $j\in  \llbracket 0,3\rrbracket$.
 		\end{theorem}
 %	As for the dissipative case, the idea of the proof is to construct from the $\varphi^j_\theta$ a $H^1_\text{loc}$ solution $u$ of the Helmholtz equation in $\Omega$ which satisfies a vanishing Dirichlet boundary condition on $\partial\Omega$ and the radiation condition at infinity. By well-posedness of this problem, we conclude that $u=0$ and deduce that $\varphi_j^\theta=0$.
 \begin{proof}
 %Here are the steps of the proof. 	
 %\textcolor{red}{Step 1.}
 Suppose that $\{\varphi^0_\theta,\varphi^1_\theta,\varphi^2_\theta,\varphi^3_\theta\}\in (L^2(\R))^4$  is a solution of \eqref{eq:syst_comp_cmplx} that satisfies \eqref{eq:hom}, and define the functions $\varphi^0,\varphi^1,\varphi^2,\varphi^3$  by
 		\begin{equation}\label{eq:def_varphij}
 		\begin{array}{|ll}
 		\varphi^j(t):= S\,\widetilde{D}_\theta\varphi_\theta^{j-1}(t),&t<-a,\\[3pt]
 		\varphi^j(t):=0,&-a\leq t\leq a,\\[3pt]
 		\varphi^j(t):= \widetilde{D}_\theta S\,\varphi_\theta^{j+1}(t),&t>a,
 		\end{array}
 		\end{equation}
for $j\in  \llbracket 0,3\rrbracket$,	where 	$\widetilde{D}_\theta$ is defined by \eqref{eq:opDthetainter}. As discussed below \eqref{eq:opDthetainter}, since the $\varphi_\theta^j$'s are in $L^2(\R)$, it follows from this definition, the definition \eqref{eq:sym} of $S$, and Lemma \ref{lem:analyticity_utheta}, that the $\varphi^j$'s have analytic   continuations from $(a,+\infty)$ (respectively $(-\infty,-a)$) to the part of the half-plane $\Re(z)>a$ with $-\pi/2+\theta<\Arg(z-a)<\pi/2$ (respectively $\Re(z)<-a$ with $-\pi/2+\theta<\Arg(-z-a)<\pi/2$). In particular, by (\ref{eq:Dtilde-D}), and recalling \eqref{eq:sym} and that $\tau_\theta(-t)=\tau_\theta(t)$ for $t<-a$,
  	\begin{equation*}%\label{eq:def_varphij-tautheta}
  \begin{array}{|ll}
  \varphi^j(\tau_\theta(t))= S \,{D}_\theta\varphi_\theta^{j-1}(t),&t<-a,\\[3pt]
  \varphi^j(\tau_\theta(t))=0,&-a\leq t \leq a,\\[3pt]
  \varphi^j(\tau_\theta(t))= {D}_\theta S\,\varphi_\theta^{j+1}(t),&t>a,
  \end{array}
  \end{equation*}
for $j\in  \llbracket 0,3\rrbracket$. Comparing this equation with \eqref{eq:syst_comp_cmplx} we see that
\begin{equation} \label{eq:same}
\varphi^j_\theta (s) = \varphi^j \big(\tau_\theta(s)\big)
  \quad\text{for } s \in \R \text{ and } j \in \llbracket0,3\rrbracket.
\end{equation}
Thus, for $s>a$ and for $s<-a$, the $\varphi^j_\theta$'s are the analytic continuations of the $\varphi^j$'s to the complex path parametrized by $\tau_\theta$. Therefore, to complete the proof of the theorem it is enough to show, for $j\in  \llbracket 0,3\rrbracket$, that $\varphi^j(t)=0$ for  $t\in \R$ with $|t|>a$, for this will imply,  by the uniqueness of analytic continuation, that $\varphi^j(\tau_\theta(s))=0$ for $s>a$ and $s<-a$, which will imply, using \eqref{eq:same},  that each $\varphi^j_\theta=0$.

To establish, for $j\in  \llbracket 0,3\rrbracket$, that $\varphi^j(t)=0$ for $t>a$ and $t<-a$, we show, via an application of Cauchy's integral theorem as in the proof of Theorem \ref{th:deformed_representation}, that uses the analyticity of the $\varphi^j$'s noted above, Lemma \ref{lem-Dtheta(t)infty}, and \eqref{eq:same}, that \eqref{eq:def_varphij} implies that the $\varphi^j$ satisfy \eqref{eq:HSMM_complexe} with $g=0$, i.e.
%Second step: using Lemma \ref{lem-Dtheta(t)infty},
%  we can then adapt the proof of Theorem \ref{th:deformed_representation} in order to deform the path of integration in \eqref{eq:def_varphij-tautheta} back to the real axis and show that
%  $\{\varphi^0,\varphi^1,\varphi^2,\varphi^3\}$ satisfies
 \begin{equation*}%\label{eq:HsMM_real}
 \begin{array}{|lcl}
 \varphi^j(t)=S\,D\,\varphi^{j-1}(t),\quad t<-a,\\[3pt]
 \varphi^j(t)=0,\quad-a<t<a,\\[3pt]
 \varphi^j(t)=D\,S\,\varphi^{j+1}(t),\quad t>a,
 \end{array}%\quad j\in \llbracket0,3\rrbracket
 \end{equation*}
for $j\in \llbracket0,3\rrbracket$, where $D$ is defined in \eqref{eq:opD_comp_exprG}, %for a real frequency
 or equivalently by \eqref{eq:opD_comp_exprG_cmplx} with $\theta=0$. (The key step in the argument is to show, as we have done in the proof of Theorem \ref{th:deformed_representation}, that
 $$
 \int_a^{+\infty} h(t-a,a\mp \tau_\theta(s))\varphi(s)\tau_\theta'(s)\, \rd s = \int_a^{+\infty}  h(t-a,a\mp s)\varphi(s)\, \rd s
 $$
 if, for every $0<\delta<M$, $\varphi$ is analytic in the domain $T_\theta^{\delta,M}$ introduced in the proof of Theorem \ref{th:deformed_representation}, and continuous in its closure, and if the behaviour of $\varphi(z)$ as $z\to 0$ and as $z\to+\infty$ is suitably constrained, satisfying the bounds of Lemma \ref{lem-Dtheta(t)infty}.)
 We note also that it follows from \eqref{eq:def_varphij} and Lemma \ref{lem-Dtheta(t)infty2} that $\varphi^j\in L^2_{\mathrm{loc}}(\R)$, for $j\in \llbracket0,3\rrbracket$. Moreover, the $\varphi_\theta^j$ are in $L^2(\R)$ and it follows, from \eqref{eq:syst_comp_cmplx} and Corollary \ref{cor-Dtheta(t)infty}, that each $\varphi_\theta^j$ decays exponentially at infinity, so that $(t^2+1)\varphi_\theta^j(t)\in L^1(\R)$, $j\in \llbracket0,3\rrbracket$. Thus we can  apply Proposition \ref{prop:FF_sol} to \eqref{eq:def_varphij}, noting the second equality in \eqref{eq:Dtilde-D}, to see that the radiation conditions \eqref{eq:RadCond} hold.
  But this is enough to conclude that $\varphi^j=0$, for $j\in \llbracket0,3\rrbracket$, by Theorem \ref{thm:uniqueHSM}.
 \end{proof}
  %the following  uniqueness theorem for the HSM formulation for real $k$  proved in \cite{bibid}:
%\vspace{0.2ex}	
% Since the $\varphi^j_\theta$'s are the analytic continuations of the $\varphi^j$'s, they also vanish indentically.

\section{The {complex-scaled} HSM method for the general case}
\label{sec-CHSMgeneralcase}
Let us now explain how to extend the  {complex-scaled} HSM method of section \ref{sec-CHSM} to solve the general problem presented in the introduction. More precisely, for a real wavenumber $k>0$ and for a function $\rho$ and a subdomain $\Omega$ of $\mathbb{R}^2$ satisfying the hypotheses described in section \ref{sec-introduction}, the objective is to derive a complex-scaled HSM formulation to compute the solution $u\in H^1_{\mathrm{loc}}(\Omega)$ of  (\ref{pb:probleme_diff_initial}) and (\ref{eq:Sommerfeld}). For the sake of simplicity, we will restrict attention to  the case $\Omega=\mathbb{R}^2$, but adding bounded obstacles contained in $\Omega_a$ is completely straightforward (see Section \ref{sec:NumGC}). As in section \ref{sec-introduction}, the source term $f$ is in $L^2(\Omega)$, with compact support that is a subset of $\Omega_a$.

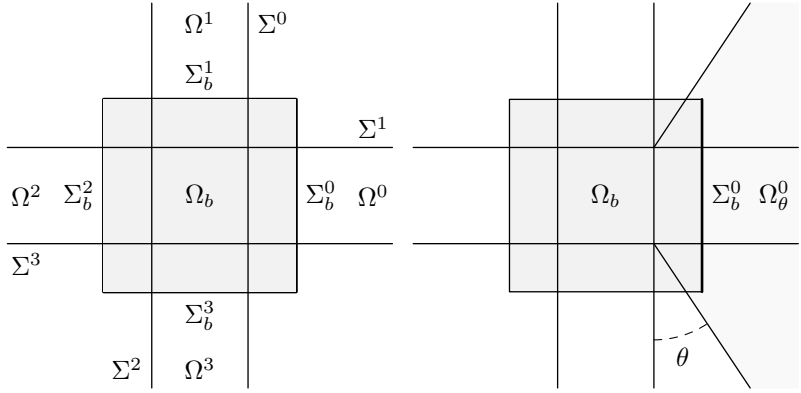
\begin{figure}[!h]
	\definecolor{ffffff}{rgb}{0.2,0.2,0.2}
	\definecolor{yqyqyq}{rgb}{0.5,0.5,0.5}
	\begin{center}
		\begin{tikzpicture}[line cap=round,line join=round,x=0.4cm,y=0.4cm,scale=.8,yshift=2]
		\begin{scope}[shift={(-15,0)}]
		\draw  [line width=1.2pt] (2.,2.) rectangle (10,10);
		\fill [color=gray!10] (2.,2.) rectangle (10,10);
		%		\draw [dashed,->] (6,6)--(7.5,6);
		%		\draw[below] (7.5,6) node {$x_1$};
		%		\draw[left] (6,7.2) node {$x_2$};
		%		\draw [dashed,->] (6,6)--(6,7.3);
		%		\draw (5.2,5.2) node {$\Omega_a$};
		\draw (6,6) node {$\Omega_b$};
		\draw[above] (13.2,8) node {$\Sigma^1$};
		\draw[right] (8,13.2) node {$\Sigma^0$};
		\draw[left] (4,-1.2) node {$\Sigma^2$};
		\draw[below] (-1.2,4) node {$\Sigma^3$};
		%\draw[right] (8,6) node {$\Sigma^0_a$};
		\draw[right] (10,6) node {$\Sigma^0_b$};
		%	\draw[above] (6,8) node {$\Sigma^1_a$};
		\draw[above] (6,10) node {$\Sigma^1_{b}$};
		%	\draw[left] (4,6) node {$\Sigma^2_a$};
		\draw[left] (2,6) node {$\Sigma^2_b$};
		%	\draw[below] (6,4) node {$\Sigma^3_a$};
		\draw[below] (6,2) node {$\Sigma^3_b$};
		\draw  (13.2,6) node {$\Omega^0$};
		\draw  (6,13.2) node {$\Omega^1$};
		\draw  (-1.2,6) node {$\Omega^2$};
		\draw  (6,-1.2) node {$\Omega^3$};
		\draw [line width=.5pt] (-2.,4)--(14.,4);
		\draw [line width=.5pt] (-2.,8)--(14.,8);
		\draw [line width=.5pt] (4,-2)--(4,14.);
		\draw [line width=.5pt] (8,-2)--(8,14.);
		\end{scope}
		\end{tikzpicture}\hspace{2mm}\begin{tikzpicture}[line cap=round,line join=round,x=0.4cm,y=0.4cm,scale=.8,yshift=2]
		\begin{scope}[shift={(-15,0)}]
		\filldraw[color=gray!5] (8,8)--(12,14)--(14,14)--(14,-2)--(12,-2)--(8,4)--(8,8);
		\fill [color=gray!10] (2.,2.) rectangle (10,10);
		\draw  [line width=.5pt] (2.,2.) rectangle (10,10);
		\draw [line width=.5pt] (8,8)--(12,14);
		\draw [line width=.5pt] (12,-2)--(8,4);
		%		\draw [dashed,->] (6,6)--(7.5,6);
		%		\draw[below] (7.5,6) node {$x_1$};
		%		\draw[left] (6,7.2) node {$x_2$};
		%		\draw [dashed,->] (6,6)--(6,7.5);
		%		\draw (5.5,5.5) node {$\Omega_a$};
		\draw (6,6) node {$\Omega_b$};
		\draw [line width=1.pt] (10,2)--(10,10.);
		\draw[right] (10,6) node {$\Sigma^0_b$};
		\draw  (13,6) node {$\Omega^0_\theta$};
		\draw [line width=.5pt] (-2.,4)--(14.,4);
		\draw [line width=.5pt] (-2.,8)--(14.,8);
		\draw [line width=.5pt] (4,-2)--(4,14.);
		\draw [line width=.5pt] (8,-2)--(8,14.);
		\draw [dashed] (8,0) arc (270:304:4) ;
		\draw  (9.2,-.7) node {$\theta$};
		\end{scope}
		\end{tikzpicture}
	\end{center}
	\caption{The notations for the general case.}\label{Fig:notations_Case2}
\end{figure}

The idea is to introduce, in addition to the lines $\Sigma^j$, $j\in\llbracket0,3\rrbracket$, a square $\Omega_b:=(-b,b)^2$ for some $b>a$. (As in section \ref{section-HSMcomplexfreq} we set $\Sigma_b:= \partial \Omega_b$ and denote the sides of $\Sigma_b$ by $\Sigma_b^j$, $j\in\llbracket0,3\rrbracket$; see Figure \ref{Fig:notations_Case2}.)
We will show how to derive a formulation of problem (\ref{pb:probleme_diff_initial},\ref{eq:Sommerfeld}) whose unknowns are the complex-scaled traces $\varphi_\theta^j$, $j\in\llbracket0,3\rrbracket$, associated to the infinite lines $\Sigma^j$, and the restriction $u_b:=u|_{\Omega_b}$ of the solution $u$ to the square $\Omega_b$. To do that, we need to make the following assumption on the parameter $\theta$:
\begin{equation}
\label{eq-thetainfpisur4}
\theta<\frac{\pi}{4}.
\end{equation}
Let us derive the equations linking the $\varphi_\theta^j$, $j\in\llbracket0,3\rrbracket$, and $u_b$. On the one hand, the $\varphi_\theta^j$ still satisfy the system of compatibility relations (\ref{eq:syst_comp_cmplx}). But, instead of (\ref{eq:syst_comp_BC_cmplx}), we have to impose equality between  $\varphi_\theta^j$ and $u_b$ on $\Sigma^j_a$:
\begin{equation*}
%\label{eq-varphijthetaegalub}
\varphi_\theta^j(t)={u_b}|_{\Sigma^j_a}(t),\quad -a<t<a,\:j\in\llbracket0,3\rrbracket.
\end{equation*}
On the other hand, we can derive a variational formulation for $u_b$ in $\Omega_b$. Since $-\Delta u_b-k^2\rho u_b=f$ in $\Omega_b$ {and $f$ is supported in $\Omega_a$}, the following Green's identity holds for all $v_b\in H^1(\Omega_b)$, where $n$ is the normal  unit vector pointing out of $\Omega_b$:
\begin{equation}
\label{FVub1}
\int_{\Omega_b}\left(\nabla u_b\cdot \overline{\nabla v_b}-k^2 \rho u_b  \overline{v_b}\right)-\int_{\Sigma_b}\frac{\partial u_b}{\partial n}\overline{v_b}\, =\int_{\Omega_{{a}}}f\overline{v_b}.
\end{equation}
The last idea is to replace  in the previous identity the normal derivative of $u_b$  on the $j$th side of the square  by an integral representation as a function of $\varphi_\theta^j$. Indeed, %if we denote by $\Sigma^j_b$ the side of the square included in $\Omega^j$,
we must have, for $j\in\llbracket0,3\rrbracket$,
$$\frac{\partial u_b}{\partial n}-\ri ku_b=\frac{\partial U^j}{\partial n}-\ri kU^j\mbox{ on }\Sigma^j_b,$$
where $U^j$  denotes the restriction of the solution $u$ to the half-plane $\Omega^j$. {(Our choice of Robin traces instead of normal derivatives is so that later we have uniqueness for the boundary value problem \eqref{pb:interiorbvp} for all $k>0$.)}%\scednote{Is this OK? Helpful? Or better to just delete this sentence?}

We have proved in Theorem \ref{th:deformed_representation} that $U^j(\bsx^j)$ has an integral representation in terms of $\varphi_\theta^j$, as soon as $\bsx^j$ belongs to the domain $\Omega^j_\theta$ defined by (\ref{eq:Omegaj_theta}).
Notice that under condition (\ref{eq-thetainfpisur4}) one has $\Sigma^j_b\subset \Omega^j_\theta$ (see Figure \ref{Fig:notations_Case2}).
%\begin{figure}[!h]
%	\definecolor{ffffff}{rgb}{0.2,0.2,0.2}
%	\definecolor{yqyqyq}{rgb}{0.5,0.5,0.5}
%	\begin{center}
%		\begin{tikzpicture}[line cap=round,line join=round,x=0.4cm,y=0.4cm,scale=1,yshift=2]
%		\begin{scope}[shift={(-15,0)}]
%		\filldraw[color=gray!5] (8,8)--(12,14)--(14,14)--(14,-2)--(12,-2)--(8,4)--(8,8);
%		\fill [color=gray!10] (2.,2.) rectangle (10,10);
%		\draw  [line width=.5pt] (2.,2.) rectangle (10,10);
%		\draw [line width=.5pt] (8,8)--(12,14);
%		\draw [line width=.5pt] (12,-2)--(8,4);
%		%		\draw [dashed,->] (6,6)--(7.5,6);
%		%		\draw[below] (7.5,6) node {$x_1$};
%		%		\draw[left] (6,7.2) node {$x_2$};
%		%		\draw [dashed,->] (6,6)--(6,7.5);
%		%		\draw (5.5,5.5) node {$\Omega_a$};
%		\draw (6,6) node {$\Omega_b$};
%		\draw [line width=1.pt] (10,2)--(10,10.);
%		\draw[right] (10,6) node {$\Gamma^0_b$};
%		\draw  (13,6) node {$\Omega^0_\theta$};
%		\draw [line width=.5pt] (-2.,4)--(14.,4);
%		\draw [line width=.5pt] (-2.,8)--(14.,8);
%		\draw [line width=.5pt] (4,-2)--(4,14.);
%		\draw [line width=.5pt] (8,-2)--(8,14.);
%		\draw [dashed] (8,0) arc (270:304:4) ;
%		\draw  (9.2,-.7) node {$\theta$};
%		\end{scope}
%		\end{tikzpicture}
%	\end{center}
%	\caption{The inclusion $\Gamma^0_b\subset \Omega^0_\theta$ for $\theta<\pi/4$.}\label{Fig:Gamma0b-subset-Omega0theta}
%\end{figure}
Consequently, one can use the formula (\ref{eq-repinOmegatheta}) to rewrite the above Robin compatibility condition on $\Sigma^j_b$. Precisely, we have
$$\begin{array}{l} \displaystyle \left(\frac{\partial u_b}{\partial n}-\ri ku_b\right)\Big|_{\Sigma^j_b}(t) =\\\displaystyle\int_{\R}
\left(\partial_1h(b-a,t-\tau_\theta(s))-\ri kh(b-a,t-\tau_\theta(s))\right)\varphi^j_\theta(s)\tau'_\theta(s)\, \rd s,\quad -b<t<b,
\end{array}$$
where  $\partial_1h$ denotes the derivative of $h$ with respect to its first variable. This leads us to define the following Dirichlet-to-Robin operator $\Lambda_\theta$. For $\psi\in L^2(\R)$ and $0<\theta<\pi/4$, define $\Lambda_\theta \psi\in L^2(-b,b)$ by
\begin{equation}
\label{eq-Lambda-theta}
\Lambda_\theta \psi(t):=\int_{\R}\lambda(b-a,t-\tau_\theta(s))\psi(s)\tau'_\theta(s)\,\rd s, \quad -b<t<b,
\end{equation}
where we have set
$$\lambda(x_1,z):=\partial_1h(x_1,z)-\ri kh(x_1,z), \quad x_1>0, \;\; z\in \C,$$
which one can easily check, using \cite[(10.6.2)]{DLMF} and the definition \eqref{eq:hpr_kernel} of $h$, %that this kernel
takes the explicit form
\begin{equation}
\label{eq-kernel-Lambda-theta}
\lambda(x_1,z)=\frac{\ri k}{2R}\left(\left[1-\ri kx_1\right]H_1^{(1)}(kR)-\frac{kx_1^2}{R}H_2^{(1)}(kR)\right),\;\; R=[x_1^2+z^2]^{1/2}.
\end{equation}
With this notation, the previous equations linking $u_b$ and the $\varphi^j_\theta$ can be written as
\begin{equation}
\label{eq-condDtN}
 \left(\frac{\partial u_b}{\partial n}-\ri ku_b\right)\Big|_{\Sigma^j_b}(t)=\Lambda_\theta \varphi^j_\theta(t) ,\quad -b<t<b, \;\; j\in\llbracket0,3\rrbracket.
\end{equation}

Our complete formulation reads as follows:
\begin{equation}\label{eq:HSMM_casgeneral}
\begin{array}{c}
\text{Find}\; \{\varphi^0_\theta,\varphi^1_\theta,\varphi^2_\theta,\varphi^3_\theta\}\in (L^2(\R))^4\;\text{and}\;u_b\in H^1(\Omega_b)\;\text{such that}\\[5pt]
\begin{array}{|lcl}
\varphi^j_\theta(t)=S\,D_\theta\,\varphi^{j-1}_\theta(t)\;\text{for}\;t<-a,\\[3pt]
\varphi^j_\theta(t)=u_b|_{\Sigma_a^j}(t)\;\text{for}\;-a<t<a,\\[3pt]
\varphi^j_\theta(t)=D_\theta\,S\,\varphi_\theta^{j+1}(t)\;\text{for}\;t>a.
\end{array} \quad j\in \llbracket0,3\rrbracket,\\[20pt]
\text{and such that, }\forall v_b\in  H^1(\Omega_b),\hfill\\
\displaystyle\int_{\Omega_b}\left(\nabla u_b\cdot \overline{\nabla v_b}-k^2 \rho u_b  \overline{v_b}\right)-\ri k\sum_{j=0}^3\int_{\Sigma^j_b} u_b  \overline{v_b}\\[6pt] \displaystyle -\sum_{j=0}^3\int_{-b}^b\Lambda_\theta \varphi^j_\theta(t)\overline{v_b}|_{\Sigma^j_b}(t)\rd t=\int_{\Omega_{{a}}}f\overline{v_b},
\end{array}
\end{equation}
where $S$, $D_\theta$, and $\Lambda_\theta$ are defined by (\ref{eq:sym}), (\ref{eq:opD_comp_exprG_cmplx}), and (\ref{eq-Lambda-theta}).

Let us denote by $ \Phi_\theta$, $\Phi(u_b)$, and
{$\widetilde\Phi_\theta$} the following elements of $(L^2(\R))^4$:
\begin{equation*}%\label{eq:inconnues_casgeneral}
\begin{array}{l}
\Phi_\theta:=\{\varphi^0_\theta,\varphi^1_\theta,\varphi^2_\theta,\varphi^3_\theta\},\;\;
\Phi(u_b):=\{u_b|_{\Sigma_a^0},u_b|_{\Sigma_a^1},u_b|_{\Sigma_a^2},u_b|_{\Sigma_a^3}\},\\[5pt]
\widetilde{\Phi}_\theta= \{\widetilde{\varphi}^0_\theta,\widetilde{\varphi}^1_\theta,\widetilde{\varphi}^2_\theta,\widetilde{\varphi}^3_\theta\}:=\Phi_\theta-\Phi(u_b).
\end{array}
\end{equation*}%\textcolor{magenta}{Je pense qu'on ne peut pas échapper à cette nouvelle inconnue}
Note that, using the second equation of (\ref{eq:HSMM_casgeneral}) and recalling \eqref{L20def}, we have $ \widetilde{\Phi}_\theta\in (L^2_0(\R))^4$. With these notations the first block of equations in (\ref{eq:HSMM_casgeneral}) can be rewritten (cf.\ \eqref{eq:HSMM_matrice_cmplx}) as $(\mathbb{I}-\mathbb{D}_\theta) \widetilde{\Phi}_\theta=\mathbb{D}_\theta\Phi(u_b)$, where $\mathbb{D}_\theta$ is defined by (\ref{eq:def_Dmat_cmplx}), so that problem (\ref{eq:HSMM_casgeneral}) can be rewritten as
\begin{equation}\label{eq:HSMM_casgeneral-variational}
\begin{array}{c}
\text{Find}\;  \widetilde{\Phi}_\theta\in (L^2_0(\R))^4\;\text{and}\;u_b\in H^1(\Omega_b)\;\text{such that}\\[5pt]
\forall   \widetilde{\Psi}\in (L^2_0(\R))^4\;\text{and}\;\forall v_b\in H^1(\Omega_b),\\[5pt]
((\mathbb{I}-\mathbb{D}_\theta) \widetilde{\Phi}_\theta, \widetilde{\Psi})_{(L^2_0(\R))^4}-(\mathbb{D}_\theta\Phi(u_b), \widetilde{\Psi})_{(L^2_0(\R))^4}\\[5pt]
\displaystyle +\int_{\Omega_b}\left(\nabla u_b\cdot \overline{\nabla v_b}-k^2 \rho u_b  \overline{v_b}\right)-\ri k\sum_{j=0}^3\int_{\Sigma^j_b} u_b  \overline{v_b}\\[6pt] \displaystyle
-\sum_{j=0}^3\int_{-b}^b\Lambda_\theta\left( \widetilde{\varphi}^j_\theta+{u_b}|_{\Sigma^j_{a}}\right)(t)\,\overline{v_b}|_{\Sigma^j_b}(t)\, \rd t
=\int_{\Omega_{{a}}}f\overline{v_b}.
\end{array}
\end{equation}
%where $\mathbb{D}_\theta$ is defined by (\ref{eq:def_Dmat_cmplx}).
Let us first prove a uniqueness result for this problem:
\begin{prop}
	\label{lem-uniqueness-generalcase}
If $f=0$ then the only solution of problem (\ref{eq:HSMM_casgeneral-variational}) is the trivial solution $u_b= 0$, $\widetilde{\Phi}_\theta= 0$.
\end{prop}
\begin{proof}
Suppose that $\widetilde{\Phi}_\theta\in (L^2_0(\R))^4$ and $u_b\in H^1(\Omega_b)$ are such that (\ref{eq:HSMM_casgeneral-variational}) holds with $f=0$. Then $(\mathbb{I}-\mathbb{D}_\theta) \widetilde{\Phi}_\theta=\mathbb{D}_\theta\Phi(u_b)$, and the second of equations \eqref{eq:HSMM_casgeneral} holds with $f=0$ and $\varphi_\theta^j:= \widetilde \varphi^j_\theta + {u_b}|_{\Sigma^j_a}$, $j\in \llbracket0,3\rrbracket$, so that $-\Delta u_b-k^2\rho u_b=0$ in $\Omega_b$ and \eqref{eq-condDtN} holds for $j\in \llbracket0,3\rrbracket$.

Let us denote by $u_\infty\in H^1_{\mathrm{loc}}(\R^2\backslash \overline{\Omega_a})$ the unique solution of \eqref{pb:probleme_Dir_real}-(\ref{eq:Sommerfeld}) with
$g={u_b}|_{\Sigma_a}\in H^{1/2}(\Sigma_a)$. Then, as we have shown in section \ref{sec:Derivation}, the vector of complex-scaled traces of $u_\infty$ satisfies \eqref{eq:HSMM_matrice_cmplx} with $\Phi_g = \Phi(u_b)$, so that, by % so that $u_\infty\in H^1_{loc}(\R^2\backslash \overline{\Omega_a})$.
Theorem \ref{the-uniqueness}, it coincides with the vector $\Phi_\theta:= \{\varphi^0_\theta,\varphi^1_\theta,\varphi^2_\theta,\varphi^3_\theta\}$. Thus, applying Theorem \ref{th:deformed_representation} (as we did above to derive \eqref{eq-condDtN} from \eqref{pb:probleme_diff_initial}-\eqref{eq:Sommerfeld}) we see that
%coincide with the $\varphi^j_\theta$'s are the complex-scaled traces of $u_\infty$,  so that
 $$\Lambda_\theta\varphi^j_\theta= \left(\frac{\partial u_\infty}{\partial n}-\ri ku_\infty\right)\Big|_{\Sigma^j_b}, \quad j\in \llbracket0,3\rrbracket.$$
 Consequently, $v:=u_b-u_\infty$ belongs to $ H^1(\Omega_b\backslash \overline{\Omega_a})$ and satisfies
\begin{equation} \label{pb:interiorbvp}
\begin{array}{|lcr}
	\Delta v +  k^2\, v= 0 & \text{in} & \Omega_b\backslash \overline{\Omega_a},\\[4pt]
v=0&\text{on}&\Sigma_a,\\[4pt]
\displaystyle \frac{\partial v}{\partial n}-\ri kv=0&\text{on}&\Sigma_b.
\end{array}
\end{equation}
But, for every $k>0$, this homogeneous problem has no solution except $v=0$. {(To see this apply Green's identity (cf.\ \eqref{FVub1}) in $\Omega_b\backslash \overline{\Omega_a}$  to deduce that $\int_{\Sigma_b} |v|^2 =0$, so that $v=\partial v/\partial n=0$ on $\Sigma_b$, and apply Holmgren's uniqueness theorem; \cite[p.~104]{ChGrLaSp:11}.)} Thus $u_b=u_\infty$ in $\Omega_b\backslash \overline{\Omega_a}$ so that the function
$$w:=\left\{\begin{array}{lcr}
u_b & \text{in} & \Omega_b,\\[4pt]
u_\infty&\text{in}&\R^2\backslash \overline{\Omega_a},
\end{array}\right.$$
is well-defined, and is a solution of the homogeneous Helmholtz equation in $\R^2$ which satisfies the radiation condition (\ref{eq:Sommerfeld}). As a consequence (e.g., \cite[Theorem 8.7]{CoKr:98}), $w= 0$ in $\R^2$, so that $u_b=0$ in $\Omega_b$ and $u_\infty=0$. Further, each $\varphi_\theta^j$ is zero, since it is  a complex trace of $u_\infty$. This completes the proof.
\end{proof}

The well-posedness of problem (\ref{eq:HSMM_casgeneral-variational}) follows from classical arguments combined with the results of the previous sections and the following lemma:
\begin{lem}
	\label{lem-DtNcompact}
	The operator $\Lambda_\theta$, defined by (\ref{eq-Lambda-theta}), is a compact operator from $L^2(\R)$ to $L^2(-b,b)$.
\end{lem}%\textcolor{magenta}{En fait, on s'en moque de la compacité ?}
\begin{proof}
	To establish this lemma, we will prove that the kernel of
	 $\Lambda_\theta$  is Hilbert-Schmidt, i.e.\ that
	\begin{equation} \label{eq:isHS}
	(t,s)\mapsto \lambda(b-a,t-\tau_\theta(s))\;\in L^2((-b,b)\times\R),
	\end{equation}
	by using similar arguments as in the proof of Proposition \ref{prop:D_realk_cmplx_moins}. On the one hand, since $b>a$ and $\theta<\pi/4$, the function $R$ which appears  in expression
	(\ref{eq-kernel-Lambda-theta}) never vanishes for $-b\leq t\leq b$ and $s\in\R$. As a consequence, the kernel is continuous on $[-b,b]\times \R$. %a smooth function of the variables $t$ and $s$.
On the other hand we have, {using \eqref{eq:Rasymp1}, \eqref{eq:Rasymp2}, and \eqref{eq:HankAsym}, the asymptotic} estimate
	$$|\lambda(b-a,t-\tau_\theta(s))|=\mathcal{O}\left(\frac{\re^{-k\sin \theta |s|}}{|s|^{1/2}}\right) \quad\text{as }|s|\rightarrow +\infty,$$
	uniformly in $t$, for $|t|\leq b$. Together, these properties prove \eqref{eq:isHS}. %which ends the proof.
	%\textcolor{magenta}{J'ai l'impression  qu'il y a trop de calculs qui se répètent sans homogénéité partout...}
\end{proof}
	
Recall that we call a sesquilinear form $a(\cdot,\cdot)$ on a Hilbert space $\mathcal{H}$  {\em compact} if the associated linear operator $A$ on $\mathcal{H}$, defined by $(A\phi,\psi)=a(\phi,\psi)$, for all $\phi,\psi\in \mathcal{H}$, is compact. (Here $(\cdot,\cdot)$ denotes the inner product on $\mathcal{H}$.) Equivalently, $a(\cdot,\cdot)$ is compact is, whenever $\phi_n\rightharpoonup0$ and $\psi_n\rightharpoonup0$ (weak convergence in $\mathcal{H}$), it holds that $a(\phi_n,\psi_n)\to 0$.

In our final theorem we show that the sesquilinear form on {the Hilbert space} $(L^2_0(\R))^4\times H^1(\Omega_b)$ which appears on the left-hand side of  (\ref{eq:HSMM_casgeneral-variational}) is the sum of coercive plus compact sesquilinear forms, and that, as a consequence of this and of the above uniqueness result,  problem (\ref{eq:HSMM_casgeneral-variational}) is well-posed. {Regarding the last sentence of the theorem, note that, given the constraint $0<\theta<\pi/4$,
\begin{equation*}% \label{eq:OmegaContained}
\Omega\subset \Omega_b \cup \bigcup_{j=0}^3 \Omega^j_\theta,
\end{equation*}
so that the solution of the original scattering problem can be recovered in the whole of $\Omega$ from the solution $(\widetilde \Phi_\theta,u_b)$ of (\ref{eq:HSMM_casgeneral-variational}).}

\begin{theorem}
	\label{the-generalcase}
{For every $\theta\in (0,\pi/4)$ and every $f\in L^2(\Omega)$ with support in $\Omega_a$,  Problem (\ref{eq:HSMM_casgeneral-variational}) has exactly one solution $(\widetilde \Phi_\theta,u_b)\in (L^2_0(\R))^4\times H^1(\Omega_b)$. Further, for some constant $c>0$ depending on $\theta$,
\begin{equation} \label{eq:stabilityGen}
\|\widetilde \Phi_\theta\|_{(L^2_0(\R))^4}+\|u_b\|_{H^1(\Omega_b)} \leq c\|f\|_{L^2(\Omega_a)},
\end{equation}
for all $f\in L^2(\Omega)$ with support in $\Omega_a$. Moreover, if $(\widetilde \Phi_\theta,u_b)\in (L^2_0(\R))^4\times H^1(\Omega_b)$ is the solution of (\ref{eq:HSMM_casgeneral-variational}) and $u\in H^1_{\mathrm{loc}}(\Omega)$ is the solution of
\eqref{pb:probleme_diff_initial}-\eqref{eq:Sommerfeld}, then $u=u_b$ in $\Omega_b$, while $u$ is given in terms of $\varphi_\theta^j$ in $\Omega_\theta^j$ by \eqref{eq-repinOmegatheta}, for $j\in \llbracket0,3\rrbracket$. }
\end{theorem}
\begin{proof}
	{As an operator on $(L^2_0(\R))^4$ we have, by Theorem \ref{th:HSMM_cmplxed}, that $\mathbb{D}_\theta=\mathbb{D}_\theta^1 + \mathbb{D}_\theta^2$, where $\|\mathbb{D}_\theta^1\|\leq 1/\sqrt{2}$ and $\mathbb{D}_\theta^2$ is compact. Thus}  we can decompose the sesquilinear form which appears on the left-hand side of  (\ref{eq:HSMM_casgeneral-variational}) as the sum of a first sesquilinear form
\begin{eqnarray*}
((\mathbb{I}-{\mathbb{D}^1_\theta}) \widetilde{\Phi}_\theta, \widetilde{\Psi})_{(L^2_0(\R))^4}+  \int_{\Omega_b}\left(\nabla u_b\cdot \overline{\nabla v_b}+ u_b  \overline{v_b}\right)
-\ri k\sum_{j=0}^3\int_{\Sigma^j_b} u_b  \overline{v_b},
\end{eqnarray*}
	which is coercive on $(L^2_0(\R))^4\times H^1(\Omega_b)$, % thanks to Theorem \ref{th:HSMM_cmplxed},
	and a second sesquilinear form
\begin{eqnarray*}
&{-(\mathbb{D}^2_\theta}\widetilde \Phi_\theta, \widetilde{\Psi})_{(L^2_0(\R))^4} -(\mathbb{D}_\theta\Phi(u_b), \widetilde{\Psi})_{(L^2_0(\R))^4} -(1+k^2)\int_{\Omega_b} u_b  \overline{v_b}\\
& \hspace{20ex}  -\sum_{j=0}^3\int_{-b}^b\Lambda_\theta\left( \widetilde{\varphi}^j_\theta+{u_b}|_{\Sigma^j_b}\right)(t)\,\overline{v_b}|_{\Sigma^j_b}(t)\, \rd t
\end{eqnarray*}
	which is compact on the same space. The proofs of compactness of the four terms of this second sesquilinear form rely on different arguments. {The first is compact because $\mathbb{D}_\theta^2$ is compact.} For the second term, we notice that the operator $u_b\longmapsto \mathbb{D}_\theta\Phi(u_b)$ is compact from $ H^1(\Omega_b)$ to $(L^2_0(\R))^4$, as it is the composition of the bounded operator $\mathbb{D}_\theta$ on $(L^2_0(\R))^4$ and the compact map $u_b\longmapsto \Phi(u_b)$. (To see the compactness of this last map, note that it can be thought of as a composition of the bounded trace map from $H^1(\Omega_b)$ to $H^{1/2}(\Sigma_b)$ and the compact embedding $H^{1/2}(\Sigma_b)\subset L^2(\Sigma_b)$.) The compactness of the third term is a consequence of the Rellich embedding theorem, that the embedding $H^1(\Omega_b)\subset L^2(\Omega_b)$ is compact. To see that the last term is compact one can use Lemma \ref{lem-DtNcompact}. Indeed, note that continuity of $\Lambda_\theta$,  combined with compactness of the map $v_b\longmapsto v_b|_{\Sigma^j_b}$  from $ H^1(\Omega_b)$ to $L^2(\Sigma^j_b)$, suffices to conclude.

Since the sesquilinear form is coercive plus compact, the Fredholm alternative holds {(e.g.\ \cite[Theorem 2.33]{McLean})}, so that   	
%Finally, the Fredholm alternative holds and the
{unique solvability of (\ref{eq:HSMM_casgeneral-variational}) and the stability bound \eqref{eq:stabilityGen} are} a consequence of Proposition \ref{lem-uniqueness-generalcase}. {We have shown the last sentence of the theorem in our derivation, earlier in this section, of (\ref{eq:HSMM_casgeneral-variational}) from \eqref{pb:probleme_diff_initial}-\eqref{eq:Sommerfeld}, using Theorem \ref{th:deformed_representation}.}
\end{proof}
% \section{The Complexified HSMM method for case b}
% This to include
% \begin{itemize}
% \item Shifted complexification (the French team are keen to do the extension to the case where the complexification, deformation of the path of integration, starts away from the corner, which seems particularly relevant in this case).
% \item Formulation of the coupled system (Inhomogeneous Helmholtz with variable coefficients in bounded domain, plus complex-scaled HSMM equations)
% \item its Fredholm property, in fact that operator is coercive + compact
% \end{itemize}

\section{Numerical implementation and results}\label{sec:numerical}
%Let us remind that Half-Space Matching formulations have been introduced for the numerical approximation of scattering problems with complex backgrounds, when usual methods are not applicable.
In this section we demonstrate, through some illustrative numerical experiments {implemented in XLiFE++ \cite{KielLune17}}, that the complex-scaled HSM formulations  \eqref{eq:HSMM_matrice_cmplx} and \eqref{eq:HSMM_casgeneral-variational} can be solved numerically to compute solutions to the scattering problems \eqref{pb:probleme_Dir_real}-\eqref{eq:Sommerfeld} and \eqref{pb:probleme_diff_initial}-\eqref{eq:Sommerfeld}, respectively.

\subsection{Numerical implementation of the deformed half-space representation}
Before considering the discretization of the HSM systems, we just want to provide an illustration of Theorem  \ref{th:deformed_representation}. More precisely, let
\begin{equation}\label{eq:Hank_validation}
	u({\bsx}):=\frac{\ri}{4}H_0^{(1)}(k\,R(x_1,x_2)),\quad \bsx\in\R^2\setminus\{0\},
\end{equation}
where $R$ is defined by \eqref{eq:R(s,t)}, so that $u$ satisfies \eqref{pb:probleme_Dir_real}-\eqref{eq:Sommerfeld} in the case that $g:=u|_{\Sigma_a}$. Then $\varphi^0:=u|_{\Sigma^0}$ is given by
$$\varphi^0(s)=\frac{\ri}{4}H_0^{(1)}(k\,R(s,a)), \quad s\in \R.$$
The corresponding complex-scaled trace $\varphi^0_\theta$, defined by \eqref{eq:parametrage_t_theta} and \eqref{eq:definition_phi_theta}, is an even function given, for $|s|< a$, by
\[
	\varphi^0_\theta(s):=\varphi^0(s),
\]
and, for $s> a$, by
\begin{eqnarray} \nonumber
	\varphi^0_\theta(s):= \varphi^0(\tau_\theta(s)) & =\frac{\ri}{4}H_0^{(1)}\left(k\,\sqrt{a^2+(a+(s-a)\re^{\ri \theta})^2}\,\right)\\ \label{eq:asympCT}
& \sim \frac{\re^{\ri({ka}+\pi/4-\theta/2)}\,\re^{\ri k (s{-a})\re^{\ri \theta}}}{2\sqrt{2\pi ks}}\quad \text{as }\;s\rightarrow +\infty,
\end{eqnarray}
by \eqref{eq:HankAsym} {and \eqref{eq:Rasymp1}}. The asymptotic behaviour  \eqref{eq:asympCT} agrees with \eqref{eq:varphiAsymp} and Proposition \ref{cor:varphiL2}, indeed demonstrates that \eqref{eq:varphiAsymp} and Proposition \ref{cor:varphiL2} are sharp.

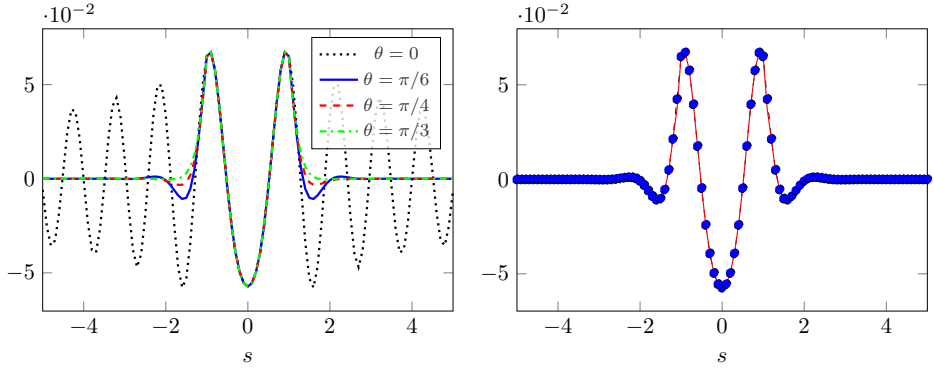
\begin{figure}
	\centering
\begin{tikzpicture}[scale=0.85]
	% axis coulb be changed in log axis or loglogaxis	
	\begin{groupplot}[group style={group size=2 by 1}, width=0.5\textwidth]
	    \nextgroupplot[
			     % ==================
			     % Main options
			     % ==================
				 width= 8 cm,
				 height= 6 cm,
				  xmin = -5,
 				  xmax = 5,
				 line width = 0.25pt, % axis line width
				 % ==================
			     % Label options
			     % ==================
				 xlabel = $s$,
				 %ylabel = $e$,
				 label style = {font=\normalsize},
				 % ==================
				 % Title options
				 % ==================
				 %title = $s \mapsto \varphi^0_\theta(s) $,
				 % ==================
				 % Legend options
				 % ==================
				 legend style = {
				 				 font=\footnotesize,
				 				 legend pos = north east, % automatic legend positioning: <north east>, <outer south west>
							      legend style={
								 			    fill opacity=0.8, % transparent legend box
  								 			   % draw=none,    % remove box surrounding the legend
  								 			   },
								 line width = .2pt, % line that boxes the legend ,
								},
				 ]
 		\addplot[
				 % ==================
				 % Line options
				 % ==================
 				black,
 				 line width = 1 pt,
				 style = dotted, % <solid>, <dashed>, <dotted>, <dashdotted>, more in documentation p. 181
				 % ==================
				 % Mark options
				 % ==================
 				 mark = none, % <none>, <> , more in documentation p. 176
  				 mark options={				 				
							   fill = black,
							   dotted, % type of line surrounding the mark: <solid>, <dashed>, ...
							   line width = .2 pt,
							  },
  				 % only marks,
				 ]
		 		table[
 		x index = 0,
 		y index = 1]
		{solexact0.txt};
		;
 		\addlegendentry{$\theta=0$}
 		\addplot[
				 % ==================
				 % Line options
				 % ==================
 				 blue,
 				 line width = 1 pt,
				 style = solid, % <solid>, <dashed>, <dotted>, <dashdotted>, more in documentation p. 181
				 % ==================
				 % Mark options
				 % ==================
 				 mark = none, % <none>, <> , more in documentation p. 176
  				 mark options={				 				
							   fill = blue,
							   solid, % type of line surrounding the mark: <solid>, <dashed>, ...
							   line width = .2 pt,
							  },
  				 % only marks,
				 ]
		 		table[
 		x index = 0,
 		y index = 1]
		{solexact30.txt};
		;
 		\addlegendentry{$\theta=\pi/6$} % legend can be defined, among other things, be <\addlegendentry{XXX}> after each plot or by <\legend{XXX, YYY}> all at once
 		\addplot[
				 % ==================
				 % Line options
				 % ==================
 				 red,
 				 line width = 1 pt,
				 style = dashed, % <solid>, <dashed>, <dotted>, <dashdotted>, more in documentation p. 181
				 % ==================
				 % Mark options
				 % ==================
 				 mark = none, % <none>, <> , more in documentation p. 176
  				 mark options={				 				
							   fill = blue,
							   solid, % type of line surrounding the mark: <solid>, <dashed>, ...
							   line width = .2 pt,
							  },
  				 % only marks,
				 ]
		 		table[
 		x index = 0,
 		y index = 1]
		{solexact45.txt};
		;
 		\addlegendentry{$\theta=\pi/4$}	 % legend can be defined, among other things, be <\addlegendentry{XXX}> after each plot or by <\legend{XXX, YYY}> all at once
 		\addplot[
				 % ==================
				 % Line options
				 % ==================
 				 green,
 				 line width = 1 pt,
				 style = dashdotted, % <solid>, <dashed>, <dotted>, <dashdotted>, more in documentation p. 181
				 % ==================
				 % Mark options
				 % ==================
 				 mark = none, % <none>, <> , more in documentation p. 176
  				 mark options={				 				
							   fill = blue,
							   solid, % type of line surrounding the mark: <solid>, <dashed>, ...
							   line width = .2 pt,
							  },
  				 % only marks,
				 ]
		 		table[
 		x index = 0,
 		y index = 1]
		{solexact60.txt};
		;
 		\addlegendentry{$\theta=\pi/3$}
			    \nextgroupplot[
			     % ==================
			     % Main options
			     % ==================
				 width= 8 cm,
				 height= 6 cm,
				  xmin = -5,
		 				  xmax = 5,
				 line width = 0.25pt, % axis line width
				 % ==================
			     % Label options
			     % ==================
				 xlabel = $s$,
				 %ylabel = $e$,
				 label style = {font=\normalsize},
				 % ==================
				 % Title options
				 % ==================
				 %title = $s \mapsto \varphi^0_\theta(s) $ and $s \mapsto \varphi^0_{\theta,{\bf h}}(s) $ for $\theta=\pi/6$,
				 % ==================
				 % Legend options
				 % ==================
				 legend style = {
				 				 font=\footnotesize,
				 				 legend pos = north east, % automatic legend positioning: <north east>, <outer south west>
							      legend style={
								 			    fill=none, % transparent legend box
		  								 			   % draw=none,    % remove box surrounding the legend
		  								 			   },
								 line width = .2pt, % line that boxes the legend ,
								},
				 ]
		
		 		\addplot[
				 % ==================
				 % Line options
				 % ==================
		 				 black,
		 				 line width = 0 pt,
				 style = dashdotted, % <solid>, <dashed>, <dotted>, <dashdotted>, more in documentation p. 181
				 % ==================
				 % Mark options
				 % ==================
		 				 mark = *, % <none>, <> , more in documentation p. 176
		  				 mark options={
							   fill =blue,
							   %solid, % type of line surrounding the mark: <solid>, <dashed>, ...
							   line width = .2 pt,
							  },
		  				 % only marks,
				 ]
		 		table[
		 		x index = 0,
		 		y index = 1]
		{solexact_numerique_30.txt};
		;
		 		%\addlegendentry{computed solution} % legend can be defined, among other things, be <\addlegendentry{XXX}> after each plot or by <\legend{XXX, YYY}> all at once
		 		\addplot[
				 % ==================
				 % Line options
				 % ==================
		 				 red,
		 				 line width = .5 pt,
				 style = solid, % <solid>, <dashed>, <dotted>, <dashdotted>, more in documentation p. 181
				 % ==================
				 % Mark options
				 % ==================
		 				 mark = none, % <none>, <> , more in documentation p. 176
		  				 mark options={
							   fill = blue,
							   solid, % type of line surrounding the mark: <solid>, <dashed>, ...
							   line width = .2 pt,
							  },
		  				 % only marks,
				 ]
		 		table[
		 		x index = 0,
		 		y index = 3]
		{solexact_numerique_30.txt};
		;
		 		%\addlegendentry{exact solution}
		\end{groupplot}
	\end{tikzpicture}
	\caption{Left:  representation of the exact complex-scaled traces $s\mapsto\varphi^0_\theta(s)$ for $\theta=0,\pi/6,\pi/4,\pi/3$. Right: comparison of the exact (red line) and computed (blue dots) complex-scaled trace $s\mapsto\varphi^0_\theta(s)$ for $\theta=\pi/6$. In both plots $a=1$ and $k=2\pi$.}
	\label{fig:varphi_theta}
	\end{figure}

\begin{figure}[htbp]
		\begin{center}
			\begin{tikzpicture}
		\node[inner sep=0pt] (source) at (0,0)
		    {\includegraphics[width=13cm]{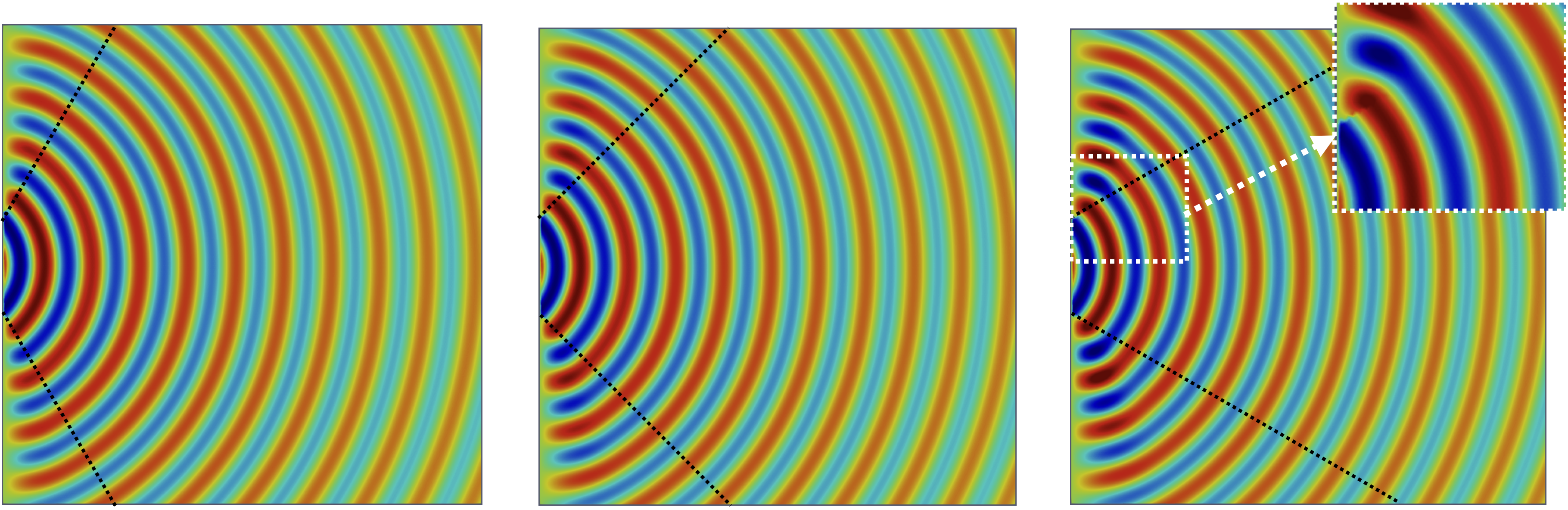}};
			\draw(-4.2,-2.3) node {$\theta=\pi/6$};
			\draw(0.3,-2.3) node {$\theta=\pi/4$};
			\draw(4.8,-2.3) node {$\theta=\pi/3$};
			% \draw[black,dashed] (-5.77,-2)--(-5.77,2);
% 			\draw[black,dashed] (-1.55,-2)--(-1.55,2);
% 			\draw[black,dashed] (2.74,-2)--(2.74,2);
			\end{tikzpicture}
	\caption{Representation of $U_\theta(\varphi^0_\theta)$ in $\Omega^0$ (see  \eqref{eq:U_theta_num}), where $\varphi^0_\theta(s)=\dfrac{\ri}{4}H^{(1)}_0(k\, R(\tau_\theta(s),a))$, with $a=1$ and $k=2\pi$. The dotted lines are part of the boundary of $\Omega^0_\theta$.}
	\label{fig:half_rep}
	\end{center}
	\end{figure}
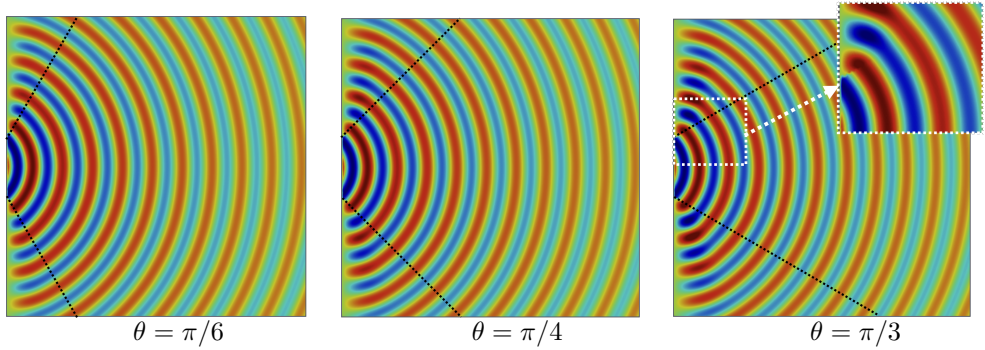

We
 represent, on Figure \ref{fig:varphi_theta}, $\varphi^0_\theta$ for four different values of $\theta$, with $a=1$ and $k=2\pi$.
We see that $\varphi^0_\theta$ is more and more rapidly decaying at infinity as $\theta$ increases, in line with \eqref{eq:asympCT}, \eqref{eq:varphiAsymp}, and Proposition \ref{cor:varphiL2}. Then we represent in Figure~\ref{fig:half_rep}  the function
\begin{equation}\label{eq:U_theta_num}
	 U_\theta(\varphi^0_\theta)({\bsx}^0) = \int_{\R}h(x_1^0-a,x_2^0-\tau_\theta(s))\,\varphi^0_\theta(s)\,\tau'_\theta(s)\, \rd s
	 \end{equation}
(cf.\ \eqref{eq:expr_Ujtheta}) in the half-space $\Omega^0$,  for  $\theta=\pi/6, \pi/4, \pi/3$, evaluating this integral accurately by standard numerical quadrature methods (namely a 5th order composite Gauss quadrature rule on a fine mesh of step length $0.1$). From Theorem  \ref{th:deformed_representation} we know that
\begin{equation}
	\label{Utheta=H10}
 U_\theta(\varphi^0_\theta)({\bsx}^0) = u(\bsx^0) := \frac{\ri}{4}H_0^{(1)}(k\, R(x_1^0,x_2^0)),\quad {\bsx}^0\in \Omega^0_\theta,
\end{equation}
where $\Omega^0_\theta:=\{{\bsx}^{0}=(x_1^0,x_2^0): \;x_1^0-a>(|x_2^0|-a)\tan\theta\}$. The boundary of $\Omega^0_\theta$ is indicated on the figures by dashed lines and, as predicted in \eqref{Utheta=H10}, $U_\theta(\varphi^0_\theta)$ coincides with $u$ in $\Omega^0_\theta$. $U_\theta(\varphi^0_\theta)$ is not equal to $u$ outside $\Omega^0_\theta$; in particular, it is easy to see from the definition that $U_\theta(\varphi^0_\theta)(a,x_2^0)=0$ for $|x_2^0|>a$. It appears at first glance   that $U_\theta(\varphi^0_\theta)$ is continuous across the dashed lines in Figure~\ref{fig:half_rep}, but an  application of the residue theorem, modifying the argument of Theorem \ref{th:deformed_representation}, shows a jump in the value of $U_\theta(\varphi^0_\theta)$ across the dashed lines of
\[
	\varphi^0(x_2^0\pm\ri (x^0_1-a))=\varphi^0_\theta(s)%\;\;\text{where}}\;\; x_2^0 = \pm(a+(x_1^0-a)\cot(\theta)= {\pm a +(s\mp a)\cos(\theta).}
\]
at the point $\bsx^0=(x_1^0,x_2^0)$ where $x_2^0 = \pm(a+(x_1^0-a)\cot(\theta))= \pm a +(s\mp a)\cos(\theta)$.
%which is nothing else than $
%	\varphi^0_\theta(s)\;\text{if}\; x_1^0-a=\pm s\sin\theta.
%$
This jump across the dashed lines is just about visible very close to $\bsx^0=(a,\pm a)$, but not visible elsewhere %We can hardly see this jump
because $\varphi^0_\theta(s)$ is exponentially decaying as $s\to \pm \infty$; see the zoom for $\theta=\pi/3 $.
\subsection{Discretization of the complex-scaled HSM formulation and validations} \label{sec:NUM1}
%Let us briefly describe how to discretize this new formulation and show some numerical results to illustrate and validate the method.\\\\
We approximate the solution ${\Phi}_{\theta}-{\Phi}_{g}$ of \eqref{eq:HSMM_matrice_cmplx} by $\widetilde{\Phi}_{\theta,{\bf h}}\in {\bf V}_{0,\bf{h}}$, where ${\bf V}_{0,\bf{h}}$ is a finite dimensional space  ${\bf V}_{0,\bf{h}}\subset(L^2_0(\R))^4$ that we specify below and $\widetilde{\Phi}_{\theta,{\bf h}}$ is the solution of the following Galerkin approximation:
\begin{equation}\label{eq:HSMM_discrete}
	\begin{array}{c}
		\text{Find}\;  \widetilde{\Phi}_{\theta,{\bf h}}\in {\bf V}_{0,\bf{h}} \text{ such that }
		\\[5pt]
		((\mathbb{I}-\mathbb{D}_\theta) \widetilde{\Phi}_{\theta,{\bf h}}, \widetilde{\Psi}_{\bf h})_{(L^2_0(\R))^4}
		= 	( \mathbb{D}_\theta{\Phi}_{g,{\bf h}}, \widetilde{\Psi}_{\bf h})_{(L^2_0(\R))^4},\quad \forall   \widetilde{\Psi}_{\bf h}\in {\bf V}_{0,\bf{h}},
	\end{array}
\end{equation}
where $\mathbb{D}_\theta$ is defined by (\ref{eq:def_Dmat_cmplx}).

To define the approximation space ${\bf V}_{0,\bf h}$, where ${\bf h}:=(h,q,T)$, let us first introduce $V_{\bf h}\subset L^2(\R)$. To construct $V_{\bf h}$ we truncate the infinite line at some distance $T>0$ and build $V_{\bf h}$ with 1D Lagrange finite elements of degree $q\geq 1$ and maximum element length
 $h$ supported on $[-T,T]$. The space ${\bf V}_{0,\bf h}$ is nothing else but $ {\bf V}_{\bf h}\cap (L^2_0(\R))^4$ where ${\bf V}_{\bf h}:=(V_{\bf h})^4$. Thus, for $\widetilde{\Phi}_{\theta,{\bf h}}=
\{\tilde\varphi^0_{\theta,{\bf h}},\tilde\varphi^1_{\theta,{\bf h}},\tilde\varphi^2_{\theta,{\bf
h}},\tilde\varphi^3_{\theta,{\bf h}}\} \in {\bf V}_{0,{\bf h}}$, each $\tilde\varphi^j_{\theta,{\bf h}}$ is
a continuous piecewise polynomial function supported in $[-T,T]$ which vanishes on $[-a,a]$.
In
\eqref{eq:HSMM_discrete}, ${\Phi}_{g,{\bf h}}\in {\bf V}_{\bf h}\cap (L^2(-a,a))^4$ is an interpolate of
${\Phi}_{g}$. Finally, we approximate ${\Phi}_{\theta}$ by ${\Phi}_{\theta,{\bf h}}=
\{ \varphi^0_{\theta,{\bf h}},\varphi^1_{\theta,{\bf h}},\varphi^2_{\theta,{\bf
h}},\varphi^3_{\theta,{\bf h}} \} \in {\bf V}_{\bf h}$,
given by
${\Phi}_{\theta,{\bf h}}:=\widetilde{\Phi}_{\theta,{\bf h}}+{\Phi}_{g,{\bf h}}$.

It is clear that the approximation space ${\bf V}_{0,\bf h}$ that we have constructed has the approximation property that, for all $\Phi\in (L_0^2(\R))^4$,
$$
\inf_{\widetilde \Psi_{\bf h}\in {\bf V}_{0,\bf h}} \|\Phi-\widetilde \Psi_{\bf h}\|_{(L_0^2(\R))^4} \to 0
$$
as $h\to 0$ and $T\to +\infty$. Thus, and since the sesquilinear form in \eqref{eq:HSMM_discrete} is coercive plus compact on $(L_0^2(\R))^4$ by Theorem \ref{th:HSMM_cmplxed}(ii) (cf.\ Theorem \ref{the-generalcase}), standard convergence results for Galerkin methods apply (e.g., \cite[Theorem 4.2.9]{sauter-schwab11}). These give that, for some $h_0>0$ and $T_0>0$, the solution $\widetilde{\Phi}_{\theta,{\bf h}}$ of \eqref{eq:HSMM_discrete} is well-defined for all $0<h\leq h_0$ and $T\geq T_0$, and a quasi-optimality error estimate holds, that, for some constant $C>0$ and all $0<h\leq h_0$ and $T\geq T_0$,{
\begin{eqnarray} \nonumber
\|\Phi_\theta-{\Phi}_{\theta,{\bf h}}\|_{(L^2(\R))^4} &\leq & C\inf_{\Psi_{\bf h}\in \Phi_{g, \bf h}+{\bf V}_{0,\bf h}} \| \Phi_{\theta}-\Psi_{\bf h}\|_{(L^2(\R))^4}\\ \nonumber
  &{\leq} & C\Big(\|\Phi_{\theta}\|_{(L^2(\R\setminus(-T,T)))^4}\\ \label{eq:QuasiOpt}
  & & \hspace{12ex} + \; \inf_{\Psi_{\bf h}\in \Phi_{g, \bf h}+{\bf V}_{0,\bf h}} \| \Phi_{\theta}-\Psi_{\bf h}\|_{(L^2(-T,T))^4}\Big).
\end{eqnarray}
This right hand side tends to zero} as $h\to 0$ and $T\to+\infty$, i.e., our Galerkin method is convergent, as long as $\| \Phi_{g}-\Phi_{g,\bf h}\|_{(L^2(-a,a))^4} \rightarrow 0$ as $h\to 0$.
%\begin{eqnarray*}
%\|(\Phi_\theta-\Phi_g)-\widetilde{\Phi}_{\theta,{\bf h}}\|_{(L^2_0(\R))^4} &\leq & C\Big(\inf_{\widetilde \Psi_{\bf h}\in {\bf V}_{0,\bf h}} \| (\Phi_{\theta}-\Phi_g)-\widetilde \Psi_{\bf h}\|_{(L_0^2(\R))^4}\\
%& & \hspace{10ex}  + \; \|\Phi_g-\Phi_{g,\bf h}\|_{(L^2(-a,a))^4}\Big) \to 0
%\end{eqnarray*}
%\begin{eqnarray*}
%\|\Phi_\theta-{\Phi}_{\theta,{\bf h}}\|_{(L^2(\R))^4} &\leq & C\Big(\inf_{\Psi_{\bf h}\in \Phi_g+{\bf V}_{0,\bf h}} \| \Phi_{\theta}-\Psi_{\bf h}\|_{(L^2(\R))^4}\\
%& & \hspace{10ex}  + \; \|\Phi_g-\Phi_{g,\bf h}\|_{(L^2(-a,a))^4}\Big) \to 0
%\end{eqnarray*}

{To implement the Galerkin method \eqref{eq:HSMM_discrete}, the}  integrals
$$
 \int_a^{T}\left(\int_{-T}^Th(\tau_\theta(t)-a,a-\tau_\theta(s))\,\varphi^j_{\theta,{\bf h}}(s)\,\tau'_\theta(s)\, \rd s\right)\tilde\psi^{j\pm 1}_{\bf h}(t)\rd t,\quad \tilde\psi^{j\pm 1}_{\bf h}\in V_{\bf h},
 $$
which appear in the variational formulation, {need to be approximated. In the results below we use a standard quadrature formula, without any specific treatment of the singularity.}

To validate the method, we consider $u$ given by \eqref {eq:Hank_validation} which is the solution of \eqref{eq:HSMM_matrice_cmplx}-\eqref{eq:Sommerfeld} with $g:=u|_{\Sigma_a}$. % the particular data  the Hankel function on the % boundary of the square:
% \[
% 	g(\bsx) = \frac{\ri}{4}H_0^{(1)}(k R(x_1,x_2)),
% \]
% so  that the exact solution of  \eqref{eq:HSMM_matrice_cmplx}-\eqref{eq:Sommerfeld} is given by $$u(\bsx) = \frac{\ri}{4}H_0^{(1)}(k R(x_1,x_2)).$$
{In the case that $\theta=\pi/6$, $a=1$, and $k=2\pi$,} we draw in the right hand side of Figure \ref{fig:varphi_theta} the exact complex-scaled trace $\varphi^0_{\theta}$ of $u$ (by symmetry, the four traces are equal in this case) and {the computed complex-scaled trace $\varphi^0_{\theta,{\bf h}}$, obtained by solving \eqref{eq:HSMM_discrete} with $h=0.1$, $q=1$, and $T=5$}. We observe a very good agreement.
%		\begin{figure}[htbp]
%		\begin{center}
%			\begin{tikzpicture}
%		\node[inner sep=0pt] (source) at (0,0)
%		    {\includegraphics[width=8cm]{real30.png}};
%			\end{tikzpicture}
%	\caption{The computed complex-scaled trace $\varphi^0_{\theta,{\bf h}}$ obtained by solving \eqref{eq:HSMM_discrete} and the exact complex scaled trace $\varphi^0_{\theta}$ for $\theta=\pi/6$.}
%	\label{fig:trace30}
%	\end{center}
%	\end{figure}
% Then, we represent in Figure \ref{fig:half_rep}  the reconstructed solution in a part of $\Omega^0_\theta$ for $\theta=\pi/6, \pi/4, \pi/3$ defined by \eqref{eq:Omegaj_theta}, by using formula \eqref{eq-repinOmegatheta}. Note that larger is $\theta$, smaller is $\Omega^0_\theta$.

\begin{figure*}[htbp]
		\centering
		\begin{tikzpicture}
			\begin{groupplot}[group style={group size=2 by 1}, width=0.5\textwidth]
			    \nextgroupplot[
		xlabel={$T$},
		ylabel={$\log\|\varphi^0_\theta - \varphi^0_{\theta, {\bf h}}\|_{L^2(\R)}$},
		xmin=1, xmax=3.2,
		ymin=-3.5, ymax=-1,
		legend pos=north east,
		ymajorgrids=true,
		grid style=dashed,
		axis lines = left,
		]
		%theta = pi/6
		\addplot[color=blue, mark=square]coordinates{
			(1.2, -1.2508808) (1.4, -1.3941059)(1.6, -1.54154346) (1.8, -1.6914057) (2, -1.842516) (2.2, -1.9940049) (2.4, -2.1451319) (2.6, -2.2950948) (2.8, -2.4428563) (3, -2.586546)
		};		
		%theta = pi/4
		\addplot[color=red, mark=triangle]coordinates{
			(1.2, -1.3475187) (1.4, -1.5357236)(1.6, -1.7314478) (1.8, -1.9315395) (2, -2.132129) (2.2, -2.3270068) (2.4, -2.5050777) (2.6, -2.649015856) (2.8, -2.745638054) (3, -2.797220253)
		};		
		%theta = pi/3
		\addplot[color=green, mark=diamond]coordinates{
			(1.2, -1.402283977) (1.4, -1.616836565)(1.6, -1.838081242) (1.8, -2.054323181) (2, -2.241640785) (2.2, -2.370061618) (2.4, -2.434585512) (2.6, -2.459083936) (2.8, -2.467576881) (3, -2.469973285)
		};

 		 \addplot[color=green, mark=diamond]coordinates{(1.8, -1.3)(2,-1.53)(2.2,-1.76 )(2.4, -2.0090)};
		 \addplot[color=green]coordinates{(2.4, -2.0090)(2.4, -1.3)(1.8,-1.3)};
		 \addplot[color=red, mark=triangle]coordinates{(1.8, -1.3)(2,-1.5)(2.2,-1.7 )(2.4, -1.8789)};
		  \addplot[color=red]coordinates{(2.4, -1.8789)(2.4, -1.3)(1.8,-1.3)};
		  \addplot[color=blue, mark=square]coordinates{(1.8, -1.3)(2,-1.43)(2.2,-1.57)(2.4, -1.7093)};
		 \addplot[color=blue]coordinates{(2.4, -1.7093)(2.4, -1.3)(1.8,-1.3)};
 		\node[] at (axis cs: 2.15,-1.2) {\tiny $1$};
 		\node[color=blue] at (axis cs: 2.67,-1.4) {\tiny $k\sin(\pi/6)$};
		\node[color=red] at (axis cs: 2.67,-1.8) {\tiny $k\sin(\pi/4)$};
		\node[color=green] at (axis cs: 2.67,-1.95) {\tiny $k\sin(\pi/3)$};
				\nextgroupplot[
						xlabel={$T$},
						xmin=1, xmax=3.2,
						ymin=-3.5, ymax=-1,
						legend pos=north east,
						ymajorgrids=true,
						grid style=dashed,
						axis lines = left,
						]
		
						%theta = pi/6
						\addplot[color=blue, mark=square]coordinates{
						(1.2, {log10(0.023203})
						(1.4, {log10(0.0126916)})
						(1.6, {log10(0.006770525)})
						(1.8, {log10(0.003560175)})
						(2, {log10(0.00186891)})
						(2.2, {log10(0.001005495)})
						(2.4, {log10(0.0005912975)})
						(2.6, {log10(0.0004184925)})
						% (2.8, {log10(0.000488785)})
% 						(3, {log10(0.00047697)})
						};
		
						%theta = pi/4
						\addplot[color=red, mark=triangle]coordinates{
						(1.2, {log10(0.01672115)})
						(1.4, {log10(0.00747425)})
						(1.6, {log10(0.00325355)})
						(1.8, {log10(0.001499185)})
						(2, {log10(0.000901445)})
						(2.2, {log10(0.000759355)})
						(2.4, {log10(0.0007336725)})
						(2.6, {log10(0.0007300825)})
						% (2.8, {log10(0.00102641)})
% 						(3, {log10(0.00102731)})
						};
		
						%theta = pi/3
						\addplot[color=green, mark=diamond]coordinates{
							(1.2, {log10(0.01371945)}) (1.4, {log10(0.005564525)}) (1.6, {log10(0.0025784)}) (1.8, {log10(0.0018241825)}) (2, {log10(0.00170818)}) (2.2, {log10(0.0016944025)}) (2.4, {log10(0.00169234)}) (2.6, {log10(0.0016927425)})
							% (2.8, {log10(0.00238169)}) (3, {log10(0.0023829)})
						};
						\legend{$\theta = \pi/6$, $\theta = \pi/4$, $\theta = \pi/3$};
		
				 		 \addplot[color=green, mark=diamond]coordinates{(1.5, -1.3)(1.7,-1.76)(1.9,-2.22)(2.1, -2.71)};
						  \addplot[color=green]coordinates{(2.1, -2.71)(2.1, -1.3)(1.5,-1.3)};
						  \addplot[color=red, mark=triangle]coordinates{(1.5, -1.3)(1.7,-1.71)(1.9,-2.12)(2.1, -2.4577)};
						  \addplot[color=red]coordinates{(2.1, -2.4577)(2.1, -1.3)(1.5,-1.3)};
						  \addplot[color=blue, mark=square]coordinates{(1.5, -1.3)(1.7,-1.58)(1.9,-1.86)(2.1, -2.1186)};
						 \addplot[color=blue]coordinates{(2.1, -2.1186)(2.1, -1.3)(1.5,-1.3)};
				 		% \node[] at (axis cs: 2.15,-1.2) {\tiny $1$};
% 				 		\node[color=blue] at (axis cs: 2.67,-1.4) {\tiny $k\sin(\pi/6)$};
% 						\node[color=red] at (axis cs: 2.67,-1.8) {\tiny $k\sin(\pi/4)$};
% 						\node[color=green] at (axis cs: 2.67,-1.95) {\tiny $k\sin(\pi/3)$};
				\end{groupplot}
		\end{tikzpicture}
		\caption{Absolute error in $\varphi^0_{\theta,{ \bf h}}$ for three different values of  $\theta$ computed with $P3$ elements ($q=3$), $a=1$, and $h=0.002$ and $k = \pi$ (left), $h=0.001$ and $k=2\pi$ (right).}
		\label{fig:abs_error_kpi}
	\end{figure*}
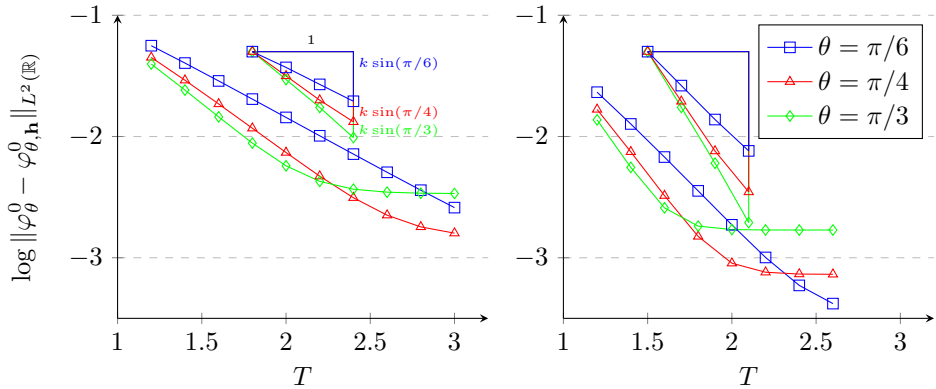

More quantitatively, {to explore the dependence of the error on $T$}, we plot in Figure \ref{fig:abs_error_kpi}, {for a fixed small value of $h$,} the error
\begin{equation} \label{eq:errorEq}
\|\varphi^0_\theta-\varphi^0_{\theta,{\bf h}}\|_{L^2(\R)}
\end{equation}
as a function of $T$ for $\theta=\pi/6, \pi/4, \pi/3$. As $\varphi^0_{\theta,{\bf h}}$ is zero outside $(-T,T)$ the error cannot be smaller than the $L^2$ norm of $\varphi^0_\theta$ on $\R\setminus(-T,T)$, which decreases like $\re^{-k\sin (\theta) T}/T^{1/2} $ as $T$ tends to $+\infty$ by \eqref{eq:asympCT} and Proposition \ref{cor:varphiL2}. {On the other hand, when $h$ is small enough so that the second term on the right hand side of \eqref{eq:QuasiOpt} is negligible, the quasi-optimality bound  \eqref{eq:QuasiOpt} and Proposition \ref{cor:varphiL2} imply that {$\re^{-k\sin (\theta) T}/T^{1/2}$} is also an upper bound for the error, {precisely} that the error \eqref{eq:errorEq} is $\leq C \re^{-k\sin (\theta) T}/T^{1/2} $, for some constant $C>0$ and all sufficiently large $T$.  And, indeed,} we observe in Figure \ref{fig:abs_error_kpi} this rate of exponential behaviour as $T$ increases, until the other sources of error become significant. %{and the graph reaches a plateau.}  }

\begin{figure}[htbp]
		\begin{center}
			\begin{tikzpicture}
		\node[inner sep=0pt] (source) at (0,0)
		    {\includegraphics[width=13cm,height=6cm]{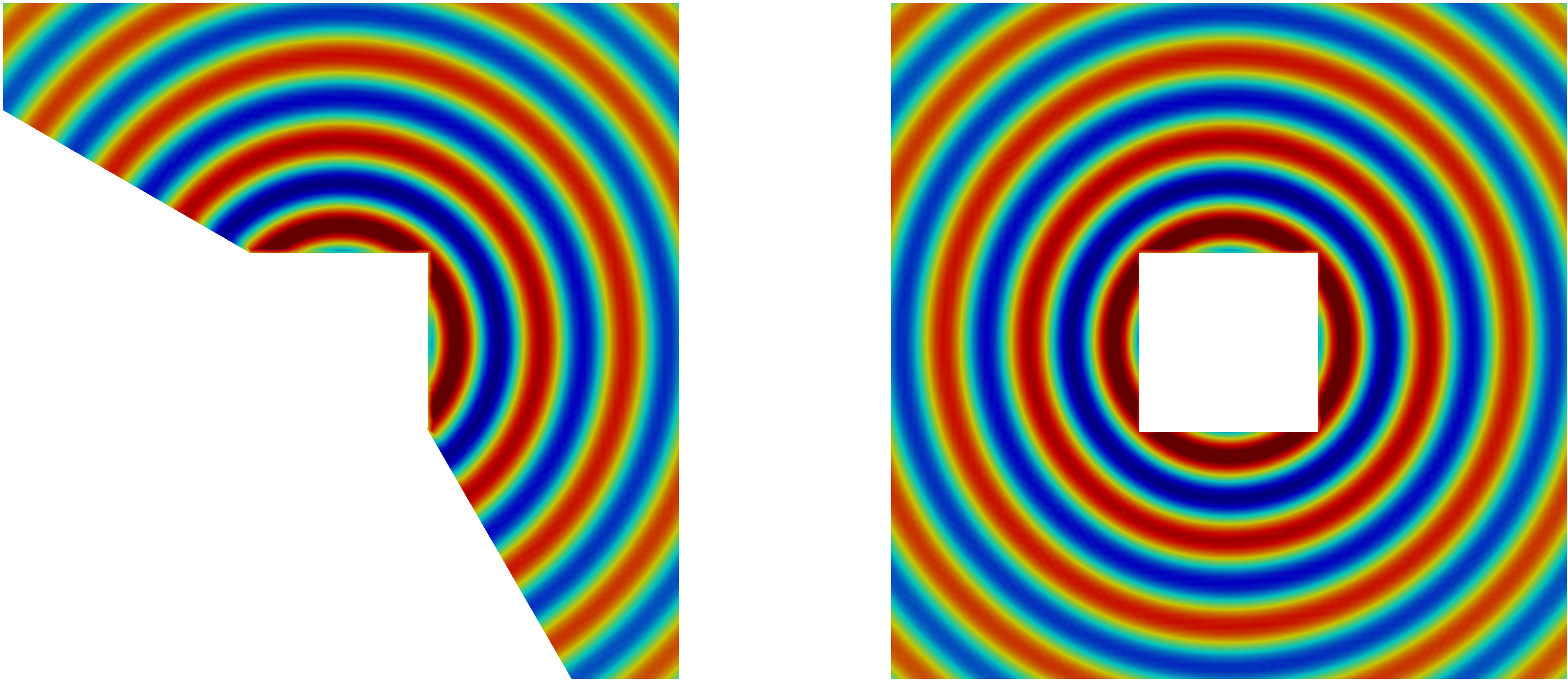}};
			\draw(-3.5,-3.2) node {$u\big|_{\Omega^0_\theta\cup \Omega^1_\theta}$};
			\draw(4.5,-3.2) node {$u\big|_{\Omega}$};
			\end{tikzpicture}
	\caption{Reconstruction of the computed solution in $\Omega^0_\theta\cup \Omega^1_\theta$ (left) and in the whole domain~(right), with $\theta=\pi/6$, $a=1$, $k=2\pi$.}
	\label{fig:rep_sol_theta30}
	\end{center}
	\end{figure}

Using a standard numerical quadrature applied to formula \eqref{eq-repinOmegatheta}, with $\varphi_\theta^j$ approximated by  $\varphi^0_{\theta,{ \bf h}}$, we represent finally in Figure~\ref{fig:rep_sol_theta30} the numerical solution in $\Omega^0_\theta\cup \Omega^1_\theta$ (left) and the whole of $\Omega$ (right) for $\theta=\pi/6$. To reconstruct the solution in the whole of $\Omega$ {several choices are possible, since $\theta<\pi/4$ is such that the reconstruction domains overlap ($\Omega^{j}_\theta$ overlaps with $\Omega^{j\pm1}_\theta$, for $j\in \llbracket 0, 3 \rrbracket$). Here} we have reconstructed the solution using the identity
\[
	u({\bsx}^j)=U_\theta(\varphi^j_{\theta})({\bsx}^j),\quad {\bsx}^j\in \Omega^j_{\pi/4},\;\;j\in \llbracket 0, 3 \rrbracket,
\]
where $U_\theta$ is defined in \eqref{eq:expr_Ujtheta}, and where we have used that $ \Omega^j_{\pi/4}\subset  \Omega^j_{\theta}$.
We  notice in Figure~\ref{fig:rep_sol_theta30} that these different representations are compatible, up to a small discretization error not visible in the plots.%, which is another sign that the method works.

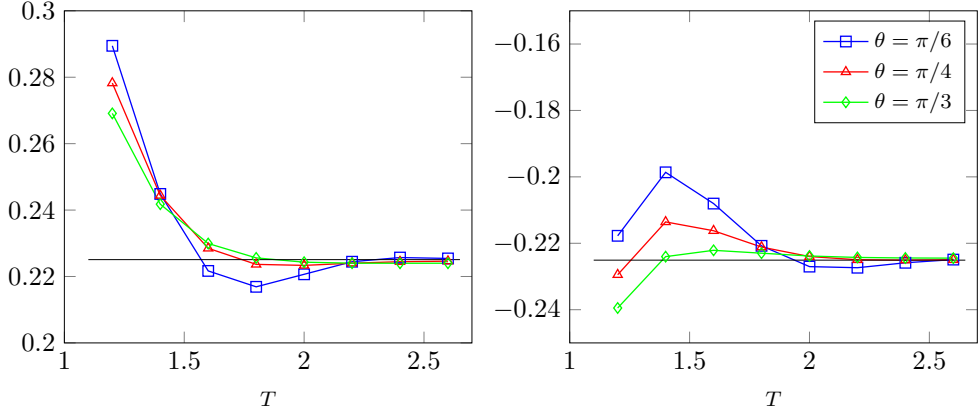
\begin{figure}[h]
			\centering
			\begin{tikzpicture}
				\begin{groupplot}[group style={group size=2 by 1}, width=0.5\textwidth]
				    \nextgroupplot[
						     % ==================
						     % Main options
						     % ==================
							 width= 7 cm,
							 height= 6 cm,
							  xmin = 1,
			 				  xmax = 2.7,
			   			      ymin = 0.2,
			   				  ymax = 0.3,
							 line width = 0.25pt, % axis line width
							 % ==================
						     % Label options
						     % ==================
							 xlabel = $T$,
							 %ylabel = $e$,
							 label style = {font=\footnotesize},
							 % ==================
							 % Title options
							 % ==================
							 %title = $s \mapsto \varphi^0_\theta(s) $,
							 % ==================
							 % Legend options
							 % ==================
							 legend style = {
							 				 font=\footnotesize,
							 				 legend pos = north east, % automatic legend positioning: <north east>, <outer south west>
										      legend style={
											 			    fill=none, % transparent legend box
			  								 			   % draw=none,    % remove box surrounding the legend
			  								 			   },
											 line width = .2pt, % line that boxes the legend ,
											},
							 ]
			 		\addplot[
							 % ==================
							 % Line options
							 % ==================
			 				blue,
			 				 line width = 0.5 pt,
							 style = solid, % <solid>, <dashed>, <dotted>, <dashdotted>, more in documentation p. 181
							 % ==================
							 % Mark options
							 % ==================
			 				 mark = square, % <none>, <> , more in documentation p. 176
			  				 mark options={				 				
										   fill = blue,
										   solid, % type of line surrounding the mark: <solid>, <dashed>, ...
										   line width = 0.5 pt,
										  },
			  				 % only marks,
							 ]
					 		table[
			 		x index = 0,
			 		y index = 1]
					{farfield30.txt};
					;
			 		\addplot[
							 % ==================
							 % Line options
							 % ==================
			 				red,
			 				 line width = 0.5 pt,
							 style = solid, % <solid>, <dashed>, <dotted>, <dashdotted>, more in documentation p. 181
							 % ==================
							 % Mark options
							 % ==================
			 				 mark = triangle, % <none>, <> , more in documentation p. 176
			  				 mark options={				 				
										   fill = blue,
										   solid, % type of line surrounding the mark: <solid>, <dashed>, ...
										   line width = 0.5 pt,
										  },
			  				 % only marks,
							 ]
					 		table[
			 		x index = 0,
			 		y index = 1]
					{farfield45.txt};
					;
			 		\addplot[
							 % ==================
							 % Line options
							 % ==================
			 				green,
			 				 line width = 0.5 pt,
							 style = solid, % <solid>, <dashed>, <dotted>, <dashdotted>, more in documentation p. 181
							 % ==================
							 % Mark options
							 % ==================
			 				 mark = diamond, % <none>, <> , more in documentation p. 176
			  				 mark options={				 				
										   fill = green,
										   solid, % type of line surrounding the mark: <solid>, <dashed>, ...
										   line width = 0.5 pt,
										  },
			  				 % only marks,
							 ]
					 		table[
			 		x index = 0,
			 		y index = 1]
					{farfield60.txt};
					;
					\addplot[color=black]coordinates{(1.1, 0.225079)(2.65, 0.225079)};
					
				    \nextgroupplot[
						     % ==================
						     % Main options
						     % ==================
							 width= 7 cm,
							 height= 6 cm,
					 		  xmin = 1,
					 			 				  xmax = 2.7,
					 			   			      ymax = -0.15,
					 			   				  ymin = -0.25,
							 line width = 0.25pt, % axis line width
							 % ==================
						     % Label options
						     % ==================
							 xlabel = $T$,
							 %ylabel = $e$,
							 label style = {font=\footnotesize},
							 % ==================
							 % Title options
							 % ==================
							 %title = $s \mapsto \varphi^0_\theta(s) $,
							 % ==================
							 % Legend options
							 % ==================
							 legend style = {
							 				 font=\footnotesize,
							 				 legend pos = north east, % automatic legend positioning: <north east>, <outer south west>
										      legend style={
											 			    fill=none, % transparent legend box
			  								 			   % draw=none,    % remove box surrounding the legend
			  								 			   },
											 line width = .2pt, % line that boxes the legend ,
											},
							 ]
			 		\addplot[
							 % ==================
							 % Line options
							 % ==================
			 				blue,
			 				 line width = 0.5 pt,
							 style = solid, % <solid>, <dashed>, <dotted>, <dashdotted>, more in documentation p. 181
							 % ==================
							 % Mark options
							 % ==================
			 				 mark = square, % <none>, <> , more in documentation p. 176
			  				 mark options={				 				
										   fill = blue,
										   solid, % type of line surrounding the mark: <solid>, <dashed>, ...
										   line width = 0.5 pt,
										  },
			  				 % only marks,
							 ]
					 		table[
			 		x index = 0,
			 		y index = 2]
					{farfield30.txt};
					;
			 		\addplot[
							 % ==================
							 % Line options
							 % ==================
			 				red,
			 				 line width = 0.5 pt,
							 style = solid, % <solid>, <dashed>, <dotted>, <dashdotted>, more in documentation p. 181
							 % ==================
							 % Mark options
							 % ==================
			 				 mark = triangle, % <none>, <> , more in documentation p. 176
			  				 mark options={				 				
										   fill = blue,
										   solid, % type of line surrounding the mark: <solid>, <dashed>, ...
										   line width = 0.5 pt,
										  },
			  				 % only marks,
							 ]
					 		table[
			 		x index = 0,
			 		y index = 2]
					{farfield45.txt};
					;
			 		\addplot[
							 % ==================
							 % Line options
							 % ==================
			 				green,
			 				 line width = 0.5 pt,
							 style = solid, % <solid>, <dashed>, <dotted>, <dashdotted>, more in documentation p. 181
							 % ==================
							 % Mark options
							 % ==================
			 				 mark = diamond, % <none>, <> , more in documentation p. 176
			  				 mark options={				 				
										   fill = green,
										   solid, % type of line surrounding the mark: <solid>, <dashed>, ...
										   line width = 0.5 pt,
										  },
			  				 % only marks,
							 ]
					 		table[
			 		x index = 0,
			 		y index = 2]
					{farfield60.txt};
					;
					\addplot[color=black]coordinates{(1.1, 0.225079)(2.65, 0.225079)};

					\legend{$\theta = \pi/6$, $\theta = \pi/4$, $\theta = \pi/3$};
					\addplot[color=black]coordinates{(1.1, -0.225079)(2.65, - 0.225079)};
					\hspace{0.3cm}
			% \begin{axis}[
% 			xlabel={$T$},
% 			ylabel={$\log|C_\infty - C_{\infty, h}|$},
% 			xmin=1.25, xmax=3.5,
% 			ymin=-3, ymax=-1,
% 			%xtick={5, 10, 15, 20, 25, 30, 35, 40},
% 			%ytick={-1.5, -2, -2.5, -3, -3.5, -4, -4.5, -5},
% 			legend pos=north east,
% 			ymajorgrids=true,
% 			grid style=dashed,
% 			axis lines = left,
% 			]
%
% 			%theta = pi/6
% 			\addplot[color=blue, mark=square]coordinates{
% 				(1.5, {log10(0.089226)}) (1.75, {log10(0.0550791)}) (2, {log10(0.0352221)}) (2.25, {log10(0.0233751)}) (2.5, {log10(0.0158714)}) (2.75, {log10(0.0107466)}) (3, {log10(0.00700471)}) (3.25, {log10(0.00418174)}) (3.5, {log10(0.00206321)}) %(3.75, {log10(0.000530392)}) (4, {log10(0.000517978)})
% 			};
%
% 			%theta = pi/4
% 			\addplot[color=red, mark=square]coordinates{
% 				(1.5, {log10(0.0679242)}) (1.75, {log10(0.0362356)}) (2, {log10(0.0199346)}) (2.25, {log10(0.0119823)}) (2.5, {log10(0.00820372)}) (2.75, {log10(0.0061939)}) (3, {log10(0.00488595)}) (3.25, {log10(0.00397273)}) (3.5, {log10(0.00339656)}) %(3.75, {log10(0.00310809)}) (4, {log10(0.00301268)})
% 			};
%
% 			%theta = pi/3
% 			\addplot[color=green, mark=square]coordinates{
% 				(1.5, {log10(0.0535831)}) (1.75, {log10(0.0234924)})(2, {log10(0.00972691)}) (2.25, {log10(0.00609228)}) (2.5, {log10(0.00668242)}) (2.75, {log10(0.00728845)}) (3, {log10(0.00748623)}) (3.25, {log10(0.0074821)}) (3.5, {log10(0.00742127)}) %(3.75, {log10(0.00736404)}) (4, {log10(0.00732558)})
% 			};
% 			\legend{$\theta = \pi/6$, $\theta = \pi/4$, $\theta = \pi/3$}
% 			\end{axis}
		\end{groupplot}
			\end{tikzpicture}
			\caption{Real (left figure) and imaginary (right figure) parts of the far-field coefficient for three different $\theta$ with $P3$ elements ($q=3$), $h=0.001$, $a=1$, and $k = 2\pi$, computed using the formula \eqref{FF_computed}. The black lines indicate the exact values.}
			\label{fig:abs_error_far-field-kpi}
		\end{figure}

Finally, we validate formula \eqref{eq-farfieldtheta} for the far-field pattern. In the present case, 	the far-field pattern is the same in all directions and one has for $j\in \llbracket0,3\rrbracket$:
\begin{equation}\label{eq:FF_calc1}
	F(\cos (j\pi/2),\sin  (j\pi/2))=\frac{1-\ri}{\sqrt{k\pi}}=\sqrt{\frac{k}{\pi}}\frac{1-\ri}2\int_\R\varphi^j_\theta(s)\tau'_\theta(s)\,\rd s.
	\end{equation}
In Figure \ref{fig:abs_error_far-field-kpi} we plot the real and imaginary part of the right hand side of \eqref{eq:FF_calc1} when $j=0$, $a=1$, and $k=2\pi$, with $\varphi^0_\theta$ approximated by $\varphi^0_{\theta,{\bf h}}$, that is we plot
\begin{equation}\label{FF_computed}
	F_{\theta,{\bf h}}:=\sqrt{\frac{k}{\pi}}\frac{1-\ri}2\int_{-T}^T\varphi^j_{\theta,{\bf h}}(s)\tau'_\theta(s)\,\rd s
	\end{equation}
as a function of $T$ for different values of $\theta$. We again observe a rapid convergence towards the exact value as $T$ increases.
% \subsection{Numerical analysis of discretisation and truncation} French team (Sonia and Yohanes were in the room!) will think about this. A key point is to analyse the truncation, to prove that the error is exponentially convergent as a function of the truncation (indeed maybe in a wavenumber explicit way??), which is a key difference from the non-complex-scaled version. It must be true that the truncated equations (in both cases a) and b)) are a small perturbation of the original equations (so that coercivity + compactness is preserved)? Except for the truncation, the numerical analysis is probably straightforward, standard Galerkin analysis?

\subsection{Numerical results for the general case} \label{sec:NumGC}
Finally, to discretize problem  \eqref{eq:HSMM_casgeneral-variational}, we combine the previous tools with a classical Lagrange finite element approximation of the 2D unknown $u_b$.
Precisely, we use the HSM method to solve the problem of diffraction of the incident plane wave
$$
u^i(\bsx) = \exp(\ri k(x_1\cos(\pi/6)+x_2\sin(\pi/6)))
$$
by a perfectly reflecting scatterer  which is  the union of a disk and a triangle, this union contained in the square $\Omega_a$ with $a=0.8$.  The scattering problem is \eqref{pb:probleme_diff_initial}-\eqref{eq:Sommerfeld} in $\Omega$, the domain exterior to the scatterer, with $f=0$, $\rho =1$, and the Dirichlet boundary condition $u=g:= -u^i$ on $\partial \Omega$. The HSM problem we solve is a variation on \eqref{eq:HSMM_casgeneral-variational} in which: i) we replace $\Omega_b$ by $\widetilde \Omega_b := \Omega \cap \Omega_b$; ii) we replace $H^1(\Omega_b)$ by an affine subspace of $H^1(\widetilde \Omega_b)$ that respects the Dirichlet boundary condition; precisely, we seek $u_b \in \{v\in H^1(\widetilde \Omega_b): v=g \mbox{ on } \partial \Omega\}$. It is easy to see, by a straightforward variation of the arguments of Proposition \ref{lem-uniqueness-generalcase} and Theorem \ref{the-generalcase}, that this modification of \eqref{eq:HSMM_casgeneral-variational} remains well-posed: uniqueness holds and the modified sesquilinear form is coercive plus compact on $(L_0^2(\R))^4\times H^1_D(\widetilde \Omega_b)$, where $H^1_D(\widetilde \Omega_b):=\{v\in H^1(\widetilde \Omega_b): v=0 \mbox{ on } \partial \Omega\}$. {As in Section \ref{sec:NUM1}, this implies convergence of Galerkin methods of numerical solution, provided the sequence of approximation spaces used is asymptotically dense in $(L_0^2(\R))^4\times H^1_D(\widetilde \Omega_b)$.}

In our numerical implementation of a Galerkin method for this HSM formulation we fix $a=0.8$ and $b=1.2$, we use P2-elements for both the 2D unknown $u_b$ and the 1D unknowns $\varphi^j_{\theta}$, $j\in\llbracket0,3\rrbracket$, with a maximum element diameter for both meshes of $h=0.05$, choose the truncation parameter $T=5$, and choose $\theta=\pi/6$. The additional integrals, involving the operator $\Lambda_\theta$, which appear in the variational formulation are approximated, like the other ones, by a standard quadrature formula.

\begin{figure}[htbp]
		\begin{center}
			\begin{tikzpicture}
		\node[inner sep=0pt] (source) at (0,0)
		    {\includegraphics[width=13cm]{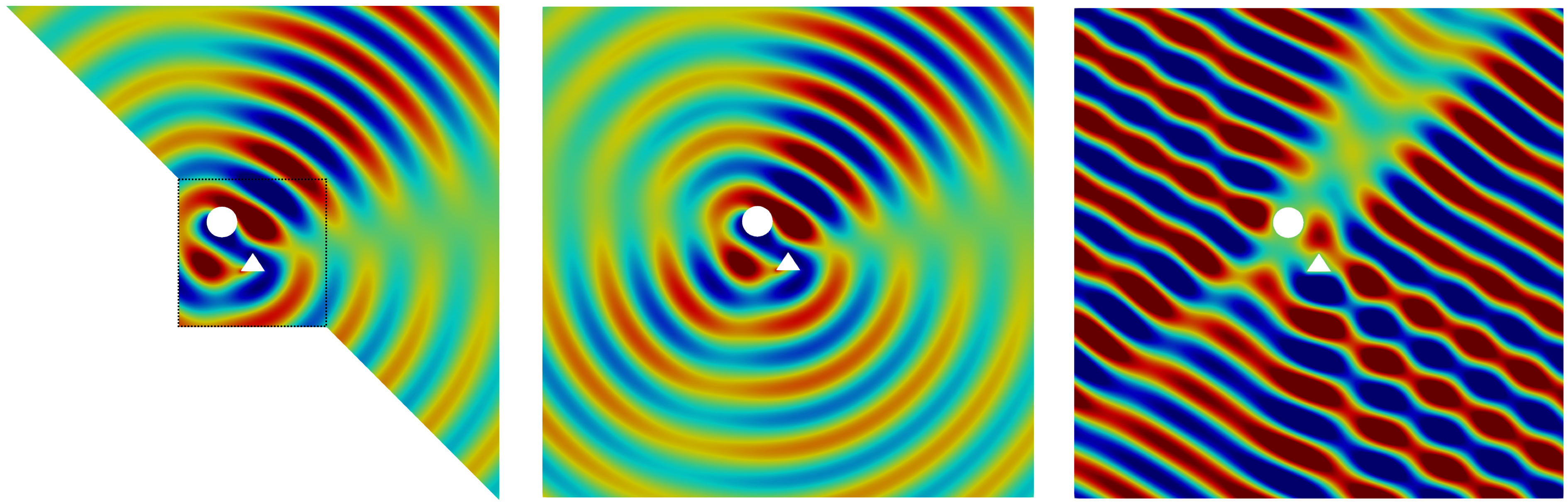}};
			%\draw(-4.5,-2.2) node {Incident field};
			\draw(-2,-2.2) node {Diffracted field};
			\draw(4.5,-2.2) node {Total field};
			\end{tikzpicture}
	\caption{The diffracted field $u$ reconstructed in $\Omega^b\cup\Omega^0_{\pi/4}\cup\Omega^1_{\pi/4}$ (left) and in the whole of $\Omega$ (middle), and the corresponding total field $u+u^i$ (right) for the case of scattering by a Dirichlet obstacle.}
	\label{fig:general_case}
	\end{center}
	\end{figure}

To reconstruct the solution everywhere several choices are possible, since $\theta<\pi/4$ and the different reconstruction domains overlap (for instance $\widetilde \Omega^b$ with $\Omega^0_\theta$, or $\Omega^0_\theta$ with $\Omega^1_\theta$, ...). We have reconstructed the solution as
\begin{eqnarray*}
	u({\bsx})&= &u_b({\bsx}),\quad {\bsx}\in\widetilde \Omega_b,\\
	u({\bsx}^j)&= &U_\theta(\varphi_\theta^j)({\bsx}^j),\quad {\bsx}^j\in \Omega^j_{\pi/4}\setminus\widetilde \Omega_b,\;\;j\in \llbracket 0, 3 \rrbracket,
\end{eqnarray*}
where $U_\theta$ is defined in \eqref{eq:expr_Ujtheta} and where we have used that $ \Omega^j_{\pi/4}\subset  \Omega^j_{\pi/6}$. See Figure \ref{fig:general_case} for the reconstruction in $\widetilde \Omega_b\cup\Omega^0_{\pi/4}\cup\Omega^1_{\pi/4}$ (left), the reconstruction in the whole domain $\Omega$ (middle), and the corresponding total field $u+u^i$ (right).
We notice that the different deformed half-space representations {(the representations in $\Omega_{\pi/4}^j\setminus\widetilde \Omega_b$, $j\in \llbracket 0, 3 \rrbracket$) are compatible between themselves, and are also compatible with the 2D solution $u_b$ in $\widetilde \Omega_b$}.

\section{Perspectives}\label{sec:Perspectives}
Our new complex-scaled HSM method has been presented in this paper for a relatively simple configuration. We expect that it can be extended easily to more complex problems for which the relevant half-space Green's functions are known sufficiently explicitly, such as acoustic scattering in stratified media, including cases where the stratification is different in different half-spaces (see e.g. \cite{Ott} for a presentation of the method in the dissipative case). The method is also expected to work well in at least some cases with infinite boundaries, for instance scattering by an infinite wedge with Dirichlet, Neumann or Robin boundary conditions. In all these cases, the complex-scaled HSM should be a convenient way to take into account possible surface/guided waves that propagate towards infinity. Elastic scattering in isotropic media can also be considered. More challenging extensions are to the cases where PMLs are observed to fail, such as anisotropic media. A potential advantage over PML of the complex-scaled HSM method in such cases is that it requires the existence of exponentially decaying analytical continuation of the traces of the solution only in a few directions (on the boundaries of a few half-planes).
\appendix
\section{Properties of the operators $\mathbb{D}$ and $\mathbb{D}_\theta$}\label{Appendix_Matrices}
The operator $\mathbb{D}$ defined by \eqref{eq:def_Dmat} can be rewritten as%\scednote{This is embarrassing, but I don't understand the $\otimes$ notation here!}
\[
		\mathbb{D}={D}\,S\otimes\mathbb{J} + S\,D\otimes \mathbb{J}^*\quad\text{with}\quad\mathbb{J}:=\left[\begin{matrix}
		 0&1&0&0\\
		 0&0&1&0\\
		 0&0&0&1\\
		 1&0&0&0
		 \end{matrix}\right],
\]
where the operator $S\in \mathcal{L}(L^2(\mathbb{R}))$ (the space of continuous linear operators on $L^2(\R)$)  is defined in \eqref{eq:sym} and the operator $D\in \mathcal{L}(L^2(\mathbb{R}))$ is defined in \eqref{eq:opD_comp_def} and \eqref{eq:DdefExt}.
Here $A \otimes \mathbb{M}$ denotes the tensor product of an operator $A \in \mathcal{L}(L^2(\mathbb{R}))$ with a $4 \times 4$ scalar matrix $\mathbb{M}$ (see, e.g., \cite[\S 12.4]{Aubin-2000}), which yields an operator of $\mathcal{L}((L^2(\mathbb{R}))^4)$. But it is actually enough to see this as a simple notation which makes the writing of the proof below easier: $A \otimes \mathbb{M}$ is the block operator matrix obtained by multiplying each scalar component of $\mathbb{M}$ by the operator $A$. One can easily verify that it satisfies the basic property $\|A \otimes \mathbb{M}\| \leq \|A\| \, \|\mathbb{M}\|.$

The operator $\mathbb{D}_\theta$, given by \eqref{eq:def_Dmat_cmplx}, has the same definition as $\mathbb{D}$ just replacing $D$ by $D_\theta$. The operators $D$ and $D_\theta$ satisfy similar properties, given respectively in Proposition \ref{prop:D_complexe} and in Propositions \ref{prop:D_realk_cmplx_plus}, \ref{prop:D_realk_cmplx_moins}, and \ref{prop:D_realk_cmplx_tot}.

We show in this appendix properties (i) and (ii) of Theorem \ref{th:HSMM_complexe} and Theorem \ref{th:HSMM_cmplxed}. These results are properties of $\mathbb{D}$ and $\mathbb{D}_\theta$, respectively, that are based on {the above} properties of $D$ and $D_\theta$. Let us give now the proof for (i) and (ii) of Theorem \ref{th:HSMM_complexe} for $\mathbb{D}$;  the same proof holds for Theorem \ref{th:HSMM_cmplxed} for $\mathbb{D}_\theta$.

As $D$ and $S$ are continuous operators on $L^2(\R)$, the continuity of $\mathbb{D}$ on $(L^2(\R))^4$ (part (i) of Theorem \ref{th:HSMM_complexe}) is obvious.

To show part (ii) of Theorem \ref{th:HSMM_complexe} we  consider $\mathbb{D}$ as an operator on $(L^2_0(\R))^4$. Let us denote by $\chi_+$ (respectively $\chi_-$) the characteristic function of $(a,+\infty)$ (respectively $(-\infty,-a)$). We have, by using \eqref{eq:decomp_L2}, that, for $\varphi\in L^2_0(\R)$
\[
\varphi=\chi_+\varphi+\chi_-\varphi,
\]
	and
	\begin{equation}\label{eq:norm_phi_chipm} \|\varphi\|_{L^2(\R)}^2=\|\chi_+\varphi\|_{L^2(\R)}^2+\|\chi_-\varphi\|_{L^2(\R)}^2.
\end{equation}
We can reformulate items (ii) and (iii) of Proposition \ref{prop:D_complexe} using these characteristic functions as follows:

\vspace{0.5ex}

\begin{itemize}
		\item[(ii)] $D\chi_+=L\chi_++K\chi_+\;\in\mathcal{L}(L^2(\R),L^2(a,+\infty))$, where $L, K\in \mathcal{L}(L^2(a,+\infty))$ are such that $\|L\|\leq 1/\sqrt{2}$ % is of norm less or equal to $1/\sqrt{2}$
and $K$ is compact.
		\item[(iii)] $D\chi_-\; \in\mathcal{L}(L^2(\R),L^2(a,+\infty))$ is compact.
		\end{itemize}

\vspace{0.5ex}

\noindent These properties, and that $S\,\chi_{{\pm}}=\chi_{{\mp}}\,S$, lead us to decompose $\mathbb{D}$ as an operator on $(L^2_0(\R))^4$ as
\[
	\mathbb{D}=\mathbb{L}+\mathbb{K}
\]
where
\[
	\mathbb{K}:=({D}\,S\chi_+ +{K}\,{S\, \chi_-})\otimes \mathbb{J} + (S\,{D}\chi_-+S\,K\chi_+)\otimes \mathbb{J}^*
\]
and
\[
	\mathbb{L}:=({L}\,{S\chi_-})\otimes \mathbb{J} + (S\,L\chi_+)\otimes \mathbb{J}^*.
\]
%\textcolor{red}{U}sing the identity ${D}S\chi_+={D}\chi_-S\chi_+$, we deduce
It follows from (ii) and (iii) that $\mathbb{K}$ is compact.
Moreover, noting that $\chi^+\,L=L$, so that $\chi^-\,S\,L=S\,L$, we have
\[
	\mathbb{L}=\chi_+\left({L}\,S\otimes \mathbb{J}\right) \,\chi_-+ \chi_-\left(S\,L\otimes \mathbb{J}^*\right)\,\chi_+.
\]
We deduce then by \eqref{eq:norm_phi_chipm} that, for all $\Phi\in (L^2_0(\R))^4$,
\begin{align*}
		 \|\mathbb{L}\Phi\|_{L^2(\R)^4}^2& =\|\left({L}\,S\otimes \mathbb{J}\right) \chi_-\Phi\|_{L^2(\R)^4}^2+\|\left(S\,L\otimes \mathbb{J}^*\right)\chi_+\Phi\|_{L^2(\R)^4}^2\\
		&\leq \frac{1}{2}\|\chi_-\Phi\|_{L^2(\R)^4}^2+\frac{1}{2}\|\chi_+\Phi\|_{L^2(\R)^4}^2=\frac{1}{2}\|\Phi\|_{L^2(\R)^4}^2,
\end{align*}
where we have used that $\|{L}\,S\otimes \mathbb{J}\|\leq \|{L}\,S\|\,\|\mathbb{J}\| \leq 1/\sqrt{2}$, and the same bound for $S\,L\otimes \mathbb{J}^*$.

{\section{Technical lemmas}\label{Appendix_R}
The lemmas in this annex  (cf.\ \cite[Lemma 4.4]{NonConvex}) concern the complex functions
\begin{equation}
\label{eq-defRandRhat}
R(\hat{z},z):=(\hat{z}^2+z^2)^{1/2},\quad \hat{R}(\hat{z},z):=(\hat{z}+z^2)^{1/2},\quad z,\,\hat{z}\in\C.
\end{equation}
Note that, as throughout the rest of the paper, all square roots in this appendix are principal square roots, i.e. square roots with argument in the range $(-\pi/2,\pi,2]$.
 \begin{lem} \label{lem:Rbound}
We have
\[
	|R(t,z)|^2\geq |\cos(\mathrm{Arg}(z))|\,\left(t^2+|z|^2\right),\quad t\in\R,\;z\in\C,\;z\neq 0.
\]
\end{lem}
\begin{proof} Let $\gamma:= \mathrm{Arg}(z)$. We have
$$
|R(t,z)|^4 = |t^2+|z|^2\re^{2\ri\gamma}|^2 = t^4+|z|^4 + 2t^2|z|^2\cos(2\gamma)=\left(t^2-|z|^2\right)^2+4t^2|z|^2\cos^2(\gamma)
$$
which yields
\[
	|R(t,z)|^4 \geq\cos^2(\gamma)\left[\left(t^2-|z|^2\right)^2+4t^2|z|^2\right]= \cos^2(\gamma)\,\left(t^2+|z|^2\right)^2.
\]
\end{proof}
\begin{lem}
	\label{lem:Im(R-z)}
	There exists a constant $C>0$ such that, for all $z\in\C$ and $t\in\R$,
		\[
		\dsp|{R}(t,z)-z|\leq C\frac{t^2}{|z|},\quad  \Re(z)>0,\,|z|\geq |t|;
		\]
	in particular
		\[
		\Im({R}(t,z))\geq \Im(z)-C\frac{t^2}{|z|},\quad  \Re(z)>0,\,|z|\geq|t|.
		\]
	Moreover, for all $A>0$, there exists a constant $C'{>0}$ such that, for  $-A\leq t\leq A$,
	\[
	\Im({R}(t,z))\geq \Im(z)-C',\quad  \Re(z)>0.
	\]
	\end{lem}
	\begin{proof}
Since the function $z\mapsto (1+z)^{1/2}-1$ is an analytic function of $z$ in the open unit disk that vanishes at the origin and is bounded in the closed disk, for some constant $C>0$,
\begin{equation}
\label{eq:sqrt1+z2}
|(1+z)^{1/2}-1|\leq C|z|,\quad |z|\leq 1.
\end{equation}
Further, the function
$$z\mapsto(t^2+z^2)^{1/2}- z\left(1+\frac{t^2}{z^2}\right)^{1/2}$$
is analytic in the domain $\{z\in\C:  \Re(z)>0\} $ and it vanishes for $z\in \R$ with $z>0$, so that it vanishes everywhere in $\{z\in\C:  \Re(z)>0\} $. This implies that
$${R}(t,z)-z=z\left[\left(1+\frac{t^2}{z^2}\right)^{1/2}-1\right],\quad \Re(z)>0.$$
Combined with (\ref{eq:sqrt1+z2}), this proves  the first inequality.

%This provides the first inequalities since
%$$\hat{R}(\hat{z},z)-\sqrt{z^2}=\sqrt{z^2}\left[\left(1+\frac{\hat{z}}{z^2}\right)^{1/2}-1\right]$$
%and
%$\sqrt{z^2}=z$ if $\Re(z)>0$ and $\sqrt{z^2}=-z$ if $\Re(z)<0$.
 The second inequality is a direct consequence of the first since
 $$ |{R}(t,z)-z|\geq \Im(z-{R}(t,z)).$$
 Finally, the second inequality yields that, for any constant $A>0$ and for all $|t|\leq A\leq|z|$,
$\Im({R}(t,z))\geq \Im(z)-CA$. But a similar inequality also holds for   $|z|\leq A$, possibly with another constant, since $|R(t,z)|\leq \sqrt{2}A$ for $|t|\leq A$ and $|z|\leq A$.
		\end{proof}
\begin{rem}\label{rem:Im(R+z)}
Note that, replacing $z$ by $-z$, one directly deduces from the previous lemma similar results valid for $\Re(z)<0$. In particular,
		\[
	\dsp|{R}(t,z)+z|\leq C\frac{t^2}{|z|},\quad  \Re(z)<0,\,|z|\geq|t|,
	\]
and
	\[
	\Im({R}(t,z))\geq -\Im(z)-C\frac{t^2}{|z|},\quad  \Re(z)<0,\,|z|\geq|t|.
	\]
\end{rem}

Analogous results can be proved in the more general case where the positive real number $t^2$ is replaced by any complex number $\hat{z}$.
 \begin{lem} \label{lem:Im(Rhat-z)}
There exists a constant $C>0$ such that, for all $z,\hat{z}\in\C$, ${z}\neq 0$,
\begin{eqnarray}\nonumber
|\hat{R}(\hat{z},z)-z|\leq C\left|\frac{\hat{z}}{z}\right|,\quad  \Re(z)\geq|\hat{z}|^{1/2},\\ \label{eq:extra}
\Im(\hat{R}(\hat{z},z))\geq \Im(z)-C\left|\frac{\hat{z}}{z}\right|,\quad  \Re(z)\geq|\hat{z}|^{1/2}.
\end{eqnarray}
Moroever, if $U$ is a bounded subset of $\C$ and %$\gamma\in\R$ is such that
$0<\gamma<\pi/2$, there exists a constant $C'>0$ such that
\[
\Im(\hat{R}(\hat{z},z))\geq \Im(z)-C',\quad  \hat{z}\in U,\; |\Arg(z)|<\gamma.
\]
\end{lem}
\begin{proof}
To proceed  as in the proof of Lemma \ref{lem:Im(R-z)}, we just have to show that the function
$$z\mapsto(\hat{z}+z^2)^{1/2}- z\left(1+\frac{\hat{z}}{z^2}\right)^{1/2}$$
is analytic in the domain $\{\Re(z)>|\hat{z}|^{1/2}\} $. For that, we have to check that the branch cuts of the functions  $z\mapsto (\hat{z}+z^2)^{1/2}$ and  $z\mapsto \left(1+ {\hat{z}}/{z^2}\right)^{1/2}$ do not intersect the domain $\Re (z)> |\hat{z}|^{1/2}$.

%It is straightforward concerning
 This is clear for the function $z\mapsto \left(1+ {\hat{z}}/{z^2}\right)^{1/2}$. %On the other hand,
The branch cut of the function $z \mapsto (\hat{z}+z^2)^{1/2}$ is the  subset of the hyperbola $\{z\in\C:\Im (z^2+\hat{z})=0\}$ where $\Re (z^2+\hat{z})\leq 0$.
%Only one branch of this hyperbola intersects the domain $\Re (z)> |\hat{z}|^{1/2}$.
The intersection of this hyperbola with the domain $\Re (z)> |\hat{z}|^{1/2}$ is the  connected  set
$$\hat{H}=\{z\in\C: \Im (z^2+\hat{z})=0,\,\Re (z)> |\hat{z}|^{1/2}\},$$
which describes a curve which is asymptotic to the real axis when $\Re (z)\rightarrow +\infty$. To conclude, we have to prove that for all $z\in \hat{H}$, $\Re (z^2+\hat{z})>0$.  This is clearly true for large values of $\Re (z)$. If it were not true for all $z\in \hat{H}$, there would exist some $z_0\in \hat{H}$ such that  $\Re (z_0^2+\hat{z})=0$. But then $z_0^2+\hat{z}=0$, so that  $|z_0|=|\hat{z}|^{1/2}$, which is impossible since $\Re (z_0)> |\hat{z}|^{1/2}$.
\end{proof}

An important application of the previous lemmas   for our purpose is the following: %\scednote{The statement of the first part reworded as discussed at last meeting.}%\scednote{OK, I'm struggling to follow this a bit, need a bit more time! Is this meant to do both the cases $\Re(z)>a$ and $\Re(z)<-a$ at the same time? {\color{red} Note added 30/9/20. It is slightly odd that $m(\theta,z)$ and $(|z-a|^2+|\tau_\theta(s)-a|^2)$ depend on $z$, but I think you want to keep this and it is fine. But I think we need to clarify what $C$ depends on, in particular (I think) that $C$ does not depend on $z$?}}\asbbednote{I think that $C$ depends only on $\theta$. I have changed the statement accordingly. }
\begin{lem}
	\label{lem-Rbound-z-tau}
For $\theta\in(0,\pi/2)$, let $\tau_\theta$ be defined by \eqref{eq:parametrage_t_theta}, and suppose that %\scednote{Late change, as I'm working through proof of Lemma 4.1: changed $($ to $[$ in range of $\theta_0$.}
$\theta_0 \in [\theta,\pi/2)$. Then, for $w\in\C$ such that $w\neq 0$ and $-\theta_0+\theta\leq \Arg(w)\leq \theta_0$, %there exists a positive constant $m(\theta,z)$ such that,
and for all $s\in\R$, it holds that
\begin{equation}\label{eq:Rbound-z-tau}
	|R(w,\tau_\theta(s)-a)|^2\geq \cos(\theta_0)(|w|^2+|\tau_\theta(s)-a|^2) \geq \cos^2(\theta_0)(|w|^2+|\tau_\theta(s)-a|^2).
	\end{equation}
%where $m(\theta,z)=\min (\cos(\Arg(w)),\cos(\theta-\Arg(w)))$.
Further, there exists a constant $C>0$, depending only on $a$ and $\theta$, such that, for all $s\in \R$, 	
%for $z\in\C$ such that $z\neq a$ and $0\leq \Arg(w)\leq\theta$, there exists a constant $C>0$ such that for all $s\in\R$,
\begin{equation}\label{eq:ImRbound-z-tau}
	\Im(R(w,\tau_\theta(s)-a))\geq [\cos(\theta-\Arg(w))]^{1/2} \Im(w)-C \geq \cos(\theta-\Arg(w)) \Im(w)-C,
\end{equation}
provided $w\neq 0$ and $0\leq \Arg(w)\leq\theta$.
\end{lem}
\begin{proof}
 Let $\gamma:= \Arg(w) \in 	[\theta-\theta_0,\theta_0]$. We have, for all $s\in\R$,
 $$ |R(w,\tau_\theta(s)-a)|=|R(|w|\re^{\ri\gamma},\tau_\theta(s)-a)|=|R(|w|,(\tau_\theta(s)-a)\re^{-\ri\gamma})|.$$
Applying Lemma \ref{lem:Rbound} we obtain the inequality
 \begin{equation}\label{eq:module_R_lemmeB4}|R(w,\tau_\theta(s)-a)|^2\geq |\cos( \hat{\gamma}(s))|(|w|^2+|\tau_\theta(s)-a|^2)
	 \end{equation}
 where $ \hat{\gamma}(s):=\Arg\left((\tau_\theta(s)-a)\re^{-\ri\gamma}\right)$.

We now consider separately the three cases $|s|\leq a$, $s>a$, and $s< -a$ to derive the two inequalities of the lemma.

Case 1. For $|s|\leq a$, since $\tau_\theta(s)=s$, we have $ |\cos( \hat{\gamma}(s))|=\cos( \gamma) \geq \cos(\theta_0)$ so that \eqref{eq:Rbound-z-tau} holds.
%	\begin{equation}\label{eq:module_R_lemmeB4_|s|leqa}
%		|R(w,\tau_\theta(s)-a)|^2\geq \cos({\gamma})(|w|^2+|\tau_\theta(s)-a|^2).
%		\end{equation}
	Further, applying the third inequality of Lemma \ref{lem:Im(R-z)} we get that there exists a constant $c>0$, dependent only on $a$, %independent of $z$,
such that
\begin{equation}\label{eq:Im_R_lemmeB4_|s|leqa}
	\Im(R(w,\tau_\theta(s)-a)) =\Im(R(s-a,w))\geq  \Im(w)-c.
\end{equation}

Case 2. For $s>a$, since $\tau_\theta(s)-a= \tilde s\re^{\ri\theta}$, where $\tilde s := s-a> 0$, and $\hat{\gamma}(s)=\theta-\gamma$ so that $ |\cos( \hat{\gamma}(s))|=\cos(\theta-\gamma) \geq \cos(\theta_0)$, we have from \eqref{eq:module_R_lemmeB4} that
%\scednote{Have rephrased the bits in red in the next equation so that the reference back to this equation from the 3rd case is clearer (at least to me).}
	\begin{equation}\label{eq:module_R_lemmeB4_sgeqa}
		|R(w,\tau_\theta(s)-a)|^2 = |R(w,\tilde s \re^{\ri \theta})|^2\geq \cos(\theta-\gamma)(|w|^2+\tilde s^2), %\cos({\theta-\gamma})(|w|^2+|\tau_\theta(s)-a|^2).
		\end{equation}
so that \eqref{eq:Rbound-z-tau} holds. To see that \eqref{eq:ImRbound-z-tau} holds, for $\gamma \in [0,\theta]$ as required in the lemma, note that
\begin{equation}\label{eq:argument_R_lemmeB4_sgeqa}
\Arg(R(w,\tau_\theta(s)-a))=\Arg(R(|w|\re^{\ri\gamma},\tilde{s}\re^{\ri \theta})) \in [\gamma,\theta].
\end{equation}
% so that 	
%	$$\Arg(R(|w|e^{\ri\gamma},\tilde{s}e^{\ri \theta}))=\gamma+\Arg(R(|w|,\tilde{s}e^{\ri(\theta-\gamma)}))$$ from which we deduce the bounds
%	\begin{equation}
%	\label{eq:argument_R_lemmeB4_sgeqa}
%\gamma\leq	\Arg(R(|w|e^{\ri\gamma},\tilde{s}e^{\ri \theta}))\leq \theta.
%	\end{equation}
Combining this with  the bound \eqref{eq:module_R_lemmeB4_sgeqa} we see that%on $|R(w,\tau_\theta(s)-a)|$, we have
%$$  \Im(R(w,\tau_\theta(s)-a)))\geq \sqrt{\cos({\theta-\gamma})}|w|\sin(\gamma)$$ which gives
\begin{eqnarray}\nonumber
	 \Im(R(w,\tau_\theta(s)-a)) & = & \Im(R(w,\tilde s \re^{\ri \theta})) = |R(w,\tilde s \re^{\ri \theta})| \, \sin\left(\Arg(R(w,\tilde s \re^{\ri \theta}))\right)\\ \label{eq:Im_R_lemmeB4_sgeqa}
&  \geq &\sqrt{\cos(\theta-\gamma)}\,|w|\sin(\gamma)=\sqrt{\cos(\theta-\gamma)}\, \Im(w).
\end{eqnarray}

Case 3. For $s<-a$,  setting $\tilde{s}:=-(s+a)>0$, we have  $\tau_\theta(s)-a=-(\tilde{s}\re^{\ri\theta}+2a)$ and $ |\cos( \hat{\gamma}(s))|=|\cos(\tilde{\theta}(s)-\gamma)|$,   where $\tilde{\theta}(\tilde{s}):=\Arg(\tilde{s}\re^{\ri\theta}+2a)$.  Since $0\leq \tilde{\theta}(\tilde{s})\leq \theta$ for $\tilde{s}>0$, and $\gamma\in [\theta-\theta_0,\theta_0]$,  we have
	\[
		|\cos( \hat{\gamma}(s))|\geq \min_{0\leq \tilde{\gamma}\leq \theta}|\cos(\gamma-\tilde{\gamma})|\geq \cos(\theta_0),% \min (\cos(\gamma),\cos(\theta-\gamma)),
	\]
	 so that \eqref{eq:Rbound-z-tau} follows from \eqref{eq:module_R_lemmeB4}.
	%\begin{equation}\label{eq:module_R_lemmeB4_sleq-a}
%	|R(w,\tau_\theta(s)-a)|^2\geq \min (\cos(\gamma),\cos(\theta-\gamma))(|w|^2+|\tau_\theta(s)-a|^2).
%	\end{equation}
	
To show \eqref{eq:ImRbound-z-tau}, for $0\leq \gamma\leq \theta$ as required, we rewrite%\scednote{Nice! I follow this step-by-step but have no idea how you thought of this!} $R(w,\tau_\theta(s)-a)$ as
\begin{eqnarray}\nonumber
		R(w,\tilde{s}\re^{\ri\theta}+2a)&=&[R(w,\tilde{s}\re^{\ri\theta})^2+ 4a^2+4a\tilde{s}\re^{\ri\theta}]^{1/2}\\ \label{eq:ASCleverRewriting}
		&=&\hat{R}(4a^2+4a\tilde{s}\re^{\ri\theta},R(w,\tilde{s}\re^{\ri\theta})),
\end{eqnarray}
where $\hat{R}$ is defined by \eqref{eq-defRandRhat}. Using \eqref{eq:module_R_lemmeB4_sgeqa} and \eqref{eq:argument_R_lemmeB4_sgeqa}
we see that%\scednote{I have removed the $\Re$ from this next equation as I couldn't see what role it played - it just made things more complicated? - AS: I completely agree}
			$$|R(w,\tilde{s}\re^{\ri\theta})|^2  \geq \cos(\theta-\gamma) \tilde s^2 \geq \cos(\theta) \tilde s^2 \quad {\mbox{ and }\quad \Arg(R(w,\tilde{s}\re^{\ri\theta})) \in [\gamma,\theta].}$$
 Thus, by \eqref{eq:extra} %Lemma \ref{lem:Im(Rhat-z)}
 applied to \eqref{eq:ASCleverRewriting}, there exist {constants $C>0$ and $\tilde s_0>0$}, depending only on $\theta$  and $a$,  such that%\scednote{Have removed $\tilde s_0$ - the following bound holds for all $\tilde s$. AS: it seems to be ok...}
	\[
		\Im(R(w,\tilde{s}\re^{\ri\theta}+2a))\geq \Im(R(w,\tilde{s}\re^{\ri\theta}))-{C}\frac{1+\tilde{s}}{\tilde{s}},
	\]
{for $\tilde s\geq \tilde s_0$.}
Hence, and using \eqref{eq:Im_R_lemmeB4_sgeqa}, it follows that, for $\tilde{s}\geq {\max(1,\tilde s_0)}$,
	\[
		\Im(R(w,\tilde{s}\re^{\ri\theta}+2a))\geq \sqrt{\cos({\theta-\gamma})}\Im(w)-2C.
	\]
On the other hand, for $\tilde{s}\leq {\max(1,\tilde s_0)}$  the last inequality of Lemma \ref{lem:Im(Rhat-z)}  gives that there exists a constant $C'>0$, depending only on $a$, such that
	\[
		\Im(R(w,\tilde{s}\re^{\ri\theta}+2a))\geq \Im(w)-C'.
	\]
Gathering the two estimates we have that, for $s<-a$,
	\begin{equation}\label{eq:Im_R_lemmeB4_sleq-a}
		 \Im(R(w,\tau_\theta(s)-a)))\geq \sqrt{\cos({\theta-\gamma})}\Im(w)-\max(C',2C).
	\end{equation}

\vspace{0.5ex}

We have shown that \eqref{eq:Rbound-z-tau} holds in each case. That \eqref{eq:ImRbound-z-tau} also holds follows from \eqref{eq:Im_R_lemmeB4_|s|leqa}, \eqref{eq:Im_R_lemmeB4_sgeqa}, and \eqref{eq:Im_R_lemmeB4_sleq-a}, trivially noting that $(\cos(t))^{1/2}\geq \cos(t)$ for $t\in \R$.
\end{proof}
%
%
% We will distinguish the three cases $|s|<a$, $s\geq a$ and $s\leq -a$.
%
% For $|s|<a$,   since $\tau_\theta(s)=s$, we deduce from Lemma \ref{lem:Rbound} that:
% $$|R(z-a,s-a)|^2\geq \cos(\gamma)(|z-a|^2+(s-a)^2),$$
% which gives the result since $\cos(\gamma)\leq m(\theta). $
%
% For $s>a$, we have
% $$ |R(z-a,\tau_\theta(s)-a)|=|R(|z-a|e^{i\gamma},(s-a)e^{i\theta})|=|R(|z-a|,(s-a)e^{i(\theta-\gamma)})|.$$
% Then we can again apply Lemma \ref{lem:Rbound} and we obtain the following inequality:
% $$|R(z-a,\tau_\theta(s)-a)|^2\geq \cos(\theta-\gamma)(|z-a|^2+|s-a|^2).$$
%
% Finally, for $s<-a$, setting $\tilde{s}=-(s+a)$, we have  $\tau_\theta(s)-a=-(\tilde{s}e^{i\theta}+2a)$. Therefore:
%$$ |R(z-a,\tau_\theta(s)-a)|=|R(|z-a|,(\tilde{s}e^{i\theta}+2a)e^{-i\gamma})|.$$
%Using again  Lemma \ref{lem:Rbound}, this implies that
% $$|R(z-a,\tau_\theta(s)-a)|^2\geq \cos(\tilde{\theta}-\gamma)(|z-a|^2+|\tilde{s}e^{i\theta}+2a|^2)$$
% where $\tilde{\theta}=\Arg(\tilde{s}e^{i\theta}+2a)$. Since $0\leq \tilde{\theta}\leq \theta$ and
% $$|\tilde{s}e^{i\theta}+2a|^2\geq \tilde{s}^2,$$
}
\section{Mapping properties of complex-scaled integral operators} \label{sec:ComplIntOper}
In this appendix we prove mapping properties of the analytic continuations into the complex plane of the single- and double-layer potential operators $\cS^j$ and $\cD^j$, defined by \eqref{eq:SDjdef} for $j\in\llbracket0,3\rrbracket$, $\phi\in L^2(\Sigma_a)$, and $|t|>a$. The proofs use the  bounds established in Appendix \ref{Appendix_R} and bounds on the relevant Hankel functions. Specifically \cite[Lemma 3.4]{DHSLMM}, for some constant $c_1>0$,
\begin{equation} \label{eq:Hbound1}
\left|\re^{-\ri z} H_1^{(1)}(z)\right| \leq c_1\left( |z|^{-1}+(1+|z|)^{-1/2}\right), \quad \Re(z)>0.
\end{equation}
Similarly, it follows from \eqref{eq:HankAsym} and \cite[\S10.2(ii)]{DLMF} that, for some constant $c_0>0$,
\begin{equation} \label{eq:Hbound0}
\left|\re^{-\ri z} H_0^{(1)}(z)\right| \leq c_0\, M\!(|z|),  \quad \Re(z)>0,
\end{equation}
where  $M\!(t) := \log(1+t^{-1}) + (1+t)^{-1/2}$, for $t>0$.

\begin{prop}  \label{lem:SDbounds}%\scednote{For the moment I have not inserted $k$ into $M(|z-a|)$. I'm reluctant to do this as: i) it makes the formula a little more complicated; ii) we've indicated already that $C$ may depend on $k$; iii) our formulas are not dimensionally correct, and that is normal in mathematics papers: e.g.\ should $|z-a|^{-1/2}$ be changed to $(k|z-a|)^{1/2}$?} \label{lem:SDbounds}
For every $\theta\in (0,\pi/2)$ there exists some constant $C>0$, that depends only on $a$, $k$, and $\theta$, such that, for $j\in\llbracket0,3\rrbracket$ and every $\phi\in L^2(\Sigma_a)$,
\begin{eqnarray} \label{eq:Sbound1}
|\cS^j\phi(z)| &\leq & C \,M\!\!\left(|z-a|\right)\,\exp(-k\Im(z))\,\|\phi\|_{L^2(\Sigma_a)},\\ \label{eq:Dbound1}
|\cD^j\phi(z)| &\leq & C \,|z-a|^{-1/2}\,\exp(-k\Im(z))\,\|\phi\|_{L^2(\Sigma_a)},
\end{eqnarray}
for $\Re(z) > a$  with $|\Arg(z-a)|\leq \theta$, while
\begin{eqnarray} \label{eq:Sbound2}
|\cS^j\phi(z)| &\leq & C \,M\!\!\left(|z+a|\right)\,\exp(k\Im(z))\,\|\phi\|_{L^2(\Sigma_a)},\\ \label{eq:Dbound2}
|\cD^j\phi(z)| &\leq & C \, |z+a|^{-1/2}\,\exp(k\Im(z))\,\|\phi\|_{L^2(\Sigma_a)},
\end{eqnarray}
for $\Re(z) < -a$  with $|\Arg(-z-a)|\leq \theta$.
\end{prop}
\begin{proof} We prove only the bounds (\ref{eq:Sbound1}-\ref{eq:Dbound1}); the proofs of (\ref{eq:Sbound2}-\ref{eq:Dbound2}) are identical. Throughout this proof $C$ will denote any positive constant depending  only on $a$, $k$, and $\theta$, not necessarily the same at each occurrence, and we assume that $\bsy^j\in \Sigma_a$ {and} that $\Re(z)>a$ with $|\Arg(z-a)|\leq \theta$. % and, without loss of generality, that $\theta\in [\pi/4,\pi/2)$.

It follows from (the analytic continuation of) \eqref{eq:SDjdef} and the Cauchy-Schwarz inequality that
$$
|\cS^j\phi(z)| \leq C I^{1/2}_\cS \|\phi\|_{L^2(\Sigma_a)}, \quad |\cD^j\phi(z)| \leq C I^{1/2}_\cD \|\phi\|_{L^2(\Sigma_a)},
$$
where
\begin{eqnarray*}
I_\cS := \int_{\Sigma_a} |\Phi(\bsx^j(z),\bsy^j)|^2 \,\rd s(\bsy^j), \quad I_\cD := \int_{\Sigma_a} \left|\frac{\partial \Phi(\bsx^j(z),\bsy^j)}{\partial n(\bsy^j)}\right|^2\, \rd s(\bsy^j).
\end{eqnarray*}
Now, from \eqref{eq:Hbound0}, and since $\bsx^j(z) := (a,z)$ and recalling the definition \eqref{eq:R(s,t)},
$$
|\Phi(\bsx^j(z),\bsy^j)| \leq C \,M\!\left(k|R(a-y_1^j,z-y_2^j)|\right)\left|\exp(\ri k R(a-y_1^j,z-y_2^j))\right|.
$$
Further, it follows from Lemma \ref{lem:Im(R-z)} that, for some constant $C_a>0$ depending only on $a$,% and $C=\re^{kC_a}$:
%\ref{lem:Im(R-z)} that  for any $C_a>0$, there exists $r_0$ such that
%$$\Im(R(a-y_1^j,z-y_2^j)\geq \Im(z)-C_a $$
%for $|z|\geq r_0$ and   all $\bsy^j\in \partial \Omega_a$. But this also holds, for some $C_a>0$, for $|z|\leq r_0$, by a simple continuity argument.  This leads to the following estimate
\begin{eqnarray} \label{eq:expbound}
|\exp(\ri k R(a-y_1^j,z-y_2^j))| \leq \exp(-k(\Im(z)-C_a))\leq C\exp(-k\Im(z)),
\end{eqnarray}
and from Lemma \ref{lem:Rbound} that
\begin{equation} \label{eq:RLowerBound}
|R(a-y_1^j,z-y_2^j)| \geq \cos(\theta)((a-y_1^j)^2+|z-a|^2)^{1/2}\geq \cos(\theta)|z-a|.
\end{equation}
Thus, noting that $M(t)$ is decreasing as $t$ increases and that, for every $c>0$, $t\mapsto M(ct)/M(t)$ is a bounded function on $t>0$, it follows that
$$
|\Phi(\bsx^j(z),\bsy^j)| \leq CM\!\left(|z-a|\right)\exp(-k\Im(z)),
$$
so that $I^{1/2}_\cS \leq CM\!\left(|z-a|\right)\exp(-k\Im(z))$ and \eqref{eq:Sbound1} follows.

Writing $n(\bsy^j) = (n_1(\bsy^j),n_2(\bsy^j))$ in the $(x_1^j,x_2^j)$ coordinate system, we see that
\begin{eqnarray} \nonumber
\left|\frac{\partial \Phi(\bsx^j(z),\bsy^j)}{\partial n(\bsy^j)}\right| &=&\\ \label{eq:normderiv}
& & \hspace{-6ex}  \frac{k}{4}\left|H_1^{(1)}(k R(a-y_1^j,z-y_2^j))\right|\frac{\left|n_1(\bsy^j)(y^j_1-a)+n_2(\bsy^j)(y_2^j-z)\right|}{\left|R(a-y_1^j,z-y_2^j)\right|}.
\end{eqnarray}
How we bound the right hand side of this equation depends on which side of $\Sigma_a$ the point $\bsy^j$ is located.
When $\bsy^j\in \Sigma_a^j\subset \Sigma_a$ the right hand side of \eqref{eq:normderiv} vanishes,  since then $y^j_1=a$ and $n_2(\bsy^j)=0$. When $\bsy^j\in \partial \Omega \setminus (\Sigma_a^j\cup\Sigma_a^{j+1})$ it holds that $y_1^j=-a$ or $y_2^j=-a$, so that $|R(\bsx^j(z),\bsy^j)|\geq C$ by  Lemma \ref{lem:Rbound}, from which it follows, applying \eqref{eq:Hbound1} and  \eqref{eq:expbound}, that the right hand side of \eqref{eq:normderiv} is bounded by $C\exp(-k\Im(z))$. Thus
$$
I_\cD \leq C\exp(-2k\Im(z)) + I_\cD^*, \mbox{ where } I_\cD^* := \int_{\Sigma_a^{j+1}} \left|\frac{\partial \Phi(\bsx^j(z),\bsy^j)}{\partial n(\bsy^j)}\right|^2\, \rd s(\bsy^j).
$$
Finally,  when $\bsy^j\in \Sigma_a^{j+1}$ it holds that $n_1(\bsy^j)=0$ and $y_2^j=a$, so that,  where $\rho:=|z-a|$, \eqref{eq:normderiv} says that%\scednote{I think, with these additions, this is now enough detail so that this argument is clear, and there isn't a need for a diagram, and indeed a diagram might be tricky to do: is this a diagram showing $\partial \Omega_a$ in $\R^2$, or showing the complex plane?}
\begin{eqnarray} \nonumber
\left|\frac{\partial \Phi(\bsx^j(z),\bsy^j)}{\partial n(\bsy^j)}\right| &=& \frac{k}{4}\left|H_1^{(1)}(k R(a-y_1^j,z-y_2^j))\right|\frac{\rho}{\left|R(a-y_1^j,z-y_2^j)\right|}\\ \label{eq:NormDivBound}
&\leq & C\left(\frac{\rho}{\rho^2+(a-y_1^j)^2}+(1+\rho)^{-1/2}\right)\exp(-k\Im(z)),
\end{eqnarray}
using \eqref{eq:Hbound1}, \eqref{eq:expbound}, and \eqref{eq:RLowerBound}. Thus
$$
\re^{2k\Im(z)}I_\cD^* \leq  \frac{C}{1+\rho} + C \int_{-a}^a  \frac{\rho^2\, \rd y^j_1}{(\rho^2+(y^j_1-a)^2)^2} %\leq  \frac{C}{1+\rho} + C\int_{-\infty}^0  \frac{\rho^2\rd s}
\leq \frac{C}{\rho},
$$
where we see the final bound by substituting $s:=(y^j_1-a)/\rho$. Thus $I^{1/2}_\cD\leq C|z-a|^{-1/2}\re^{-k\Im(z)}$ and the bound \eqref{eq:Dbound1} follows.
\end{proof}

The final result of this appendix, and arguments we make elsewhere in the paper, depend on a mapping property of the classical double-layer potential operator on the boundary of a quadrant. For the convenience of the reader we state this key and well-known result and sketch its proof in the following proposition.

\begin{prop} \label{prop:DLPquad} Let $h_0$ denote the kernel function $h$ defined in \eqref{eq:hpr_kernel} in the case that $k=0$, so that $h_0(x_1,x_2)=x_1/(\pi(x_1^2+x_2^2))$, for $x_1>0$ and $x_2\in \R$.
For  $\phi\in L^2(0,+\infty)$ define $\mathfrak{D}\phi:(0,+\infty)\to\C$ by
$$
\mathfrak{D}\phi(s) := \int_0^{+\infty} h_0(s,t)\, \phi(t)\, \rd t=\frac{1}{\pi}\int_0^{+\infty}\frac{s}{s^2+t^2}\, \phi(t)\, \rd t, \quad s>0.
$$
Then $\mathfrak{D}\phi\in L^2(0,+\infty)$ and the mapping $\mathfrak{D}: L^2(0,+\infty)\to L^2(0,+\infty)$ is bounded, with norm $\|\mathfrak{D}\|=1/\sqrt{2}$.
\end{prop}
\begin{proof} This result can be proved by making use of the equivalence of \eqref{eq:opD_comp_exprF} and \eqref{eq:opD_comp_exprG} in the static case $k=0$ as in \cite{BbFT18}, or directly via Mellin transform methods (cf.~\cite{Mitrea02,PerfektPutinar17}). Equivalently, we observe as in \cite{Ch:84} that the mapping $\mathfrak{I}:L^2(0,+\infty)\to L^2(\R)$, given by $\mathfrak{I}\phi(t) = \phi(\re^{-t})\re^{-t/2}$, $t\in\R$, is unitary, as is the Fourier transform operator $\mathfrak{F}:L^2(\R)\to L^2(\R)$, $\phi\mapsto \widehat \phi$,  given by \eqref{eq:FTdef}. Further \cite{Ch:84}, for $\phi\in L^2(\R)$,
$$
\mathfrak{I}\mathfrak{D}\mathfrak{I}^{-1}\phi(s) = \int_\R \kappa(s-t)\phi(t)\, \rd t, \quad s\in \R,
$$
where
$\kappa(\tau) := \re^{\tau/2}/(\pi(\re^{2\tau}+1))$, for $\tau\in \R$, and \cite{Ch:84}
$$
\widehat\kappa(\xi) = \frac{1}{\sqrt{2\pi}}\, \frac{\sinh(\pi(\xi-\ri/2)/2)}{\sinh(\pi(\xi-\ri/2))}, \quad \xi\in\R.
$$
Thus, for $\phi\in L^2(\R)$,
$$
\mathfrak{F}\mathfrak{I}\mathfrak{D}\mathfrak{I}^{-1}\mathfrak{F}^{-1}\widehat \phi = \sqrt{2\pi}\,\widehat\kappa\, \widehat \phi,
$$
so that \cite{Ch:84}, since $\widehat\kappa\in L^\infty(\R)$, $\mathfrak{D}$ is bounded on $L^2(0,+\infty)$ with
\[
\|\mathfrak{D}\|=\|\mathfrak{F}\mathfrak{I}\mathfrak{D}\mathfrak{I}^{-1}\mathfrak{F}^{-1}\| = \sqrt{2\pi}\,\|\widehat\kappa\|_{L^\infty(\R)} = \sqrt{2\pi}\, |\widehat \kappa(0)|=1/\sqrt{2}.
\]
\end{proof}

\begin{rem} \label{rem:D0bounded} Let $D_0$ denote the operator $D$, given by \eqref{eq:opD_comp_exprG} and \eqref{eq:DdefExt}, in the case that $k=0$, so that, for $\phi\in L^2(\R)$,  $D_0\phi(t)= 0$, for $t\leq a$, while
$$
D_0\phi(t) = \int_\R h_0(t-a,s-a)\, \phi(s)\, \rd s=\frac{1}{\pi}\int_\R \frac{t-a}{(t-a)^2+(s-a)^2}\, \phi(t)\, \rd t, \quad t>a,
$$
where $h_0$ is as defined in Proposition \ref{prop:DLPquad}.
 Then it is clear from the above proposition that, as an operator on $L^2(a,+\infty)$, $D_0$ is bounded with norm $\|D_0\|=\|\mathfrak{D}\|=1/\sqrt{2}$, and that $D_0$ is also bounded as an operator on $L^2(\R)$, with norm $\|D_0\|=2\|\mathfrak{D}\|=\sqrt{2}$.
\end{rem}

\begin{prop} \label{prop:SDbounded}
For $0<\theta<\pi/2$ and $j\in\llbracket0,3\rrbracket$, $\cS^j_\theta$ and $\cD^j_\theta$ are continuous operators from $L^2(\Sigma_a)$ to $L^2(-\infty,-a)\oplus L^2(a,+\infty)$, where, for $\phi\in L^2(\Sigma_a)$,
$$
\cS^j_\theta\phi(s) := \cS^j\phi(\tau_\theta(s)), \quad \cD^j_\theta\phi(s) := \cD^j\phi(\tau_\theta(s)), \quad |s|>a.
$$
\end{prop}
\begin{proof} It is clear from the bounds \eqref{eq:Sbound1} and \eqref{eq:Sbound2} that $\cS^j_\theta$ maps $L^2(\Sigma_a)$ continuously to $L^2(-\infty,-a)\oplus L^2(a,+\infty)$. The analogous bounds \eqref{eq:Dbound1} and \eqref{eq:Dbound2} do not quite imply that $\cD^j_\theta\phi\in L^2(-\infty,-a)\oplus L^2(a,+\infty)$ for each $\phi\in L^2(\Sigma_a)$, since $(s-a)^{-1/2}\re^{-k\sin(\theta)(s-a)}$ is in $L^1(a,+\infty)$ but not in $L^2(a,+\infty)$. To see that $\cD^j_\theta:L^2(\Sigma_a)\to  L^2(a,+\infty)$ and is continuous we argue as in the proof of Proposition \ref{lem:SDbounds}, in particular using \eqref{eq:normderiv} and \eqref{eq:NormDivBound}, and recalling that, except when $\bsy^j\in \Sigma_a^{j+1}$, the right hand side of \eqref{eq:normderiv} is $\leq C\exp(-k\Im(z))$.  These bounds imply that,
for $\phi\in L^2(\Sigma_a)$ and $s>a$,
\begin{eqnarray*}
|\cD^j_\theta\phi(s)| &\leq &C\re^{-k\sin(\theta)(s-a)}\int_{\Sigma_a}|\phi(\bsy^j)|\, \rd s(\bsy^j) +  C\int_{\Sigma_a^{j+1}}\frac{|s-a|\, |\phi(\bsy^j)|}{(s-a)^2+(y^j_1-a)^2}\, \rd s(\bsy^j),
\end{eqnarray*}
where, throughout the proof, $C>0$ denotes a constant that depends only on $\theta$, $a$, and $k$. Thus
\begin{eqnarray*}
|\cD^j_\theta\phi(s)| &\leq &C\re^{-k\sin(\theta)(s-a)}\|\phi\|_{L^2(\Sigma_a)} + C D_0\psi(s),
\end{eqnarray*}
where $D_0$ is as in Remark \ref{rem:D0bounded} above, while $\psi\in L^2(-a,a)$ denotes the restriction of $|\phi|$ to $\Sigma^{j+1}_a$, precisely $\psi(y^j_1):= |\phi((y^j_1,a))|$, for $-a< y^j_1< a$, while $\psi(y^j_1):= 0$, for $|y^j_1|\geq a$. Since (Remark \ref{rem:D0bounded}) $D_0$ is a bounded operator on $L^2(\R)$, it follows that
%\scednote{30th 9th 2020. This doesn't now look right to me. I don't think that this is the DLP operator for $k=0$. Probably it should be, but there is an algebra error. }
%$$
%D_0\psi(s) := \int_{-a}^a  \frac{ (a-s)\,\psi(s)\,\rd  y^j_1}{(s-a)^2+(y^j_1-a)^2}, \quad s>a.
%$$
%But $D_0$ is nothing but (a multiple of) the static ($k=0$) version of the operator $D$ of Proposition \ref{prop:D_complexe}, considered as a mapping from $L^2(-a,a)$ to $L^2(a,\infty)$, and this operator is continuous, as noted in the proof of that proposition. Thus
$$
\|\cD^j_\theta\phi\|_{L^2(a,+\infty)} \leq C\|\phi\|_{L^2(\Sigma_a)} + C \|D_0\|\, \|\psi\|_{L^2(-a,a)} \leq C\|\phi\|_{L^2(\Sigma_a)}.
$$
Arguing in the same way, we see that  $\cD^j_\theta$ is also continuous as a mapping from  $L^2(\Sigma_a)$ to  $L^2(-\infty,-a)$, and the proof is complete.
\end{proof}

%\scednote{Just to flag - and I'm probably most at fault - that we have three notations for the boundary of $\Omega$ in the case that $\Omega = \R^2\setminus \overline{\Omega_a}$, namely $\partial \Omega$, $\partial \Omega_a$, $\Sigma_a$.}

%\bibliographystyle{siamplain}
%\bibliography{HsMM}

\end{document}